# DETERMINISTIC EQUIVALENTS FOR CERTAIN FUNCTIONALS OF LARGE RANDOM MATRICES[1]

By Walid Hachem, Philippe Loubaton and Jamal Najim

*Ecole Nationale Supérieure d'Electricité, Ecole Nationale Supérieure des Télécommunications, Université de Marne-La-Vallée and CNRS*

*Dedicated to Vyacheslav L. Girko on the occasion of his 60th birthday.*

Consider an $N \times n$ random matrix $Y_n = (Y_{ij}^n)$ where the entries are given by $Y_{ij}^n = \frac{\sigma_{ij}(n)}{\sqrt{n}} X_{ij}^n$, the $X_{ij}^n$ being independent and identically distributed, centered with unit variance and satisfying some mild moment assumption. Consider now a deterministic $N \times n$ matrix $A_n$ whose columns and rows are uniformly bounded in the Euclidean norm. Let $\Sigma_n = Y_n + A_n$. We prove in this article that there exists a deterministic $N \times N$ matrix-valued function $T_n(z)$ analytic in $\mathbb{C} - \mathbb{R}^+$ such that, almost surely,

$$\lim_{n \to +\infty, N/n \to c} \left( \frac{1}{N} \operatorname{Trace}(\Sigma_n \Sigma_n^T - z I_N)^{-1} - \frac{1}{N} \operatorname{Trace} T_n(z) \right) = 0.$$

Otherwise stated, there exists a deterministic equivalent to the empirical Stieltjes transform of the distribution of the eigenvalues of $\Sigma_n \Sigma_n^T$. For each $n$, the entries of matrix $T_n(z)$ are defined as the unique solutions of a certain system of nonlinear functional equations. It is also proved that $\frac{1}{N} \operatorname{Trace} T_n(z)$ is the Stieltjes transform of a probability measure $\pi_n(d\lambda)$, and that for every bounded continuous function $f$, the following convergence holds almost surely

$$\frac{1}{N} \sum_{k=1}^{N} f(\lambda_k) - \int_0^\infty f(\lambda) \pi_n(d\lambda) \underset{n \to \infty}{\longrightarrow} 0,$$

where the $(\lambda_k)_{1 \le k \le N}$ are the eigenvalues of $\Sigma_n \Sigma_n^T$. This work is motivated by the context of performance evaluation of multiple in-

Received July 2005; revised November 2006.

[1]Supported by the Fonds National de la Science (France) via the ACI program "Nouvelles Interfaces des Mathématiques," project MALCOM 205 and by NEWCOM European network.

puts/multiple output (MIMO) wireless digital communication channels. As an application, we derive a deterministic equivalent to the mutual information:

$$C_n(\sigma^2) = \frac{1}{N} \mathbb{E} \log \det \left( I_N + \frac{\Sigma_n \Sigma_n^T}{\sigma^2} \right),$$

where $\sigma^2$ is a known parameter.

## 1. Introduction.

*The model.* Consider an $N \times n$ random matrix $Y_n$ where the entries are given by

$$(1.1) \qquad\qquad Y_{ij}^n = \frac{\sigma_{ij}(n)}{\sqrt{n}} X_{ij}^n,$$

where $(\sigma_{ij}(n), 1 \le i \le N, 1 \le j \le n)$ is a bounded sequence of real numbers (i.e., $\sup_n \max_{(i,j)} |\sigma_{ij}(n)| = \sigma_{\max} < +\infty$) called a variance profile and the real random variables $X_{ij}^n$ are centered, independent and identically distributed (i.i.d.) with finite $4 + \varepsilon$ moment. Consider a real deterministic $N \times n$ matrix $A_n = (A_{ij}^n)$ whose columns $(\mathbf{a}_k^n)_{1 \le k \le n}$ and rows $(\tilde{\mathbf{a}}_\ell^n)_{1 \le \ell \le N}$ satisfy

$$(1.2) \qquad\qquad \sup_{n \ge 1} \max_{k,\ell} (\|\mathbf{a}_k^n\|, \|\tilde{\mathbf{a}}_\ell^n\|) < +\infty,$$

where $\| \cdot \|$ stands for the Euclidean norm. Denote by $\Sigma_n = Y_n + A_n$. This model has two interesting features: the random variables are independent but not i.i.d. since the variance may vary and $A_n$, the centering perturbation of $Y_n$, can have a very general form. Denote by $\mathrm{Tr}(C)$ the trace of matrix $C$. The purpose of this article is to study the behavior of

$$\frac{1}{N} \mathrm{Tr}(\Sigma_n \Sigma_n^T - zI_N)^{-1}, \qquad z \in \mathbb{C} - \mathbb{R}^+,$$

that is, the Stieltjes transform of the empirical eigenvalue distribution of $\Sigma_n \Sigma_n^T$ when $n \to +\infty$ and $N \to +\infty$ in such a way that $\frac{N}{n} \to c$, $0 < c < +\infty$. Under Condition (1.2), the convergence of the empirical distribution of the eigenvalues of $\Sigma_n \Sigma_n^T$ may fail to happen (see, e.g., Section 3.1). We shall prove that there exists a deterministic matrix-valued function $T_n(z)$, whose entries are defined by a system of equations, such that

$$\frac{1}{N} \mathrm{Tr}(\Sigma_n \Sigma_n^T - zI_N)^{-1} - \frac{1}{N} \mathrm{Tr}\, T_n(z) \underset{n \to \infty, N/n \to c}{\longrightarrow} 0, \qquad z \in \mathbb{C} - \mathbb{R}^+.$$

It is also proved that $\frac{1}{N} \mathrm{Tr}\, T_n(z)$ is the Stieltjes transform of a probability measure $\pi_n(d\lambda)$ and that the following convergence holds true:

$$(1.3) \qquad\qquad \frac{1}{N} \sum_{k=1}^{N} f(\lambda_k) - \int_0^\infty f(\lambda) \pi_n(d\lambda) \underset{n \to \infty}{\longrightarrow} 0,$$

where $f$ is a continuous and bounded function and where the $(\lambda_k)_{1 \le k \le N}$ are the eigenvalues of $\Sigma_n \Sigma_n^T$. The advantage of considering $\frac{1}{N} \operatorname{Tr} T_n(z)$ as a deterministic approximation instead of $\mathbb{E} \frac{1}{N} \operatorname{Tr}(\Sigma_n \Sigma_n^T - z I_N)^{-1}$ (which is deterministic as well) lies in the fact that $T_n(z)$ is in general far easier to compute than $\mathbb{E} \frac{1}{N} \operatorname{Tr}(\Sigma_n \Sigma_n^T - z I_N)^{-1}$ whose computation relies on Monte Carlo simulations. These Monte Carlo simulations become increasingly heavy as the size of matrix $\Sigma_n$ increases.

*Motivations.* This problem is mainly motivated by the context of performance analysis of multiple input/multiple output (MIMO) digital communication systems. The performance of these systems depends on the so-called channel matrix $H_n$ whose entries $(H_{i,j}^n, 1 \le i \le N, 1 \le j \le n)$ represent the gains between transmit antenna $j$ and receive antenna $i$. Matrix $H_n$ is often modeled as a realization of a random matrix. In certain contexts, the Gram matrix $H_n H_n^*$ is unitarily equivalent to a matrix $(A_n + Y_n)(A_n + Y_n)^*$ where $Y_n$ is a random matrix given by (1.1) and $A_n$ is a possibly full rank deterministic matrix satisfying (1.2). A fundamental result in information theory (see [26]) states that the mutual information $C_n$, which is a basic performance index of the channel, can be expressed in terms of the singular values of the channel matrix $H_n$:

$$C_n(\sigma^2) = \frac{1}{N} \mathbb{E} \left[ \log \det \left( I_N + \frac{H_n H_n^*}{\sigma^2} \right) \right],$$

where $\sigma^2$ represents the variance of an additive noise corrupting the received signals. We shall show in Section 4 how to approximate such a mutual information with the help of the deterministic equivalent $T_n(z)$. Related work can be found in [7, 22, 26].

*About the literature.* If $Z_n$ is a zero-mean $N \times n$ random matrix, the asymptotics of the spectrum of the $N \times N$ Gram random matrices $Z_n Z_n^T$ have been widely studied: Marčenko and Pastur [19], Yin [28], Silverstein et al. [5, 24, 25] for i.i.d. entries, Girko [9], Khorunzhy, Khorunzhenko and Pastur [17], Boutet de Monvel, Khorunzhy and Vasilchuk [4] (see also Shlyakhtenko [23]) for non-i.i.d. entries. In the centered case, it turns out that very often the empirical distribution of the eigenvalues converges toward a limiting distribution.

The case where matrix $Z_n$ has nonzero mean has been comparatively less studied. Among the related works, we first mention [5] which studies the eigenvalue asymptotics of the matrix $(R_n + Y_n)(R_n + Y_n)^T$ in the case where the matrices $Y_n$ and $R_n$ are independent random matrices, $Y_n$ has i.i.d. entries and the empirical distribution of $R_n R_n^T$ converges to a nonrandom distribution. It is shown there that the eigenvalue distribution converges almost surely toward a deterministic distribution whose Stieltjes transform is uniquely defined by a certain functional equation. The case $(Y_n + \Delta_n)(Y_n + \Delta_n)^T$ where $Y_n$ is given by (1.1) and $\Delta_n$ is deterministic pseudo-diagonal has been studied in [13] where it is shown that under suitable assumptions the eigenvalue distribution converges. In the general case $\Sigma_n = Y_n + A_n$, the convergence of the empirical distribution of the eigenvalues of $\Sigma_n \Sigma_n^T$ may fail to happen even if the variance profile exists in some limit and the spectral distribution of $A_n A_n^T$ converges (see, e.g., Section 3.1).

Girko proposed in [10], Chapter 7, to study a deterministic equivalent of some functionals of the eigenvalues of $\Sigma_n \Sigma_n^T$ in the case where the following condition holds for $A_n$:

$$(1.4) \qquad \sup_n \max_i \sum_{j=1}^n |A_{i,j}^n| < +\infty \quad \text{and} \quad \sup_N \max_j \sum_{i=1}^N |A_{i,j}^n| < +\infty.$$

Girko showed that the entries of the resolvent $(\Sigma_n \Sigma_n^T - zI)^{-1}$ have the same asymptotic behavior as the entries of a certain deterministic holomorphic $N \times N$ matrix-valued function $T_n(z)$ characterized by a nonlinear system of $(n+N)$ coupled functional equations. Condition (1.4) is however rather restrictive. In particular, it does not hold in the context of wireless MIMO system in which (1.2) is actually much more relevant. In this paper, we thus extend some of the results of Girko [10] to the case where $A_n$ satisfies (1.2). For this, we do not follow the approach of Girko [10] based on the use of Cramér's rule but rather take the approach of Dozier and Silverstein [5] as a starting point. This approach not only allows to extend the result of [10] to deterministic matrices satisfying (1.2), but also provides a simpler proof. It also enables us to prove that the deterministic equivalent $\frac{1}{N} \operatorname{Tr} T_n(z)$ of the Stieltjes transform $\frac{1}{N} \operatorname{Tr}(\Sigma_n \Sigma_n^T - zI)^{-1}$ is itself the Stieltjes transform of a probability measure, which is a result of interest from a practical point of view [see (1.3), e.g.].

*Outline of the paper.* The main results of the paper are Theorem 2.4 [existence of the deterministic matrix-valued function $T_n(z)$], Theorem 2.5 (deterministic approximation of the Stieltjes transform) and Theorem 4.1 (deterministic approximation of the mutual information of a wireless channel). Theorems 2.4 and 2.5 are stated in Section 2. Section 3 is devoted to some examples and remarks while Section 4 is devoted to applications to digital communication; in particular, Theorem 4.1 is stated in Section 4. Proof of Theorem 2.4 is established in Section 5 while Theorem 2.5 is proved to the main extent in Section 6. Technical results and some of the proofs are collected in the Appendix.

## 2. The deterministic equivalent: main results.

2.1. *Notation, assumptions and preliminary results.* Let $N = N(n)$ be a sequence of integers such that $\lim_{n \to \infty} \frac{N}{n} = c \in (0, \infty)$. We denote by $\mathbf{i}$ the complex number $\sqrt{-1}$ and by $\mathrm{Im}(z)$ the imaginary part of $z \in \mathbb{C}$. Consider an $N \times n$ random matrix $Y_n$ where the entries are given by

$$Y_{ij}^n = \frac{\sigma_{ij}(n)}{\sqrt{n}} X_{ij}^n,$$

where $X_{ij}^n$ and $(\sigma_{ij}(n))$ are defined below.

ASSUMPTION A-1. The random variables $(X_{ij}^n; 1 \le i \le N, 1 \le j \le n, n \ge 1)$ are real, independent and identically distributed. They are centered with $\mathbb{E}(X_{ij}^n)^2 = 1$. Moreover, there exists $\varepsilon > 0$ such that

$$\mathbb{E}|X_{ij}^n|^{4+\varepsilon} < \infty,$$

where $\mathbb{E}$ denotes expectation.

REMARK 2.1. We can (and will) assume in the proofs that $\varepsilon < 4$ without lack of generality.

ASSUMPTION A-2. There exists a nonnegative finite real number $\sigma_{\max}$ such that the family of real numbers $(\sigma_{ij}(n), 1 \le i \le N, 1 \le j \le n, n \ge 1)$ satisfies

$$\sup_{n \ge 1} \max_{\substack{1 \le i \le N \\ 1 \le j \le n}} |\sigma_{ij}(n)| \le \sigma_{\max}.$$

Denote by $\mathrm{var}(Z)$ the variance of the random variable $Z$. Since $\mathrm{var}(Y_{ij}^n) = \sigma_{ij}^2(n)/n$, the family $(\sigma_{ij}(n))$ will be called a variance profile.

Let $A_n = (A_{ij}^n)$ be an $N \times n$ real deterministic matrix. We introduce the $N \times n$ matrix

$$\Sigma_n = Y_n + A_n.$$

For every matrix $M$, we will denote by $M^T$ (resp. $M^*$) its transpose (resp. its Hermitian adjoint), by $\mathrm{Tr}(M)$ [resp. $\det(M)$] its trace (resp. its determinant if $M$ is square) and by $F^{MM^T}$, the empirical distribution function of the eigenvalues of $MM^T$. The matrix $I_n$ will always refer to the $n \times n$ identity.

Let $A$ be an $n \times n$ matrix with complex entries. If $A$ is Hermitian and positive semi-definite (i.e., $z^* A z \ge 0$ for every $z \in \mathbb{C}^n$), we write $A \ge 0$. If $A$ and $B$ are Hermitian matrices, $A \ge B$ means $A - B \ge 0$.

We define the matrix $\mathrm{Im}(A)$ by

$$\mathrm{Im}(A) = \frac{1}{2\mathbf{i}}(A - A^*).$$

Let $\mathcal{B}$ be the Borel $\sigma$-field on $\mathbb{R}$. An $n \times n$ matrix-valued measure $M$ on $\mathcal{B}$ is a matrix-valued function on $\mathcal{B}$, the entries of which are complex measures. A positive matrix-valued measure $M$ is such that for every $\Delta \in \mathcal{B}$, $M(\Delta) \geq 0$ [i.e., $M(\Delta)$ is a Hermitian and positive semi-definite matrix]. In this case, the diagonal entries are nonnegative measures and for every $z \in \mathbb{C}^n$, $z^*Mz$ is a scalar nonnegative measure. For more details on matrix-valued measures, the reader is referred to the standard textbook [20].

In this article, we shall deal with several norms. Denote by $\|\cdot\|$ the Euclidean norm for vectors. In the case of matrices, the norm $\|\cdot\|_{\mathrm{sp}}$ will refer to the spectral norm and if $Z$ is a complex-valued random variable, denote by $\|Z\|_p = (\mathbb{E}|Z|^p)^{1/p}$ for $p > 0$. If $\mathcal{X}$ is a topological space, we endow it with its Borel $\sigma$-algebra $\mathcal{B}(\mathcal{X})$ and we denote by $\mathcal{P}(\mathcal{X})$ the set of probability measures over $\mathcal{X}$.

We denote by $\mathbb{C}^- = \{z \in \mathbb{C}, \mathrm{Im}(z) < 0\}$ and by $\mathbb{C}^+ = \{z \in \mathbb{C}, \mathrm{Im}(z) > 0\}$.

ASSUMPTION A-3.   Denote by $\mathbf{a}_k^n$ the $k$th column of $A_n$, by $\tilde{\mathbf{a}}_\ell^n$ its $\ell$th row and by

$$\mathbf{a}_{\max,n} = \max(\|\mathbf{a}_k^n\|, \|\tilde{\mathbf{a}}_\ell^n\|; 1 \leq k \leq n, 1 \leq \ell \leq N).$$

We assume that

$$\mathbf{a}_{\max} = \sup_{n \geq 1} \mathbf{a}_{\max,n} < \infty.$$

We will frequently drop $n$ and simply write $\sigma_{ij}$, $A$, $Y$, $\mathbf{a}_k$, and so on.

If $(\alpha_1, \ldots, \alpha_k)$ is a finite sequence of real numbers, we denote by $\mathrm{diag}(\alpha_1, \ldots, \alpha_k)$ the $k \times k$ diagonal matrix whose diagonal elements are the $\alpha_i$'s. Let

$$(2.1) \quad D_j = \mathrm{diag}(\sigma_{ij}^2(n); 1 \leq i \leq N) \quad \text{and} \quad \tilde{D}_i = \mathrm{diag}(\sigma_{ij}^2(n); 1 \leq j \leq n).$$

We will denote by $D_j^{1/2} = \mathrm{diag}(\sigma_{ij}, i \leq N)$ and $\tilde{D}_i^{1/2} = \mathrm{diag}(\sigma_{ij}, j \leq n)$. Let $\mu$ be a probability measure over $\mathbb{R}$. Its Stieltjes transform $f$ is defined by

$$f(z) = \int_{\mathbb{R}} \frac{\mu(d\lambda)}{\lambda - z}, \qquad z \in \mathbb{C}^+.$$

We list below the main properties of the Stieltjes transforms that will be needed in the sequel.

PROPOSITION 2.2.   *The following properties hold true:*
   *1. Let $f$ be the Stieltjes transform of $\mu$, then*

– the function $f$ is analytic over $\mathbb{C}^+$,
– if $z \in \mathbb{C}^+$ then $f(z) \in \mathbb{C}^+$,
– the function $f$ satisfies: $|f(z)| \leq \frac{1}{\mathrm{Im}(z)}$ and $\mathrm{Im}(\frac{1}{f(z)}) \leq -\mathrm{Im}(z)$ over $\mathbb{C}^+$,
– if $\mu(-\infty, 0) = 0$ then its Stieltjes transform $f$ is analytic over $\mathbb{C} - \mathbb{R}^+$.

Moreover, $z \in \mathbb{C}^+$ implies $zf(z) \in \mathbb{C}^+$, and the following controls hold:

$$(2.2) \qquad |f(z)| \leq \begin{cases} |\mathrm{Im}(z)|^{-1}, & \text{if } z \in \mathbb{C} - \mathbb{R}, \\ |z|^{-1}, & \text{if } z \in (-\infty, 0), \\ (\mathrm{dist}(z, \mathbb{R}^+))^{-1}, & \text{if } z \in \mathbb{C} - \mathbb{R}^+, \end{cases}$$

where $\mathrm{dist}$ stands for the Euclidean distance.

2. Conversely, let $f$ be a function analytic over $\mathbb{C}^+$ such that $f(z) \in \mathbb{C}^+$ if $z \in \mathbb{C}^+$. If $\lim_{y \to +\infty} -\mathbf{i}y f(\mathbf{i}y) = 1$, then $f$ is the Stieltjes transform of a probability measure $\mu$ and the following inversion formula holds:

$$\mu([a, b]) = \lim_{\eta \to 0^+} \frac{1}{\pi} \int_a^b \mathrm{Im}\, f(\xi + \mathbf{i}\eta)\, d\xi,$$

whenever $a$ and $b$ are continuity points of $\mu$. If moreover $zf(z) \in \mathbb{C}^+$ for $z$ in $\mathbb{C}^+$ then, $\mu(\mathbb{R}^-) = 0$. In particular, $f$ is given by

$$f(z) = \int_0^\infty \frac{\mu(d\lambda)}{\lambda - z}.$$

and has an analytic continuation on $\mathbb{C} - \mathbb{R}^+$.

3. Let $F$ be an $n \times n$ matrix-valued function over $\mathbb{C}^+$ such that

– the function $F$ is analytic,
– the matrices $\mathrm{Im}\, F(z)$ and $\mathrm{Im}\, zF(z)$ satisfy

$$\mathrm{Im}\, F(z) \geq 0 \quad \text{and} \quad \mathrm{Im}\, zF(z) \geq 0 \qquad \text{for } z \in \mathbb{C}^+.$$

Then there exists an $n \times n$ matrix $C \geq 0$ and a positive $n \times n$ matrix-valued measure $\mu$ over $\mathbb{R}^+$ such that

$$F(z) = C + \int_{\mathbb{R}^+} \frac{\mu(d\lambda)}{\lambda - z} \qquad \text{with } \mathrm{Tr} \int_{\mathbb{R}^+} \frac{\mu(d\lambda)}{1 + \lambda} < \infty.$$

4. Let $\mathbb{P}_n$ and $\mathbb{P}$ be probability measures over $\mathbb{R}$ and denote by $f_n$ and $f$ their Stieltjes transforms. Then

$$\left( \forall z \in \mathbb{C}^+, f_n(z) \xrightarrow[n \to \infty]{} f(z) \right) \quad \Longleftrightarrow \quad \mathbb{P}_n \xrightarrow[n \to \infty]{\mathcal{D}} \mathbb{P},$$

where $\mathcal{D}$ stands for convergence in distribution.

As item 3 is perhaps less standard than the other items, we give a proof in Appendix A. These results can be found in several papers without proof

(see, e.g., [3], pages 64–65).

In the sequel we shall denote by $\mathcal{S}(\mathbb{R}^+)$ the set of Stieltjes transforms of probability measures over $\mathbb{R}^+$. Proposition 2.2 part 2 gives necessary and sufficient conditions for $f$ to be in $\mathcal{S}(\mathbb{R}^+)$.

We list below some useful properties of matrices that will be of constant use later:

PROPOSITION 2.3.    *1. Let $B$ be an $N \times n$ matrix and let $C$ be an $n \times k$ matrix, then*

$$\|BC\|_{\mathrm{sp}} \leq \|B\|_{\mathrm{sp}} \|C\|_{\mathrm{sp}}.$$

*2. Let $B$ be an $N \times n$ matrix, let $\mathbf{y}$ and $\mathbf{x}$ be respectively $N \times 1$ and $n \times 1$ vectors, then*

$$\|B\|_{\mathrm{sp}} = \sup\{|\mathbf{y}^T B \mathbf{x}|, \|\mathbf{y}\| = \|\mathbf{x}\| = 1\}.$$

*In particular,*

$$|B_{ij}| \leq \|B\|_{\mathrm{sp}}.$$

*3. If $B$ and $C$ are square and invertible matrices, then*

$$B^{-1} - C^{-1} = -B^{-1}(B - C)C^{-1}.$$

*If $B$ and $C$ are resolvents (see below), then this identity is known as the resolvent identity.*

*4. Let $X$ be an $N \times n$ matrix and denote by $Q(z) = (Q_{ij}(z))$ the resolvent of $XX^T$, that is, $Q(z) = (XX^T - zI)^{-1}$.*

(a)  *We have $\|Q(z)\|_{\mathrm{sp}} \leq |\mathrm{Im}(z)|^{-1}$ and in particular $|Q_{ij}(z)| \leq |\mathrm{Im}(z)|^{-1}$.*

(b)  *If $D$ is an $N \times N$ diagonal matrix with nonnegative diagonal elements, then*

$$b(z) = -\frac{1}{z(1 + (1/n)\,\mathrm{Tr}\,DQ(z))} \in \mathcal{S}(\mathbb{R}^+).$$

Although very simple, the identity in Proposition 2.3 part 3 is extremely useful and will be of constant use in the sequel.

2.2.  *The deterministic equivalent: existence and asymptotic behavior.*  Let us introduce the following resolvents:

$$Q_n(z) = (\Sigma_n \Sigma_n^T - zI_N)^{-1} = (q_{ij}(z))_{1 \leq i,j \leq N}, \qquad z \in \mathbb{C} - \mathbb{R}^+,$$

$$\tilde{Q}_n(z) = (\Sigma_n^T \Sigma_n - zI_n)^{-1} = (\tilde{q}_{ij}(z))_{1 \leq i,j \leq n}, \qquad z \in \mathbb{C} - \mathbb{R}^+.$$

The function $f_n(z) = \frac{1}{N}\,\mathrm{Tr}\,Q_n(z)$ is the Stieltjes transform of the empirical distribution of the eigenvalues of $\Sigma_n \Sigma_n^T$. The aim of this paper is to get some insight on $f_n$ by the means of a deterministic equivalent. We will often drop subscripts $n$. In Theorem 2.4, we prove the existence of matrix-valued deterministic functions $T(z)$ and $\tilde{T}(z)$ which satisfy a system of $N+n$ functional equations. In Theorem 2.5, we prove that asymptotically,

$$\frac{1}{N}\operatorname{Tr}Q(z) \sim \frac{1}{N}\operatorname{Tr}T(z) \quad \text{and} \quad \frac{1}{n}\operatorname{Tr}\tilde{Q}(z) \sim \frac{1}{n}\operatorname{Tr}\tilde{T}(z)$$

in a sense to be defined.

THEOREM 2.4. *Let $A_n$ be an $N \times n$ deterministic matrix. The deterministic system of $N+n$ equations,*

$$(2.3) \qquad \psi_i(z) = \frac{-1}{z(1+(1/n)\operatorname{Tr}\tilde{D}_i\tilde{T}(z))} \qquad \text{for } 1 \le i \le N,$$

$$(2.4) \qquad \tilde{\psi}_j(z) = \frac{-1}{z(1+(1/n)\operatorname{Tr}D_jT(z))} \qquad \text{for } 1 \le j \le n,$$

*where*

$$(2.5) \qquad \begin{aligned} &\Psi(z) = \operatorname{diag}(\psi_i(z), 1 \le i \le N), \\ &\tilde{\Psi}(z) = \operatorname{diag}(\tilde{\psi}_j(z), 1 \le j \le n), \\ &T(z) = (\Psi^{-1}(z) - zA\tilde{\Psi}(z)A^T)^{-1}, \\ &\tilde{T}(z) = (\tilde{\Psi}^{-1}(z) - zA^T\Psi(z)A)^{-1}, \end{aligned}$$

*admits a unique solution $(\psi_1, \ldots, \psi_N, \tilde{\psi}_1, \ldots, \tilde{\psi}_n)$ in $\mathcal{S}(\mathbb{R}^+)^{N+n}$.*

*Moreover, there exist a positive $N \times N$ matrix-valued measure $\mu = (\mu_{ij})$ and a positive $n \times n$ matrix-valued measure $\tilde{\mu} = (\tilde{\mu}_{ij})$ such that*

$$\mu(\mathbb{R}^+) = I_N, \qquad \tilde{\mu}(\mathbb{R}^+) = I_n \quad \text{and}$$

$$T(z) = \int_{\mathbb{R}^+}\frac{\mu(d\lambda)}{\lambda - z}, \qquad \tilde{T}(z) = \int_{\mathbb{R}^+}\frac{\tilde{\mu}(d\lambda)}{\lambda - z}$$

*for $z \in \mathbb{C} - \mathbb{R}^+$. In particular, $\frac{1}{N}\operatorname{Tr}T(z)$ and $\frac{1}{n}\operatorname{Tr}\tilde{T}(z)$ are Stieltjes transforms of probability measures.*

Proof of Theorem 2.4 is postponed to Section 5. It relies on an iteration scheme and on properties of matrix-valued Stieltjes transforms.

THEOREM 2.5. *Assume that Assumptions A-1, A-2 and A-3 hold, then the following limits hold true almost everywhere:*

$$\lim_{n\to\infty, N/n\to c}\left(\frac{1}{N}\operatorname{Tr}Q(z) - \frac{1}{N}\operatorname{Tr}T(z)\right) = 0 \qquad \forall z \in \mathbb{C} - \mathbb{R}^+,$$

$$\lim_{n\to\infty, N/n\to c}\left(\frac{1}{n}\operatorname{Tr}\tilde{Q}(z) - \frac{1}{n}\operatorname{Tr}\tilde{T}(z)\right) = 0 \qquad \forall z \in \mathbb{C} - \mathbb{R}^+,$$

REMARK 2.6 (Limiting behavior). There are two well-known cases where the empirical distribution of the eigenvalues of $\Sigma_n \Sigma_n^T$ converges toward a limit expressed in terms of its Stieltjes transform: The case where the variance profile $\sigma_{ij} = \sigma$ is a constant [5] and the case where the centering matrix $A$ (which can be rectangular) has elements equal to zero outside the diagonal [13] (with $F^{A_n A_n^T}$ converging to a probability measure in both cases). Interestingly, one can obtain discretized versions of the limiting equations in [5] and [13] by combining (2.3)–(2.5).

COROLLARY 2.7. *Assume that Assumptions* A-1, A-2 *and* A-3 *hold. Denote by* $\mathbb{P}_n$ *and* $\pi_n$ *the probability measures whose Stieltjes transform are respectively* $\frac{1}{N}\operatorname{Tr}Q_n$ *and* $\frac{1}{N}\operatorname{Tr}T_n$. *Then the following limit holds true almost everywhere:*

$$\int_0^\infty f(\lambda)\mathbb{P}_n(d\lambda) - \int_0^\infty f(\lambda)\pi_n(d\lambda) \xrightarrow[n\to\infty]{} 0,$$

*where* $f:\mathbb{R}^+ \to \mathbb{R}$ *is a continuous and bounded function. The same results hold for the probability measures related to* $\frac{1}{n}\operatorname{Tr}\tilde{Q}_n$ *and* $\frac{1}{n}\operatorname{Tr}\tilde{T}_n$.

PROOF OF COROLLARY 2.7. We rely on the fact that $(\mathbb{P}_n)$ is almost surely tight and that $(\pi_n)$ is tight. In order to prove this, it is sufficient to prove

$$\sup_n \left( \int \lambda\mathbb{P}_n(d\lambda) \right) < \infty \qquad \text{a.s.} \quad \text{and} \quad \sup_n \left( \int \lambda\pi_n(d\lambda) \right) < \infty.$$

One has

$$\begin{aligned}
\int \lambda\mathbb{P}_n(d\lambda) &= \frac{1}{N}\operatorname{Tr}\Sigma_n\Sigma_n^T \\
&\leq \frac{2}{N}\operatorname{Tr}A_nA_n^T + \frac{2}{N}\operatorname{Tr}Y_nY_n^T \\
&\overset{(a)}{\leq} 2\mathbf{a}_{\max}^2 + \frac{2\sigma_{\max}^2}{Nn}\sum_{i,j}X_{ij}^2,
\end{aligned}$$

where (a) follows from Assumption A-3. The almost sure tightness of $(\mathbb{P}_n)$ follows then from the fact that $\frac{1}{Nn}\sum_{i,j}X_{ij}^2$ converges almost surely to one. The fact that $\sup_n(\int \lambda\pi_n(d\lambda)) < \infty$ is established in Lemma C.1.

Let us now consider a realization for which the conclusion of Theorem 2.5 holds, and $(\mathbb{P}_n)$ is tight. Then for any subsequence $n'$ of the natural numbers, there exists a further subsequence $n''$ for which both $\mathbb{P}_{n''}$ and $\pi_{n''}$ converge in distribution to, say $\mathbb{P}^*$ and $\pi^*$. Then their Stieltjes transforms also converge and, from Theorem 2.5, to each other. Therefore $\mathbb{P}^* = \pi^*$ and the conclusion holds on $(n'')$. Since this is true for any subsequence $(n')$, the conclusion of Corollary 2.7 holds for all $n$, for this realization, which occurs with probability one. □

We are indebted to the referee for the proof of Corollary 2.7 which is much simpler than the original one.

REMARK 2.8 (Concentration and martingale arguments). If one is interested in proving that the empirical measure is close to its expectation, one can rely on concentration arguments (see, e.g., [11]) at least when the entries are Gaussian, bounded or satisfy the Poincaré inequality. One can also rely on martingale arguments, regardless the nature of the entries (as long as they are independent and satisfy some mild moment assumptions, see [10], Chapter 16, and also [6]). The purpose here is to provide a "computable" deterministic equivalent which is not expressed in terms of expectations. In fact, although expectations can be computed by Monte Carlo methods, these methods quickly imply a huge amount of computations when the size of the matrix models increases.

2.3. *Outline of the proof of Theorem* 2.5. The proof relies on the introduction of intermediate quantities which are the stochastic counterparts of $\Psi$, $\tilde{\Psi}$, $T$ and $\tilde{T}$. Denote by

$$
(2.6) \quad
\begin{aligned}
&b_i(z) = \frac{-1}{z(1 + (1/n)\operatorname{Tr}\tilde{D}_i\tilde{Q}(z))} \qquad \text{for } 1 \leq i \leq N, \\
&B(z) = \operatorname{diag}(b_i(z), 1 \leq i \leq N),
\end{aligned}
$$

$$
(2.7) \quad
\begin{aligned}
&\tilde{b}_j(z) = \frac{-1}{z(1 + (1/n)\operatorname{Tr} D_j Q(z))} \qquad \text{for } 1 \leq j \leq n, \\
&\tilde{B}(z) = \operatorname{diag}(\tilde{b}_j(z), 1 \leq j \leq n),
\end{aligned}
$$

and by

$$
(2.8) \quad
\begin{aligned}
&R(z) = (B^{-1}(z) - zA\tilde{B}(z)A^T)^{-1} \quad \text{and} \\
&\tilde{R}(z) = (\tilde{B}^{-1}(z) - zA^T B(z)A)^{-1}.
\end{aligned}
$$

The introduction of the quantities $R, \tilde{R}, T$ and $\tilde{T}$ can be traced back to the work of Girko. We first prove that $\frac{1}{N}\operatorname{Tr} Q(z) \sim \frac{1}{N}\operatorname{Tr} R(z)$ and $\frac{1}{N}\operatorname{Tr}\tilde{Q}(z) \sim \frac{1}{N}\operatorname{Tr}\tilde{R}(z)$ in Lemma 6.1. These computations, quite involved, are along the same lines as the computations by Dozier and Silverstein in [5]. We then prove that $\frac{1}{N}\operatorname{Tr} R(z) \sim \frac{1}{N}\operatorname{Tr} T(z)$ and $\frac{1}{N}\operatorname{Tr}\tilde{R}(z) \sim \frac{1}{N}\operatorname{Tr}\tilde{T}(z)$ in Lemma 6.6.

**3. Examples and remarks.**

3.1. *An example where the spectral measure of $\Sigma_n \Sigma_n^T$ does not converge.* Let $Y_n$ be a $2n \times 2n$ matrix such that

$$Y_n = \begin{pmatrix} W_n & 0 \\ 0 & 0 \end{pmatrix},$$

where $W_n$ is $n \times n$ and $W_{ij}^n = \frac{X_{ij}}{\sqrt{n}}$, the $X_{ij}$'s being i.i.d. Matrix $Y_n$ can be interpreted as a matrix with a variance profile $\sigma^2(x,y)$ in the sense that

$$
\begin{aligned}
Y_{ij}^n &= \sigma\left(\frac{i}{2n}, \frac{j}{2n}\right) \frac{X_{ij}}{\sqrt{n}} \\
&\text{where } \sigma(x,y) = \begin{cases} 1, & \text{if } (x,y) \in [0,1/2] \times [0,1/2], \\ 0, & \text{else.} \end{cases}
\end{aligned}
$$

(3.1)

Consider now the following $2n \times 2n$ deterministic matrices $B_n$ and $\tilde{B}_n$ defined as

$$B_n = \begin{pmatrix} I_n & 0 \\ 0 & 0 \end{pmatrix} \quad \text{and} \quad \tilde{B}_n = \begin{pmatrix} 0 & 0 \\ 0 & I_n \end{pmatrix}.$$

Then both $F^{B_n B_n^T}$ and $F^{\tilde{B}_n \tilde{B}_n^T}$ converge toward $\frac{1}{2}\delta_0(dx) + \frac{1}{2}\delta_1(dx)$. It is well known that the limiting distribution of $F^{W_n W_n^T}$ and $F^{(W_n+I_n)(W_n+I_n)^T}$ exist and are nonrandom probability distributions. The first limiting probability distribution is known as Marčenko–Pastur's distribution, and we denote it by $\mathbb{P}_{\mathrm{MP}}$. Denote by $\mathbb{P}_{\mathrm{cub}}$ the limiting distribution (This probability distribution is sometimes referred to as the "cubic law" since the limiting Stieltjes transform of the normalized trace of the eigenvalues of $(W_n+I_n)(W_n+I_n)^T$ satisfies a cubic equation (see [10], Chapter 2 for details), hence the notation $\mathbb{P}_{\mathrm{cub}}$.) of $F^{(W_n+I_n)(W_n+I_n)^T}$ (for details, see [9, 13, 19, 25]). The following holds true:

PROPOSITION 3.1. *Let $\Upsilon_n = (Y_n + B_n)$ and let $\tilde{\Upsilon}_n = (Y_n + \tilde{B}_n)$. Then*

$$F^{\Upsilon_n \Upsilon_n^T} \xrightarrow[n \to \infty]{} \tfrac{1}{2}\mathbb{P}_{\mathrm{cub}} + \tfrac{1}{2}\delta_0 \qquad a.s.$$

$$F^{\tilde{\Upsilon}_n \tilde{\Upsilon}_n^T} \xrightarrow[n \to \infty]{} \tfrac{1}{2}\mathbb{P}_{\mathrm{MP}} + \tfrac{1}{2}\delta_1 \qquad a.s.$$

*Consider in particular the $2n \times 2n$ random matrix defined by*

$$\Sigma_n = \begin{cases} \Upsilon_n, & \text{if } n = 2p, \\ \tilde{\Upsilon}_n, & \text{if } n = 2p+1, \end{cases}$$

*that is, $\Sigma_n = (Y_n + A_n)$ where $A_n$ is alternatively equal to $B_n$ and to $\tilde{B}_n$, then $F^{\Sigma_n \Sigma_n^T}$ does not admit a limiting distribution.*

Proof. Let $(W_n + I_n)(W_n + I_n)^T = U_n \Delta_n U_n^T$ where $\Delta_n = \text{diag}(\lambda_i, 1 \leq i \leq n)$ and $U_n$ is orthogonal. Then $\frac{1}{n} \sum \delta_{\lambda_i} \to \mathbb{P}_{\text{cub}}$. Now

$$\Upsilon_n \Upsilon_n^T = \begin{pmatrix} U_n & 0 \\ 0 & I_n \end{pmatrix} \begin{pmatrix} \Delta_n & 0 \\ 0 & 0 \end{pmatrix} \begin{pmatrix} U_n^T & 0 \\ 0 & I_n \end{pmatrix},$$

and $F^{\Upsilon_n \Upsilon_n^T} = \frac{1}{2n}(\sum_{i=1}^n \delta_{\lambda_i} + n\delta_0)$ which converges toward $\frac{1}{2}\mathbb{P}_{\text{cub}} + \frac{1}{2}\delta_0$. One can prove similarly the second statement of the proposition. The absence of convergence for $F^{\Sigma_n \Sigma_n^T}$ follows. $\square$

REMARK 3.2. In view of Proposition 3.1, $F^{\Sigma_n \Sigma_n^T}$ does not have a limit since it oscillates between two different limiting distributions, despite the fact that:

(i) the variance profile is the sampling of a given function,
(ii) $F^{AA^T}$ converges toward $\frac{1}{2}\delta_0 + \frac{1}{2}\delta_1$.

The general model $\Sigma_n = (Y_n + A_n)$ where $A_n$ is diagonal has been studied in [13]. In particular, an assumption is required in order to insure the convergence of $F^{\Sigma_n \Sigma_n^T}$. In the setting of the present example [where $Y_n$ is defined by (3.1) and $A_n$ is alternatively equal to $B_n$ or $\tilde{B}_n$], this assumption writes: Denote by $(a_{ii}, 1 \leq i \leq 2n)$ the diagonal elements of $A_n$, then if the empirical (deterministic) distribution

$$H_n = \frac{1}{2n} \sum_{i=1}^{2n} \delta_{(i/2n, a_{ii}^2)}$$

converges toward a probability distribution with a compact support, then $F^{\Sigma_n \Sigma_n^T}$ converges toward a deterministic distribution which is properly described via its Stieltjes transform.

This assumption expresses a "coupled convergence" between the variance profile and the diagonal elements of $A_n = \text{diag}(a_{ii})$ and one can easily check that it is not fulfilled here. In fact, denote by $(b_{ii}, 1 \leq i \leq 2n)$ the diagonal elements of matrix $B_n$ and by $(\tilde{b}_{ii}, 1 \leq i \leq 2n)$ those of matrix $\tilde{B}_n$. Then if $n = 2p$,

$$H_n = \frac{1}{2n} \sum_{i=1}^{2n} \delta_{(i/2n, b_{ii}^2)} \xrightarrow[n \to \infty]{} 1_{[0,1/2]}(x)\,dx \otimes \delta_1(dy) + 1_{[1/2,1]}(x)\,dx \otimes \delta_0(dy),$$

while if $n = 2p + 1$, we have

$$H_n = \frac{1}{2n} \sum_{i=1}^{2n} \delta_{(i/2n, \tilde{b}_{ii}^2)} \xrightarrow[n \to \infty]{} 1_{[0,1/2]}(x)\,dx \otimes \delta_0(dy) + 1_{[1/2,1]}(x)\,dx \otimes \delta_1(dy).$$

This, in particular, does not allow the convergence of $H_n$. We believe that this is the main reason for which $F^{\Sigma \Sigma^T}$ does not converge in the present example.

3.2. *The case of a separable variance profile.* In this section, we consider the case where the variance profile is separable, that is,

$$\sigma_{ij}^2 = d_i \tilde{d}_j, \qquad 1 \le i \le N, 1 \le j \le n.$$

Note that this occurs quite frequently in the context of MIMO wireless systems. The system of $N + n$ equations defining $\Psi(z)$ and $\tilde{\Psi}(z)$ takes a much simpler form: In fact it reduces to a pair of equations. To see this, we denote by $D$ and $\tilde{D}$ the diagonal matrices $D = \mathrm{diag}(d_1, \ldots, d_N)$ and $\tilde{D} = \mathrm{diag}(\tilde{d}_1, \ldots, \tilde{d}_n)$. Note that $D_j = \tilde{d}_j D$ and $\tilde{D}_i = d_i \tilde{D}$. Denote by $\delta(z)$ and $\tilde{\delta}(z)$ the functions defined on $\mathbb{C} - \mathbb{R}^+$ by

$$(3.2) \qquad \begin{cases} \delta(z) = \dfrac{1}{n} \mathrm{Tr}(D T(z)), \\ \tilde{\delta}(z) = \dfrac{1}{n} \mathrm{Tr}(\tilde{D} \tilde{T}(z)). \end{cases}$$

It is obvious that $\Psi(z) = -\frac{1}{z}(I + \tilde{\delta} D)^{-1}$ and $\tilde{\Psi}(z) = -\frac{1}{z}(I + \delta \tilde{D})^{-1}$. Therefore, the solutions of the system of $N + n$ equations (2.5) can be expressed in terms of functions $\delta(z)$ and $\tilde{\delta}(z)$. Using the fact that $T(z) = (\Psi(z)^{-1} - zA\tilde{\Psi}(z)A^T)^{-1}$ and $\tilde{T}(z) = (\tilde{\Psi}(z)^{-1} - zA^T\Psi(z)A)^{-1}$ and the definition (3.2), we get that $\delta(z)$ and $\tilde{\delta}(z)$ are solutions of the following system of two equations:

$$(3.3) \qquad \begin{cases} \delta(z) = \dfrac{1}{n} \mathrm{Tr}[D(-z(I + D\tilde{\delta}) + A(I + \tilde{D}\delta)^{-1}A^T)^{-1}], \\ \tilde{\delta}(z) = \dfrac{1}{n} \mathrm{Tr}[\tilde{D}(-z(I + \tilde{D}\delta) + A^T(I + D\tilde{\delta})^{-1}A)^{-1}]. \end{cases}$$

Of course, Theorem 2.4 implies that this system has a unique solution in the class of functions $(\delta, \tilde{\delta})$ where $\delta$ and $\tilde{\delta}$ are Stieltjes transforms of positive measures $(\mu, \tilde{\mu})$ on $\mathbb{R}^+$ such that $\mu(\mathbb{R}^+) = \frac{1}{n}\mathrm{Tr}(D)$ and $\tilde{\mu}(\mathbb{R}^+) = \frac{1}{n}\mathrm{Tr}(\tilde{D})$. Moreover at a given point $z = -\sigma^2$, it can be shown that the following system:

$$(3.4) \qquad \begin{cases} \kappa = \dfrac{1}{n} \mathrm{Tr}[D(\sigma^2(I + D\tilde{\kappa}) + A(I + \tilde{D}\kappa)^{-1}A^T)^{-1}], \\ \tilde{\kappa} = \dfrac{1}{n} \mathrm{Tr}[\tilde{D}(\sigma^2(I + \tilde{D}\kappa) + A^T(I + D\tilde{\kappa})^{-1}A)^{-1}]. \end{cases}$$

has a unique pair of solutions $(\kappa, \tilde{\kappa})$, $\kappa > 0, \tilde{\kappa} > 0$. This, of course, is of practical interest in order to compute the approximant (4.3).

REMARK 3.3 (Limiting behavior, followed). The case with a separable variance profile illustrates quite well the kind of assumptions one needs in order to obtain a limiting behavior for the probability distribution of the eigenvalues of $\Sigma_n \Sigma_n^T$: If $\frac{1}{N} \mathrm{Tr}\, Q(z)$ converges then so does $\frac{1}{N} \mathrm{Tr}\, T(z)$ by Theorem 2.5. This in turn should imply the convergence of the right-hand side of (3.3). But, unless $D = I$ and $\tilde{D} = I$ or $A$ pseudo-diagonal, such a convergence relies on an intricate assumption between the limiting distribution of $D$, $\tilde{D}$ and the spectral properties of $A$. Otherwise stated, it is far from obvious to obtain a limiting behavior in the noncentered case even when the variance profile is separable.

3.3. *Additional information about the bias of the deterministic equivalent in a simple case.* Let $X_n$ be an $N \times n$ matrix with i.i.d. entries and denote by $\Delta_n = \text{diag}(\lambda_1, \ldots, \lambda_N)$ where $\frac{1}{N} \sum_{i=1}^{N} \delta_{\lambda_i^2} \to H(d\lambda)$. We consider the matrix

$$Y_n = \frac{1}{\sqrt{n}} \Delta_n X_n.$$

Assume that $\frac{N}{n} \to c \in (0, 1)$ then it is well known (see, e.g., [25]) that the spectral distribution of $Y_n Y_n^*$ converges toward a probability distribution $\mu$ whose Stieltjes transform $f$ satisfies

$$(3.5) \qquad f(z) = \int \frac{H(d\lambda)}{\lambda(1 - c - czf(z)) - z}, \qquad z \in \mathbb{C} - \mathbb{R}^+.$$

We shall prove, in this simple setting, that the deterministic approximation $f_n(z) = \frac{1}{N} \sum_{i=1}^{N} \psi_i(z)$ of the empirical Stieltjes transform $m_n(z) = \frac{1}{N} \text{Tr}(Y_n Y_n^* - zI)^{-1}$ satisfies the following discretized version of (3.5), where $H(d\lambda)$ is replaced by $H_n(d\lambda) = \frac{1}{N} \sum_{i=1}^{N} \delta_{\lambda_i^2}(d\lambda)$ and $c$, by $c_n = \frac{N}{n}$:

$$(3.6) \qquad f_n(z) = \frac{1}{N} \sum_{i=1}^{N} \frac{1}{\lambda_i^2(1 - c_n - c_n z f_n(z)) - z}, \qquad z \in \mathbb{C} - \mathbb{R}^+.$$

In this simple case, the variance profile $\sigma(i, j)$ is equal to $\lambda_i$ and from equations (2.3) and (2.4) we can write:

$$\psi_i(z) = \frac{-1}{z(1 + (\lambda_i^2/n) \sum_{j=1}^{n} \tilde{\psi}_j(z))}, \qquad \tilde{\psi}_j(z) = \frac{-1}{z(1 + (1/n) \sum_{i=1}^{N} \lambda_i^2 \psi_i(z))}$$

or equivalently

$$\psi_i(z) \left(1 + \frac{\lambda_i^2}{n} \sum_{j=1}^{n} \tilde{\psi}_j(z)\right) = -\frac{1}{z}, \qquad \tilde{\psi}_j(z) \left(1 + \frac{1}{n} \sum_{i=1}^{N} \lambda_i^2 \psi_i(z)\right) = -\frac{1}{z}.$$

Eliminating (after immediate manipulation) $\sum_i \lambda_i^2 \psi_i \sum_j \tilde{\psi}_j$ in both equations yields

$$\frac{1}{n} \sum_{j=1}^{n} \tilde{\psi}_j(z) = \frac{c_n}{N} \sum_{i=1}^{N} \psi_i(z) + (1 - c_n)\left(-\frac{1}{z}\right),$$

where $c_n = \frac{N}{n}$. The deterministic approximation $f_n(z)$ thus satisfies (3.6).

One can notice that the centering element in the CLT established by Bai and Silverstein in [2] is based on this kind of deterministic approximation. Moreover, the limiting bias

$$\beta_n = \mathbb{E}m_n(z) - f_n(z)$$

is computed in ([2], Section 4) when $n \to \infty$. According to these results, it turns out that:

1. The bias $\beta_n$ is of order $\frac{K}{n}$ if the entries of $X_n$ are real, $\mathbb{E}X_{ij}^4 = 3$.
2. The bias $\beta_n$ is of order $o(\frac{1}{n})$ if the entries of $X_n$ are complex, $\mathbb{E}|X_{ij}|^4 = 2$ and $\mathbb{E}X_{ij}^2 = 0$.

### 3.4. *The complex case and the case of a Gaussian stationary field.*

*The complex case.* In the case where the entries of $X_n$ and $A_n$ are complex, the interest lies in the spectrum of $(Y_n + A_n)(Y_n + A_n)^*$. All the results established in the real case hold in the complex case, with the following slight adaptations:

1. In Assumption A-1, the random variable $X_{ij}^n$ is complex. Replace $\mathbb{E}X_{ij}^2 = 1$ by $\mathbb{E}|X_{ij}|^2 = 1$ and $\mathbb{E}X_{ij}^2 = 0$.
2. In Assumption A-3, the entries of $A_n$ are complex.
3. In Theorem 2.4, $T$ and $\tilde{T}$ must be replaced by the following:

$$(3.7) \quad T(z) = (\Psi^{-1}(z) - zA\tilde{\Psi}(z)A^*)^{-1}, \qquad \tilde{T}(z) = (\tilde{\Psi}^{-1}(z) - zA^*\Psi(z)A)^{-1}.$$

4. In Theorem 2.5, $Q_n$ and $\tilde{Q}_n$ denote

$$Q_n(z) = (\Sigma_n\Sigma_n^* - zI_N)^{-1}, \qquad \tilde{Q}_n(z) = (\Sigma_n^*\Sigma_n - zI_n)^{-1}.$$

Similar adaptations must be made in the remainder of the article.

*The case of a Gaussian stationary field.* Following [12], it is possible to rely on Theorem 2.5 to provide a deterministic equivalent in the case where matrix $Y_n$ is replaced by $Z_n = (Z_{j_1j_2}^n)$ where

$$(3.8) \quad Z_{j_1j_2}^n = \frac{1}{\sqrt{n}} \sum_{(k_1,k_2)\in\mathbb{Z}^2} h(k_1, k_2)U(j_1 - k_1, j_2 - k_2),$$

where $h$ is a deterministic complex summable sequence and $(U(j_1, j_2); (j_1, j_2) \in \mathbb{Z}^2)$ is a sequence of independent complex Gaussian random variables with zero mean and unit variance, that is, $\mathbb{E}U(j_1, j_2)^2 = 0$, $\mathbb{E}|U(j_1, j_2)|^2 = 1$ [otherwise stated, $U(j_1, j_2) = V + \mathbf{i}W$, where $V$ and $W$ are independent and $\mathcal{N}(0, \frac{1}{\sqrt{2}})$-distributed]. Denote by

$$\Phi(t_1, t_2) = \sum_{(\ell_1,\ell_2)\in\mathbb{Z}^2} h(\ell_1, \ell_2)e^{2\pi\mathbf{i}(\ell_1 t_1 - \ell_2 t_2)}.$$

Function $|\Phi|$ is in particular real and bounded due to the fact that $h$ is a summable sequence. Denote by $F_p$ the $p \times p$ Fourier matrix $F_p = (F^p_{j_1,j_2})_{0 \leq j_1, j_2 < p}$ defined by

$$F^p_{j_1,j_2} = \frac{1}{\sqrt{p}} \exp 2\mathbf{i}\pi \left( \frac{j_1 j_2}{p} \right).$$

We can now state:

**THEOREM 3.4.** *Let $\Sigma_n = Z_n + B_n$ where $Z_n$'s entries are given by (3.8) and $B_n$ is a complex $N \times n$ matrix which satisfies Assumption A-3. Denote by $Q(z) = (\Sigma_n \Sigma^*_n - z I_N)^{-1}$ and $\tilde{Q}_n(z) = (\Sigma^*_n \Sigma_n - z I_n)^{-1}$. Then*

$$\lim_{n \to \infty, N/n \to c} \left( \frac{1}{N} \operatorname{Tr} Q(z) - \frac{1}{N} \operatorname{Tr} T(z) \right) = 0$$

*and*

$$\lim_{n \to \infty, N/n \to c} \left( \frac{1}{n} \operatorname{Tr} \tilde{Q}(z) - \frac{1}{n} \operatorname{Tr} \tilde{T}(z) \right) = 0$$

*for every $z \in \mathbb{C} - \mathbb{R}^+$ where $T$ and $\tilde{T}$ are given by (3.7), where the variance profile $\sigma_{ij}(n)$ must be replaced by $|\Phi(\frac{i}{N}, \frac{j}{n})|$, and matrix $A$ by $F_N B F^*_n$. The convergence holds in probability.*

Lemma 3.3 in [12] [whose assumption is fulfilled due to Assumption A-3] together with Proposition 2.2 part 4 immediately yield the desired result.

**4. Application to digital communication: approximation of the mutual information of a wireless MIMO channel.** The mutual information is the maximum number of bits per second per Hertz per antenna that can be transmitted reliably on a MIMO channel with channel matrix $H_n$. It is equal to

$$(4.1) \qquad C_n(\sigma^2) = \frac{1}{N} \mathbb{E} \log \det \left( I_N + \frac{H_n H^*_n}{\sigma^2} \right),$$

where $\sigma^2$ represents the variance of an additive noise corrupting the received signals. The mutual information $C_n(\sigma^2)$ is related to $\frac{1}{N} \operatorname{Tr}(H_n H^*_n + \sigma^2 I_N)^{-1}$ by the formula

$$\frac{\partial C_n}{\partial \sigma^2} = \frac{1}{N} \mathbb{E} \operatorname{Tr}(H_n H^*_n + \sigma^2 I_N)^{-1} - \frac{1}{\sigma^2}$$

or equivalently by

$$C_n(\sigma^2) = \int_{\sigma^2}^{+\infty} \left( \frac{1}{\omega} - \frac{1}{N} \mathbb{E} \operatorname{Tr}(H_n H^*_n + \omega I_N)^{-1} \right) d\omega,$$

which follows from (4.1) by Fubini's theorem. In fact, $\frac{1}{\omega} - \frac{1}{N}\operatorname{Tr}(H_n H_n^* + \omega I)^{-1} = \frac{1}{\omega} - \int_0^\infty \frac{\mathbb{P}_n^H(d\lambda)}{\lambda + \omega} \geq 0$, where $\mathbb{P}_n^H$ stands for the empirical distribution of the eigenvalues of $H_n H_n^*$. In certain cases, channel matrix $H_n$ is unitarily equivalent to $\Sigma_n = Y_n + A_n$ which has complex entries. Without loss of generality, we shall work with matrix $\Sigma_n$ with real entries in order to remain consistent with the general exposition. Showing that $\frac{1}{N}\operatorname{Tr}(\Sigma_n \Sigma_n^T + \omega I_N)^{-1} \simeq \frac{1}{N}\operatorname{Tr} T_n(-\omega)$ for the deterministic matrix-valued function $T_n(z)$ defined in Theorem 2.4 allows one to approximate $C_n(\sigma^2)$ by $\overline{C}_n(\sigma^2) = \int_{\sigma^2}^{+\infty}(\frac{1}{\omega} - \frac{1}{N}\operatorname{Tr} T_n(-\omega))\,d\omega$. This approximant can be written more explicitly.

THEOREM 4.1. *Assume that Assumptions* A-1, A-2 *and* A-3 *hold and denote by*

$$(4.2) \qquad \overline{C}_n(\sigma^2) = \int_{\sigma^2}^{+\infty}\left(\frac{1}{\omega} - \frac{1}{N}\operatorname{Tr} T_n(-\omega)\right)d\omega,$$

*where $T$ is given by Theorem* 2.4. *Then the following limit holds true:*

$$C_n(\sigma^2) - \overline{C}_n(\sigma^2) \underset{n \to +\infty, N/n \to c}{\longrightarrow} 0,$$

*where $\sigma^2 \in \mathbb{R}^+$. Moreover,*

$$
\begin{aligned}
(4.3) \qquad \overline{C}_n(\sigma^2) &= \frac{1}{N}\log\det\left[\frac{\Psi(-\sigma^2)^{-1}}{\sigma^2} + A\tilde{\Psi}(-\sigma^2)A^T\right] \\
&\quad + \frac{1}{N}\log\det\frac{\tilde{\Psi}(-\sigma^2)^{-1}}{\sigma^2} - \frac{\sigma^2}{nN}\sum_{i,j}\sigma_{ij}^2 T_{ii}(-\sigma^2)\tilde{T}_{jj}(-\sigma^2).
\end{aligned}
$$

Proof of Theorem 4.1 is postponed to Appendix C. In certain cases, the study of the behavior of $\overline{C}_n(\sigma^2)$ is simpler than the behavior of $C_n(\sigma^2)$, and allows one to get some insight on the behavior of the mutual information of certain MIMO wireless channels (see, e.g., [7] for preliminary results).

REMARK 4.2. Equation (4.3) has already been established in the zero mean case ($A_n = 0$): the centered case with no variance profile has been studied by Verdú and Shamai [27], the centered case with a variance profile by Sengupta and Mitra [22] (with a separable variance profile) and Tulino and Verdú ([26], Theorem 2.44).

**5. Proof of Theorem 2.4.** We will first prove the uniqueness of the solutions to the system (2.3) and (2.4) within the class $\mathcal{S}(\mathbb{R}^+)$ and then prove the existence of such solutions. The following proposition will be useful.

PROPOSITION 5.1 (Matrix-valued Stieltjes transforms).   *Let $A_n$ be an $N \times n$ real matrix. Let $\gamma_i \in \mathcal{S}(\mathbb{R}^+), 1 \le i \le N$, and $\tilde{\gamma}_j \in \mathcal{S}(\mathbb{R}^+), 1 \le j \le n$. Denote by*

$$\Gamma = \mathrm{diag}(\gamma_i, 1 \le j \le N), \qquad \tilde{\Gamma} = \mathrm{diag}(\tilde{\gamma}_j, 1 \le i \le n),$$

$$\Upsilon = (\Gamma^{-1} - zA\tilde{\Gamma}A^T)^{-1}, \qquad \tilde{\Upsilon} = (\tilde{\Gamma}^{-1} - zA^T\Gamma A)^{-1}.$$

*Then:*

1. *The matrix-valued functions $\Upsilon$ and $\tilde{\Upsilon}$ are analytic over $\mathbb{C} - \mathbb{R}^+$.*
2. *If $z \in \mathbb{C}^+$, then $\mathrm{Im}\,\Upsilon(z) \ge 0$ and $\mathrm{Im}\,z\Upsilon(z) \ge 0$. Moreover, there exists a positive $N \times N$ matrix-valued measure $\mu = (\mu_{ij})$ and a positive $n \times n$ matrix-valued measure $\tilde{\mu} = (\tilde{\mu}_{ij})$ such that*

$$\mu(\mathbb{R}^+) = I_N, \qquad \tilde{\mu}(\mathbb{R}^+) = I_n \quad and$$

$$\Upsilon(z) = \int_{\mathbb{R}^+} \frac{\mu(d\lambda)}{\lambda - z}, \qquad \tilde{\Upsilon}(z) = \int_{\mathbb{R}^+} \frac{\tilde{\mu}(d\lambda)}{\lambda - z}$$

   *for $z \in \mathbb{C} - \mathbb{R}^+$. In particular, $\frac{1}{N}\mathrm{Tr}\,\Upsilon(z)$ and $\frac{1}{n}\mathrm{Tr}\,\tilde{\Upsilon}(z)$ belong to $\mathcal{S}(\mathbb{R}^+)$.*
3. *The following inequalities hold true:*

$$\forall z \in \mathbb{C}^+ \qquad \Upsilon(z)\Upsilon^*(z) \le \frac{I_N}{(\mathrm{Im}\,z)^2} \quad and \quad \tilde{\Upsilon}(z)\tilde{\Upsilon}^*(z) \le \frac{I_n}{(\mathrm{Im}\,z)^2}.$$

4. *Let $D_j$ and $\tilde{D}_i$ be defined by* (2.1). *Denote by*

$$\gamma_i^{(2)}(z) = \frac{-1}{z(1 + (1/n)\,\mathrm{Tr}(\tilde{D}_i\tilde{\Upsilon}(z)))},$$

$$\tilde{\gamma}_j^{(2)}(z) = \frac{-1}{z(1 + (1/n)\,\mathrm{Tr}(D_j\Upsilon(z)))},$$

   *where $1 \le i \le N$ and $1 \le j \le n$. Then $\gamma_i^{(2)}$ and $\tilde{\gamma}_j^{(2)}$ are analytic over $\mathbb{C} - \mathbb{R}^+$ and belong to $\mathcal{S}(\mathbb{R}^+)$.*

PROOF.   We only prove the stated results for $\Upsilon$, the adaptation to $\tilde{\Upsilon}$ being straightforward.

Since $\gamma_i \in \mathcal{S}(\mathbb{R}^+)$, we have $\gamma_i(z) \ne 0$ over $\mathbb{C} - \mathbb{R}^+$. In particular, $\Gamma^{-1} - zA\tilde{\Gamma}A^T$ is analytic over $\mathbb{C} - \mathbb{R}^+$. In order to prove that $\Upsilon$ is analytic, it is sufficient to prove that $\det(\Gamma^{-1} - zA\tilde{\Gamma}A^T) \ne 0$ for every $z \in \mathbb{C} - \mathbb{R}^+$ or equivalently, to prove that if $h \in \mathbb{C}^N$ is such that $(\Gamma^{-1} - zA\tilde{\Gamma}A^T)h = 0$, then $h$ must be equal to zero. Let $h$ be such that $(\Gamma^{-1} - zA\tilde{\Gamma}A^T)h = 0$, then

(5.1)
$$\begin{cases} h^*(\Gamma^{-1} - zA\tilde{\Gamma}A^T)h = 0, \\ h^*(\Gamma^{-1*} - z^*A\tilde{\Gamma}^*A^T)h = 0. \end{cases}$$

We study separately the case $z \in (-\infty, 0)$ and the cases $z \in \mathbb{C}^+$ and $z \in \mathbb{C}^-$. If $z \in (-\infty, 0)$, then $\Gamma^{-1}(z) \geq |z| I_N$ and $-z\tilde{\Gamma}(z) \geq 0$ and the first equation in (5.1) yields:

$$0 = h^*(\Gamma^{-1}(z) - zA\tilde{\Gamma}(z)A^T)h \geq \|h\|^2 |z|$$

which implies $h = 0$. If $z \in \mathbb{C}^+$, then $\text{Im}(z\tilde{\Gamma}(z)) \geq 0$ and $\text{Im}\,\Gamma^{-1}(z) \leq -\text{Im}(z)I_n$ by Proposition 2.2 part 1. The system (5.1) yields $h^*(\text{Im}(\Gamma^{-1}(z)) - A\,\text{Im}(z\tilde{\Gamma}(z)) \times A^T)h = 0$. Necessarily,

$$0 = h^*(\text{Im}(\Gamma^{-1}(z)) - A\,\text{Im}(z\tilde{\Gamma}(z))A^T)h \leq -(\text{Im}\,z)\|h\|^2,$$

which yields $h = 0$. If $z \in \mathbb{C}^-$, the same kind of arguments yields $h = 0$. Therefore, $\Upsilon$ is analytic over $\mathbb{C} - \mathbb{R}^+$ and part 1 is proved.

In order to apply Proposition 2.2 part 3, we now check that $\text{Im}\,\Upsilon(z) \geq 0$ and $\text{Im}\,z\Upsilon(z) \geq 0$ for $z \in \mathbb{C}^+$:

$$\text{Im}\,\Upsilon(z) = \frac{\Upsilon(z) - \Upsilon^*(z)}{2\mathbf{i}}$$

$$\overset{\text{(a)}}{=} \Upsilon(z)\left(\frac{\Gamma^{-1*}(z) - \Gamma^{-1}(z)}{2\mathbf{i}} + A\left(\frac{z\tilde{\Gamma}(z) - z^*\tilde{\Gamma}^*(z)}{2\mathbf{i}}\right)A^T\right)\Upsilon^*(z)$$

$$= \Upsilon(z)(-\text{Im}\,\Gamma^{-1}(z) + A\,\text{Im}(z\tilde{\Gamma}(z))A^T)\Upsilon^*(z) \geq 0,$$

where (a) follows from the resolvent identity and the last inequality follows from the fact that $\gamma_i$ and $\tilde{\gamma}_j$ belong to $\mathcal{S}(\mathbb{R}^+)$ and from Proposition 2.2 part 1.

One can check similarly that $\text{Im}(z\Upsilon(z)) \geq 0$. Therefore Proposition 2.2 part 3 yields $\Upsilon(z) = C + \int_{\mathbb{R}^+} \frac{\mu(d\lambda)}{\lambda - z}$, where $C = (C_{ij})$ is an $N \times N$ nonnegative matrix and $\mu$ is an $N \times N$ matrix-valued measure. Let us now prove that

$$(5.2) \qquad \lim_{y \to +\infty} \mathbf{i}y\Upsilon(\mathbf{i}y) = -I_N.$$

First write

$$(5.3) \quad \Upsilon(z) = (\Gamma^{-1}(z) - zA\tilde{\Gamma}(z)A^T)^{-1} = (I_N - \Gamma(z)Az\tilde{\Gamma}(z)A^T)^{-1}\Gamma(z).$$

Since $\Gamma$ and $\tilde{\Gamma}$ are diagonal matrices whose entries belong to $\mathcal{S}(\mathbb{R}^+)$, we have

$$(5.4) \qquad \lim_{y \to +\infty} \mathbf{i}y\Gamma(\mathbf{i}y) = -I_N, \qquad \lim_{y \to +\infty} \mathbf{i}y\tilde{\Gamma}(\mathbf{i}y) = -I_n.$$

Therefore, it is straightforward to show

$$(5.5) \qquad \|\Gamma(\mathbf{i}y)A(\mathbf{i}y\tilde{\Gamma}(\mathbf{i}y))A^T\|_{\text{sp}} \underset{y \to +\infty}{\longrightarrow} 0.$$

Plugging (5.4) and (5.5) into (5.3), we get (5.2).

In order to prove that $C = 0$ and $\mu(\mathbb{R}^+) = I_N$, we introduce the complex functions: $f_i(z) = \Upsilon_{ii}(z), 1 \leq i \leq N$. These functions are analytic over $\mathbb{C}^+$ and satisfy: $\forall z \in \mathbb{C}^+, \operatorname{Im} f_i(z) \geq 0, \operatorname{Im} z f_i(z) \geq 0$. Moreover (5.2) implies that $\lim_{y \to +\infty} -\mathbf{i} y f_i(\mathbf{i} y) = 1$; therefore Proposition 2.2 part 2 yields the existence and uniqueness of $\nu_i \in \mathcal{P}(\mathbb{R}^+)$ such that

$$f_i(z) = \int_{\mathbb{R}^+} \frac{\nu_i(d\lambda)}{\lambda - z}.$$

On the other hand,

$$f_i(z) = \Upsilon_{ii}(z) = C_{ii} + \int_{\mathbb{R}^+} \frac{\mu_{ii}(d\lambda)}{\lambda - z}, \qquad 1 \leq i \leq N.$$

From this, we deduce that $\mu_{ii} = \nu_i \in \mathcal{P}(\mathbb{R}^+)$ and that $C_{ii} = 0$. This implies that $C = 0$ since $C \geq 0$. As $\mu$ is a positive matrix-valued measure, $|\mu_{k\ell}| \leq \frac{1}{2}(\mu_{kk} + \mu_{\ell\ell})$, $|\mu_{k\ell}|(\mathbb{R}^-) = 0$ and $|\mu_{k\ell}|(\mathbb{R}^+) < \infty$. Assume that $k \neq \ell$. Equation (5.2) yields

$$(5.6) \qquad \lim_{y \to +\infty} \mathbf{i} y \Upsilon_{k\ell}(\mathbf{i} y) = 0.$$

But

$$\mathbf{i} y \Upsilon_{k\ell}(\mathbf{i} y) = -\int_{\mathbb{R}^+} \frac{y^2}{\lambda^2 + y^2} \mu_{k\ell}(d\lambda) + \mathbf{i} \int_{\mathbb{R}^+} \frac{y\lambda}{\lambda^2 + y^2} \mu_{k\ell}(d\lambda).$$

In particular, (5.6) implies that $\lim_{y \to +\infty} \int_{\mathbb{R}^+} \frac{y^2}{\lambda^2 + y^2} \mu_{k\ell}(d\lambda) = 0$. But the dominated convergence theorem implies that $\lim_{y \to +\infty} \int_{\mathbb{R}^+} \frac{y^2}{\lambda^2 + y^2} \mu_{k\ell}(d\lambda) = \mu_{k\ell}(\mathbb{R}^+)$. Necessarily $\mu_{k\ell}(\mathbb{R}^+) = 0$, $\mu(\mathbb{R}^+) = I_N$ and part 2 is proved for $z \in \mathbb{C}^+$. Since both $\Upsilon$ and $\int \frac{\mu(d\lambda)}{\lambda - z}$ are analytic over $\mathbb{C} - \mathbb{R}^+$ and coincide over $\mathbb{C}^+$, they are equal. Taking the trace, we immediately obtain that $\frac{1}{N} \operatorname{Tr} \Upsilon(z) \in \mathcal{S}(\mathbb{R}^+)$ and part 2 is proved.

In order to prove part 3, we shall prove the two following inequalities:

$$(5.7) \qquad \forall z \in \mathbb{C}^+ \qquad \frac{\Upsilon(z) - \Upsilon(z)^*}{2\mathbf{i}} \leq \frac{I_N}{\operatorname{Im}(z)} \quad \text{and}$$

$$\frac{\Upsilon(z) - \Upsilon(z)^*}{2\mathbf{i}} \geq \Upsilon(z) \Upsilon^*(z)(\operatorname{Im}(z)).$$

The first one is straightforward:

$$\frac{\Upsilon(z) - \Upsilon(z)^*}{2\mathbf{i}} = \operatorname{Im} \Upsilon(z) = \operatorname{Im}(z) \int_{\mathbb{R}^+} \frac{\mu(d\lambda)}{|\lambda - z|^2}.$$

Since $\mu(A) \geq 0$ for every Borel set $A$, so is $\int f(\lambda) \mu(d\lambda)$ whenever $f \geq 0$. Applying this property to $f(\lambda) = \frac{1}{\operatorname{Im}(z)^2} - \frac{1}{|\lambda - z|^2} \geq 0$ and using the fact that $\mu(\mathbb{R}^+) = I_N$, we obtain

$$\int_{\mathbb{R}^+} \frac{\mu(d\lambda)}{|\lambda - z|^2} \leq \frac{I_N}{\operatorname{Im}(z)^2},$$

which yields $\frac{\Upsilon(z) - \Upsilon(z)^*}{2\mathbf{i}} \leq \frac{I_N}{\operatorname{Im}(z)}$. Using the resolvent identity, we get

$$\frac{\Upsilon(z) - \Upsilon(z)^*}{2\mathbf{i}} = \Upsilon(z) \Big( \frac{\Gamma^{-1}(z)^* - \Gamma^{-1}(z)}{2\mathbf{i}} + A \operatorname{Im}(z\tilde{\Gamma}(z)) A^T \Big) \Upsilon^*(z).$$

Let $z \in \mathbb{C}^+$, then

$$\frac{\Gamma^{-1}(z)^* - \Gamma^{-1}(z)}{2\mathbf{i}} = -\operatorname{diag}\Big( \operatorname{Im}\Big( \frac{1}{\gamma_i(z)} \Big), 1 \leq i \leq N \Big) \geq \operatorname{Im}(z) I$$

by Proposition 2.2 part 1. This and the fact that $\operatorname{Im}(z\tilde{\Gamma}(z)) \geq 0$ yields the second inequality in (5.7) and the proof of part 3 is complete.

Concerning part 4, we only prove the statement for $\gamma_i^{(2)} \in \mathcal{S}(\mathbb{R}^+)$ (since the same arguments yield the desired results for $\tilde{\gamma}_j^{(2)}$). For this we rely on the characterization given in Proposition 2.2 part 2. The quantity $z(1 + \frac{1}{n} \operatorname{Tr}(\tilde{D}_i\tilde{\Upsilon}(z)))$ is analytic over $\mathbb{C}^+$, and its imaginary part is greater than $\operatorname{Im}(z)$ there. Therefore, it does not vanish over $\mathbb{C}^+$, and $\gamma_i^{(2)}$ is analytic over $\mathbb{C}^+$. Similar arguments yield that $\gamma_i^{(2)}$ is analytic over $\mathbb{C}^-$ and $(-\infty, 0)$. Let us now compute $\operatorname{Im} \gamma_i^{(2)}(z)$:

$$\frac{\gamma_i^{(2)}(z) - \gamma_i^{(2)}(z)^*}{2\mathbf{i}} = \frac{\operatorname{Im}(z)}{|z(1 + (1/n) \operatorname{Tr}(\tilde{D}_i\tilde{\Upsilon}(z)))|^2}$$

$$+ \frac{\operatorname{Tr}[\tilde{D}_i(z\tilde{\Upsilon}(z) - z^*\tilde{\Upsilon}(z)^*)/(2\mathbf{i})]}{n|z(1 + (1/n) \operatorname{Tr}(\tilde{D}_i\tilde{\Upsilon}(z)))|^2} \geq 0$$

for $z \in \mathbb{C}^+$ by Proposition 5.1 part 2. Similarly, one can prove that $\operatorname{Im} z\gamma_i^{(2)}(z) \geq 0$ for $z \in \mathbb{C}^+$. Finally,

$$\mathbf{i}y\gamma_i^{(2)}(\mathbf{i}y) = \frac{-1}{1 + (1/n) \operatorname{Tr}(\tilde{D}_i\tilde{\Upsilon}(\mathbf{i}y))} \xrightarrow[y\to\infty]{} -1$$

since $|\operatorname{Tr}(\tilde{D}_i\tilde{\Upsilon}(\mathbf{i}y))| \leq n\sigma_{\max}^2 \|\Upsilon(\mathbf{i}y)\|_{\operatorname{sp}}$ and $\|\Upsilon(\mathbf{i}y)\|_{\operatorname{sp}} \leq \frac{1}{y}$. The proof of part 4 is completed and Proposition 5.1 is proved. $\square$

PROPOSITION 5.2 (Uniqueness of solutions). *Assume that $(\psi_i, \tilde{\psi}_j; 1 \leq i \leq N, 1 \leq j \leq n) \in \mathcal{S}(\mathbb{R}^+)^{N+n}$ and $(\phi_i, \tilde{\phi}_j; 1 \leq i \leq N, 1 \leq j \leq n) \in \mathcal{S}(\mathbb{R}^+)^{N+n}$ are two sets of solutions to the system* (2.3) *and* (2.4)*. Then these solutions are equal.*

PROOF.  Denote by

$$\Phi(z) = \operatorname{diag}(\phi_i(z), 1 \leq i \leq N), \qquad \tilde{\Phi}(z) = \operatorname{diag}(\tilde{\phi}_i(z), 1 \leq j \leq n),$$

$$T^\phi(z) = (\Phi^{-1}(z) - zA\tilde{\Phi}(z)A^T)^{-1}, \qquad \tilde{T}^\phi(z) = (\tilde{\Phi}^{-1}(z) - zA^T\Phi(z)A)^{-1},$$

$$\Psi(z) = \mathrm{diag}(\psi_i(z), 1 \le i \le N), \qquad \tilde{\Psi}(z) = \mathrm{diag}(\tilde{\psi}_i(z), 1 \le j \le n),$$

$$T^\psi(z) = (\Psi^{-1}(z) - zA\tilde{\Psi}(z)A^T)^{-1}, \qquad \tilde{T}^\psi(z) = (\tilde{\Psi}^{-1}(z) - zA^T\Psi(z)A)^{-1}.$$

Let us compute $\psi_i - \phi_i$ for $z \in \mathbb{C}^+$.

$$\psi_i(z) - \phi_i(z) = \frac{-1}{z(1 + (1/n)\,\mathrm{Tr}\,\tilde{D}_i\tilde{T}^\psi(z))} + \frac{1}{z(1 + (1/n)\,\mathrm{Tr}\,\tilde{D}_i\tilde{T}^\phi(z))}.$$

If $z \in \mathbb{C}^+$ then

$$|\psi_i(z)|^{-1} \ge \left| z\left(1 + \frac{1}{n}\,\mathrm{Tr}\,\tilde{D}_i\tilde{T}^\psi(z)\right) \right|$$

$$\ge \mathrm{Im}\left( z\left(1 + \frac{1}{n}\,\mathrm{Tr}\,\tilde{D}_i\tilde{T}^\psi(z)\right) \right) \overset{\text{(a)}}{\ge} \mathrm{Im}(z),$$

where (a) follows from Proposition 5.1 part 2. The same argument yields a similar bound for $\phi_i(z)$, that is $|\phi_i^{-1}(z)| \ge \mathrm{Im}(z)$, and

$$(5.8) \qquad |\psi_i(z) - \phi_i(z)| \le \frac{|z|}{(\mathrm{Im}\,z)^2}\left| \frac{1}{n}\mathrm{Tr}\tilde{D}_i(\tilde{T}^\phi - \tilde{T}^\psi) \right|.$$

A standard computation involving the resolvent identity yields

$$\mathrm{Tr}\,\tilde{D}_i(\tilde{T}^\phi - \tilde{T}^\psi)$$

$$= -\mathrm{Tr}\,\tilde{D}_i\tilde{T}^\phi(\tilde{\Phi}^{-1} - \tilde{\Psi}^{-1} - zA^T(\Phi - \Psi)A)\tilde{T}^\psi$$

$$= -\mathrm{Tr}\,\tilde{D}_i\tilde{T}^\phi(\tilde{\Phi}^{-1} - \tilde{\Psi}^{-1})\tilde{T}^\psi + z\,\mathrm{Tr}\,\tilde{D}_i\tilde{T}^\phi(A^T(\Phi - \Psi)A)\tilde{T}^\psi$$

$$= \mathrm{Tr}\,\tilde{D}_i\tilde{T}^\phi\tilde{\Phi}^{-1}(\tilde{\Phi} - \tilde{\Psi})\tilde{\Psi}^{-1}\tilde{T}^\psi + \sum_{k=1}^N z\,\mathrm{Tr}\,\tilde{D}_i\tilde{T}^\phi\tilde{\mathbf{a}}_k^T\tilde{\mathbf{a}}_k(\phi_k - \psi_k)\tilde{T}^\psi$$

(recall that $\tilde{\mathbf{a}}_k$ is the $k$th row of $A$). We introduce the following maximum:

$$\mathbf{M} = \sup\{|\phi_i - \psi_i|, |\tilde{\phi}_j - \tilde{\psi}_j|, 1 \le i \le N, 1 \le j \le n\}.$$

We also define $d_n$ by $d_n = \max(\frac{n}{N}, \frac{N}{n})$. Note that $d_n \ge 1$. Then,

$$\left| \frac{1}{n}\mathrm{Tr}\,\tilde{D}_i(\tilde{T}^\phi - \tilde{T}^\psi) \right| \le \mathbf{M} d_n \sigma_{\max}^2 \|\tilde{T}^\phi\|_{\mathrm{sp}} \|\tilde{T}^\psi\|_{\mathrm{sp}} \{\|\tilde{\Phi}^{-1}\|_{\mathrm{sp}} \|\tilde{\Psi}^{-1}\|_{\mathrm{sp}} + |z|\mathbf{a}_{\max,n}^2\}.$$

First notice that Proposition 5.1 part 3 yields $\|T^\phi\|_{\mathrm{sp}} \le \frac{1}{\mathrm{Im}\,z}$ (same inequality for $T^\psi$, $\tilde{T}^\phi$ and $\tilde{T}^\psi$). Since $\frac{1}{\psi_i(z)} = -z(1 + \frac{1}{n}\mathrm{Tr}\,\tilde{D}_i\tilde{T}^\psi(z))$, we have

$$\|\Psi^{-1}\|_{\mathrm{sp}} \le |z|(1 + d_n\sigma_{\max}^2\|\tilde{T}^\psi(z)\|_{\mathrm{sp}}) \le |z|\left(1 + d_n\frac{\sigma_{\max}^2}{\mathrm{Im}\,z}\right).$$

One can show similarly that $\max(\|\Phi^{-1}\|_{\mathrm{sp}}, \|\tilde{\Psi}^{-1}\|_{\mathrm{sp}}, \|\tilde{\Phi}^{-1}\|_{\mathrm{sp}}) \le |z|(1 + d_n\frac{\sigma_{\max}^2}{\mathrm{Im}\,z})$. Therefore,

$$(5.9) \quad \left| \frac{1}{n}\mathrm{Tr}\,\tilde{D}_i(\tilde{T}^\phi - \tilde{T}^\psi) \right| \le \frac{\sigma_{\max}^2 d_n|z|}{(\mathrm{Im}\,z)^2}\left\{ |z|\left(1 + d_n\frac{\sigma_{\max}^2}{\mathrm{Im}\,z}\right)^2 + \mathbf{a}_{\max,n}^2 \right\}\mathbf{M}(z).$$

Plugging (5.9) into (5.8), we obtain

$$(5.10) \qquad |\psi_i(z) - \phi_i(z)| \leq \mathcal{E}(z)\mathbf{M}(z),$$

$$\text{where } \mathcal{E}(z) = \frac{|z|^2 d_n \sigma_{\max}^2}{(\operatorname{Im} z)^4}\left\{|z|\left(1 + d_n \frac{\sigma_{\max}^2}{\operatorname{Im} z}\right)^2 + \mathbf{a}_{\max,n}^2\right\}.$$

The same kind of computation yields

$$(5.11) \qquad |\tilde{\psi}_i(z) - \tilde{\phi}_i(z)| \leq \mathcal{E}(z)\mathbf{M}(z).$$

Finally, gathering (5.10) and (5.11) we obtain

$$\mathbf{M}(z) \leq \mathcal{E}(z)\mathbf{M}(z).$$

Let $z$ be such that $|\frac{z}{\operatorname{Im} z}| \leq 2$. For $\operatorname{Im} z$ large enough, one has $\mathcal{E}(z) \leq \frac{1}{2}$, which implies that $\mathbf{M}(z) = 0$. The analyticity of $\phi_i, \psi_i, \tilde{\phi}_j, \tilde{\psi}_j$ implies that $\phi_i = \psi_i$ for $1 \leq i \leq N$ and $\tilde{\phi}_j = \tilde{\psi}_j$ for $1 \leq j \leq n$ on $\mathbb{C} - \mathbb{R}^+$. Proposition 5.2 is proved. $\square$

PROPOSITION 5.3 (Existence of solutions). *There exists $(\psi_i, \tilde{\psi}_j; 1 \leq i \leq N, 1 \leq j \leq n) \in \mathcal{S}(\mathbb{R}^+)^{N+n}$ satisfying* (2.3) *and* (2.4).

PROOF. We construct the desired solution by induction. Let

$$\psi_i^{(0)}(z) = \tilde{\psi}_j^{(0)}(z) = -\frac{1}{z} \qquad \text{for } 1 \leq i \leq N, \ 1 \leq j \leq n.$$

Then $\psi_i^{(0)}$ and $\tilde{\psi}_j^{(0)}$ belong to $\mathcal{S}(\mathbb{R}^+)$. For $p \geq 0$, let

$$\psi_i^{(p+1)}(z) = \frac{-1}{z(1 + (1/n)\operatorname{Tr}(\tilde{D}_i \tilde{T}^{(p)}(z)))} \qquad \text{for } 1 \leq i \leq N,$$

$$\tilde{\psi}_j^{(p+1)}(z) = \frac{-1}{z(1 + (1/n)\operatorname{Tr}(D_j T^{(p)}(z)))} \qquad \text{for } 1 \leq j \leq n,$$

where

$$\Psi^{(p)}(z) = \operatorname{diag}(\psi_i^{(p)}(z), 1 \leq i \leq N),$$

$$\tilde{\Psi}^{(p)}(z) = \operatorname{diag}(\tilde{\psi}_j^{(p)}(z), 1 \leq j \leq n),$$

$$T^{(p)}(z) = (\Psi^{(p)-1}(z) - zA\tilde{\Psi}(z)^{(p)}A^T)^{-1},$$

$$\tilde{T}^{(p)}(z) = (\tilde{\Psi}^{(p)-1}(z) - zA^T\Psi^{(p)}(z)A)^{-1}.$$

By Proposition 5.1 part 4, $(\psi_i^{(p)}, \tilde{\psi}_j^{(p)}; 1 \leq i \leq N, 1 \leq j \leq n)$ are analytic over $\mathbb{C} - \mathbb{R}^+$ and belong to $\mathcal{S}(\mathbb{R}^+)$ for all $p \geq 0$. Denote by

$$\mathbf{M}^{(p)} = \sup\{|\psi_i^{(p+1)} - \psi_i^{(p)}|, |\tilde{\psi}_j^{(p+1)} - \tilde{\psi}_j^{(p)}|, 1 \leq i \leq N, 1 \leq j \leq n\}.$$

The same computations as in the proof of Proposition 5.2 yield

$$\mathbf{M}^{(p+1)}(z) \le \mathcal{E}(z)\mathbf{M}^{(p)}(z),$$

$$\text{where } \mathcal{E}(z) = \frac{|z|^2 d_n \sigma_{\max}^2}{(\operatorname{Im} z)^4}\left\{|z|\left(1 + d_n \frac{\sigma_{\max}^2}{\operatorname{Im} z}\right)^2 + \mathbf{a}_{\max,n}^2\right\}.$$

Let $z \in \mathbb{C}^+$ be such that $|\frac{z}{\operatorname{Im} z}| \le 2$. For $\operatorname{Im} z$ large enough, $\psi_i^{(p)}(z)$ and $\tilde{\psi}_j^{(p)}(z)$ are Cauchy sequences. Denote by $\psi_i(z)$ and $\tilde{\psi}_j(z)$ the corresponding limits. On the other hand, $(\psi_i^{(p)})_p$ is a normal family over $\mathbb{C} - \mathbb{R}^+$ (see, e.g., [21]) since $\psi_i^{(p)}$ is bounded on every compact set included in $\mathbb{C} - \mathbb{R}^+$ uniformly in $p$ [see the third inequality of (2.2)]. Therefore one can extract, by the normal family theorem, a converging subsequence whose limit is analytic over $\mathbb{C} - \mathbb{R}^+$. Since the limit of any converging subsequence is equal to $\psi_i(z)$ in the domain $\{z \in \mathbb{C}^+, |z|/|\operatorname{Im} z| \le 2, \operatorname{Im} z \text{ large enough}\}$, $\psi_i^{(p)}(z)$ converges toward an analytic function on $\mathbb{C} - \mathbb{R}^+$ [that we still denote by $\psi_i(z)$]. One can apply the same arguments for $\tilde{\psi}_j^{(p)}(z)$.

We now prove that $\psi_i$ and $\tilde{\psi}_j$ satisfy (2.3) and (2.4) where $\Psi$, $\tilde{\Psi}$, $T$ and $\tilde{T}$ are well defined. Let $\Psi(z) = \operatorname{diag}(\psi_i(z), 1 \le i \le N)$ and $\tilde{\Psi}(z) = \operatorname{diag}(\tilde{\psi}_i(z), 1 \le j \le n)$. The convergence of $(\psi_i^{(p)})_p$ immediately yields that

$$\operatorname{Im} \psi_i^{(p)}(z) \ge 0 \Rightarrow \operatorname{Im} \psi_i(z) \ge 0,$$

$$\operatorname{Im} z\psi_i^{(p)}(z) \ge 0 \Rightarrow \operatorname{Im} z\psi_i(z) \ge 0,$$

$$|\psi_i^{(p)}(z)| \le \frac{1}{\operatorname{Im} z} \Rightarrow |\psi_i(z)| \le \frac{1}{\operatorname{Im} z}$$

on $\mathbb{C}^+$. It remains to prove that $\lim_{y\to\infty} -\mathbf{i}y\psi_i(\mathbf{i}y) = 1$. The same arguments as in the proof of Proposition 5.1 yield that both

$$T(z) = (\Psi^{-1}(z) - zA\tilde{\Psi}(z)A^T)^{-1} \quad \text{and} \quad \tilde{T}(z) = (\tilde{\Psi}^{-1}(z) - zA^T\Psi(z)A)^{-1}$$

are well defined on $\mathbb{C} - \mathbb{R}^+$. Moreover,

$$T^{(p)}(z) \underset{p\to\infty}{\longrightarrow} T(z) \quad \text{and} \quad \tilde{T}^{(p)}(z) \underset{p\to\infty}{\longrightarrow} \tilde{T}(z)$$

on $\mathbb{C} - \mathbb{R}^+$. We can deduce that

$$\psi_i(z) = \frac{-1}{z(1 + (1/n)\operatorname{Tr}(\tilde{D}_i\tilde{T}(z)))} \qquad \text{for } 1 \le i \le N,$$

$$\tilde{\psi}_j(z) = \frac{-1}{z(1 + (1/n)\operatorname{Tr}(D_j T(z)))} \qquad \text{for } 1 \le j \le n$$

on $\mathbb{C}^+$, and hence on $\mathbb{C} - \mathbb{R}^+$. Therefore, $T^{(p)}(z)T^{(p)}(z)^* \le \frac{I_n}{(\operatorname{Im} z)^2}$ implies that $T(z)T(z)^* \le \frac{I_n}{(\operatorname{Im} z)^2}$ and one can easily prove that

$$\lim_{y\to\infty} -\mathbf{i}y\psi_i(\mathbf{i}y) = 1$$

using the previous representation of $\psi_i$. Hence, $\psi_i$ belongs to $\mathcal{S}(\mathbb{R}^+)$ by Proposition 2.2 part 2. We can prove similarly that $\tilde{\psi}_j \in \mathcal{S}(\mathbb{R}^+)$. Proposition 5.3 is proved. $\quad\square$

The proof of Theorem 2.4 immediately follows from Propositions 5.1, 5.2 and 5.3.

**6. Proof of Theorem 2.5.** For the reader's ease, we recall some of the notation previously introduced.

The resolvents:

$$Q(z) = (\Sigma\Sigma^T - zI)^{-1}, \qquad \tilde{Q}(z) = (\Sigma^T\Sigma - zI)^{-1}.$$

The stochastic intermediate quantities:

$$b_i(z) = \frac{-1}{z(1 + (1/n)\operatorname{Tr}(\tilde{D}_i\tilde{Q}(z)))}, \qquad B(z) = \operatorname{diag}(b_i(z), 1 \le i \le N),$$

$$\tilde{b}_j(z) = \frac{-1}{z(1 + (1/n)\operatorname{Tr}(D_jQ(z)))}, \qquad \tilde{B}(z) = \operatorname{diag}(\tilde{b}_j(z), 1 \le j \le n),$$

$$R(z) = (B^{-1}(z) - zA\tilde{B}(z)A^T)^{-1}, \qquad \tilde{R}(z) = (\tilde{B}^{-1}(z) - zA^TB(z)A)^{-1}.$$

The deterministic quantities:

$$\psi_i(z) = \frac{-1}{z(1 + (1/n)\operatorname{Tr}(\tilde{D}_i\tilde{T}(z)))}, \qquad \tilde{\psi}_j(z) = \frac{-1}{z(1 + (1/n)\operatorname{Tr}(D_jT(z)))},$$

$$\Psi(z) = \operatorname{diag}(\psi_i(z), 1 \le i \le N), \qquad \tilde{\Psi}(z) = \operatorname{diag}(\tilde{\psi}_j(z), 1 \le j \le n),$$

$$T(z) = (\Psi^{-1}(z) - zA\tilde{\Psi}(z)A^T)^{-1}, \qquad \tilde{T}(z) = (\tilde{\Psi}^{-1}(z) - zA^T\Psi(z)A)^{-1}.$$

6.1. *Evaluation of the differences* $\frac{1}{N}\operatorname{Tr}Q(z) - \frac{1}{N}\operatorname{Tr}R(z)$ *and* $\frac{1}{n}\operatorname{Tr}\tilde{Q}(z) - \frac{1}{n}\operatorname{Tr}\tilde{R}(z)$.

LEMMA 6.1. *Let* $U_n = \operatorname{diag}(u_i^n, 1 \le i \le N)$ *be a sequence of* $N \times N$ *diagonal matrices and* $\tilde{U}_n = \operatorname{diag}(\tilde{u}_i^n, 1 \le i \le n)$, *a sequence of* $n \times n$ *diagonal matrices. Assume that there exists* $K_1$ *and* $\tilde{K}_1$ *such that*

$$\sup_{n \ge 1}\max_{1 \le i \le N}|u_i^n| \le K_1 < \infty \quad and \quad \sup_{n \ge 1}\max_{1 \le i \le n}|\tilde{u}_i^n| \le \tilde{K}_1 < \infty.$$

*Then, if* $z \in \mathbb{C}^+$,

(6.1)
$$\mathbb{E}\left|\frac{1}{N}\operatorname{Tr}(Q(z) - R(z))U\right|^{2+\varepsilon/2} \le \frac{K_2}{n^{1+\varepsilon/4}} \quad and$$

$$\mathbb{E}\left|\frac{1}{n}\operatorname{Tr}(\tilde{Q}(z) - \tilde{R}(z))\tilde{U}\right|^{2+\varepsilon/2} \le \frac{\tilde{K}_2}{n^{1+\varepsilon/4}}.$$

*In particular, for each $z \in \mathbb{C}^+$,*

$$\frac{1}{N} \operatorname{Tr}(Q(z) - R(z))U \to 0 \quad and \quad \frac{1}{n} \operatorname{Tr}(\tilde{Q}(z) - \tilde{R}(z))\tilde{U} \to 0$$

*almost surely as $N \to \infty$ and $N/n \to c > 0$.*

The computations in the proof of Lemma 6.1, along the same lines as the computations in [5], require some adaptation due to the fact that $Y_n$ has a variance profile. In this section, we state two intermediate results in Proposition 6.3 and Proposition 6.4. A sketch of the proof of Proposition 6.4 is given while the full proof is postponed to Appendix B.

We first need to introduce related matrix quantities when one column/row is removed. Denote by $\mathbf{y}_j$, $\mathbf{a}_j$ and $\boldsymbol{\xi}_j$ the $j$th column of $Y$, $A$ and $\Sigma$ respectively, and by $A^{(j)}$ and $\Sigma^{(j)}$ the $N \times (n-1)$ matrices that remain after deleting the $j$th column from matrices $A$ and $\Sigma$, respectively. Also, let

$$Q^{(j)}(z) = (\Sigma^{(j)}\Sigma^{(j)T} - zI)^{-1} \quad \text{and} \quad \tilde{Q}^{(j)}(z) = (\Sigma^{(j)T}\Sigma^{(j)} - zI)^{-1}.$$

Let $\tilde{D}_i^{(j)}$ be the $(n-1) \times (n-1)$ matrix where column $j$ and row $j$ have been removed:

$$\tilde{D}_i^{(j)} = \operatorname{diag}(\sigma_{i,\ell}^2, 1 \le \ell \le n, \ell \ne j).$$

Let

(6.2)     $$b_i^{(j)}(z) = \frac{-1}{z(1 + (1/n)\operatorname{Tr}(\tilde{D}_i^{(j)}\tilde{Q}^{(j)}))}, \qquad 1 \le i \le N,$$

$$B^{(j)}(z) = \operatorname{diag}(b_i^{(j)}(z), 1 \le i \le N),$$

(6.3)     $$\tilde{b}_\ell^{(j)}(z) = \frac{-1}{z(1 + (1/n)\operatorname{Tr}(D_\ell Q^{(j)}(z)))}, \qquad 1 \le \ell \le n, \ell \ne j,$$

$$\tilde{B}^{(j)}(z) = \operatorname{diag}(\tilde{b}_\ell^{(j)}(z), 1 \le \ell \le n, \ell \ne j),$$

$$R^{(j)}(z) = (B^{(j)-1}(z) - zA^{(j)}\tilde{B}^{(j)}(z)A^{(j)T})^{-1},$$

$$\tilde{R}^{(j)}(z) = (\tilde{B}^{(j)-1}(z) - zA^{(j)T}B^{(j)}(z)A^{(j)})^{-1}.$$

Note that any random matrix with superscript $(j)$ is independent of $\mathbf{y}_j$. This fact will be frequently used in the sequel. By the matrix inversion lemma (see [16], Appendix A, see also [15], Section 0.7.4), we have

(6.4)     $$Q = Q^{(j)} - \frac{1}{\alpha_j}Q^{(j)}\boldsymbol{\xi}_j\boldsymbol{\xi}_j^T Q^{(j)} \qquad \text{where } \alpha_j = 1 + \boldsymbol{\xi}_j^T Q^{(j)}\boldsymbol{\xi}_j.$$

Recall that $\boldsymbol{\xi}_j = \mathbf{a}_j + \mathbf{y}_j$. Following [5], we introduce the following random variables:

$$\rho_j = \mathbf{a}_j^T Q^{(j)} \mathbf{a}_j, \qquad \hat{\rho}_j = \mathbf{a}_j^T R U Q^{(j)} \mathbf{a}_j, \qquad \hat{\gamma}_j = \mathbf{y}_j^T R U Q^{(j)} \mathbf{a}_j.$$

$$\omega_j = \mathbf{y}_j^T Q^{(j)} \mathbf{y}_j, \qquad \hat{\omega}_j = \mathbf{y}_j^T R U Q^{(j)} \mathbf{y}_j,$$

$$\beta_j = \mathbf{a}_j^T Q^{(j)} \mathbf{y}_j, \qquad \hat{\beta}_j = \mathbf{a}_j^T R U Q^{(j)} \mathbf{y}_j,$$

We note that $\alpha_j = 1 + \omega_j + 2\beta_j + \rho_j$.

REMARK 6.2 (Some qualitative facts).    We can roughly split the previous random variables into three classes:

– The vanishing terms: The terms $\beta_j$, $\hat{\beta}_j$ and $\hat{\gamma}_j$ vanish whenever $n$ goes to infinity; the main reason for this to hold true follows from $\mathbb{E}\mathbf{y}_j = 0$ (this is formally proved in Proposition B.1).

– The quadratic forms based on $\mathbf{y}_j$: The behavior of the terms $\omega_j$ and $\hat{\omega}_j$ can be deduced from the following well-known result:

$$x^T A x \sim \frac{1}{n} \operatorname{Tr} A, \qquad \text{where } x^T = \frac{1}{\sqrt{n}}(X_1, \ldots, X_n),$$

the $X_i$ being i.i.d., as long as $x$ and $A$ are independent [see (B.1) for the full statement]. Concerning $\hat{\omega}_j$, $R$ is not independent of $\mathbf{y}_j$ however we shall see in the course of the proof that $R$ can be replaced by $R^{(j)}$, which is independent of $\mathbf{y}_j$. Hence the previous result can also be applied to $\hat{\omega}_j$.

– The mysterious terms: Nothing is known about the asymptotic behavior of the terms $\rho_j$ and $\hat{\rho}_j$, except that these terms are bounded. Fortunately, these terms are always combined with vanishing terms, as we shall see in the sequel.

PROPOSITION 6.3.    *The following expression holds true:*

$$\frac{1}{N} \operatorname{Tr}(Q - R) U = Z_1^n + Z_2^n + Z_3^n + Z_4^n + Z_5^n,$$

*where*

$$Z_1^n = \frac{1}{N} \sum_{j=1}^n \frac{\hat{\gamma}_j}{\alpha_j},$$

$$Z_2^n = \frac{1}{N} \sum_{j=1}^n \frac{\hat{\beta}_j}{\alpha_j}(1 - z\tilde{b}_j(\rho_j + \beta_j)),$$

$$Z_3^n = \frac{1}{N} \sum_{j=1}^n \frac{\hat{\rho}_j}{\alpha_j}(z\tilde{b}_j \beta_j),$$

$$Z_4^n = \frac{1}{N} \sum_{j=1}^n \frac{\hat{\rho}_j}{\alpha_j}(1 + z\tilde{b}_j(1 + \omega_j)),$$

$$Z_5^n = \frac{1}{N} \sum_{j=1}^n \frac{\hat{\omega}_j}{\alpha_j} - \frac{1}{N} \operatorname{Tr}(B^{-1} + zI)RUQ.$$

PROOF.   Let us develop $(R - Q)U$:

$$(R(z) - Q(z))U$$
$$= Q(z)(Q^{-1}(z) - R^{-1}(z))R(z)U$$
$$= Q(z)(\Sigma\Sigma^T - zI - B^{-1}(z) + zA\tilde{B}(z)A^T)R(z)U$$
$$= Q(z)\Bigg\{ \sum_{j=1}^n ((1 + z\tilde{b}_j(z))\mathbf{a}_j\mathbf{a}_j^T + \mathbf{a}_j\mathbf{y}_j^T + \mathbf{y}_j\mathbf{a}_j^T + \mathbf{y}_j\mathbf{y}_j^T)$$
$$- (B^{-1}(z) + zI)\Bigg\}R(z)U.$$

In particular,

$$(6.5) \quad \frac{1}{N} \operatorname{Tr}(R(z) - Q(z))U = \frac{1}{N} \sum_{j=1}^n W_j - \frac{1}{N} \operatorname{Tr}(B^{-1} + zI)R(z)UQ(z),$$

where

$$W_j = (1 + z\tilde{b}_j)\mathbf{a}_j^T RUQ\mathbf{a}_j \;+\; \mathbf{y}_j^T RUQ\mathbf{a}_j \;+\; \mathbf{a}_j^T RUQ\mathbf{y}_j \;+\; \mathbf{y}_j^T RUQ\mathbf{y}_j.$$

Using (6.4) and the fact that $\alpha_j = 1 + \omega_j + 2\beta_j + \rho_j$, a straightforward (though lengthy) computation yields:

$$W_j = \frac{\hat{\gamma}_j}{\alpha_j} + \frac{\hat{\beta}_j}{\alpha_j}(1 - z\tilde{b}_j(\rho_j + \beta_j)) + \frac{\hat{\rho}_j}{\alpha_j}(1 + z\tilde{b}_j(1 + \omega_j + \beta_j)) + \frac{\hat{\omega}_j}{\alpha_j}.$$

Plugging this into (6.5), we obtain the desired result.   □

PROPOSITION 6.4.   *There exists a nonnegative constant $K$ such that*

$$\mathbb{E}|Z_1^n|^4 \le \frac{K}{n^2},$$

$$\mathbb{E}|Z_3^n|^{4+\varepsilon} \le \frac{K}{n^{2+\varepsilon/2}},$$

$$\mathbb{E}|Z_i^n|^{2+\varepsilon/2} \le \frac{K}{n^{1+\varepsilon/4}}, \qquad i = 2, 4, 5.$$

We first give a sketch of the proof, based on the qualitative facts given in Remark 6.2. The full proof is given in Appendix B.

REMARK 6.5 (Useful upper bounds). By Proposition 2.3 part 4(b), $b_i$, $\tilde{b}_\ell$, $b_i^{(j)}$ and $\tilde{b}_\ell^{(j)}$ belong to $\mathcal{S}(\mathbb{R}^+)$ with probability one. Therefore, as a consequence of Proposition 5.1,

$$\|R(z)\|_{\mathrm{sp}} \le \frac{1}{|\mathrm{Im}\, z|} \quad \text{and} \quad \|R^{(j)}(z)\|_{\mathrm{sp}} \le \frac{1}{|\mathrm{Im}\, z|}$$

with probability one. Furthermore,

$$(6.6) \qquad \tilde{q}_{jj}(z) = -\frac{1}{z(1 + \boldsymbol{\xi}_j^T Q^{(j)}(z) \boldsymbol{\xi}_j)} = -\frac{1}{z\alpha_j(z)},$$

where (6.6) can be derived by applying the results in Sections 0.7.3 and 0.7.4 in [15], for instance. Therefore, due to Proposition 2.3 part 4(a), we have

$$(6.7) \qquad \left| \frac{1}{\alpha_j(z)} \right| \le \left| \frac{z}{\mathrm{Im}\, z} \right|.$$

SKETCH OF PROOF OF PROPOSITION 6.4. We loosely explain why the random variables $Z_i^n$ ($1 \le i \le 5$) go to zero. We first look at $Z_1^n$. The term $\hat{\gamma}_j$ vanishes uniformly in the sense that one can prove that $\sup_{1 \le j \le n} \|\hat{\gamma}_j\|_4 \le \frac{K}{\sqrt{n}}$ (cf. Proposition B.1). Since $|\alpha_j|^{-1}$ is bounded (see Remark 6.5), Minkowski's inequality implies that $\|Z_1^n\|_4 = O(\frac{1}{\sqrt{n}})$.

The terms $Z_2^n$ and $Z_3^n$ can be handled in the same way: The terms $\hat{\beta}_j$ and $\beta_j$ vanish uniformly (cf. Proposition B.1). Furthermore, the terms $\rho_j$, $\hat{\rho}_j$, $\tilde{b}_j$ and $\alpha_j^{-1}$ are uniformly bounded. By consequence, $Z_2^n$ and $Z_3^n$ vanish for large $n$.

The analysis of $Z_4^n$ relies on arguments about quadratic forms. Denote by $\mathbf{x}_j = (X_{1,j}^n, \ldots, X_{N,j}^n)^T$. Then $\mathbf{y}_j = \frac{1}{\sqrt{n}} D_j^{1/2} \mathbf{x}_j$ and the quadratic form $\omega_j$ can be written:

$$\omega_j = \mathbf{y}_j^T Q^{(j)} \mathbf{y}_j = \frac{1}{n} \mathbf{x}_j^T D_j^{1/2} Q^{(j)} D_j^{1/2} \mathbf{x}_j.$$

As $\mathbf{x}_j$ and $D_j^{1/2} Q^{(j)} D_j^{1/2}$ are independent, and $D_j^{1/2} Q^{(j)} D_j^{1/2}$ is bounded, $\sup_{1 \le j \le n} \|\omega_j - \frac{1}{n} \mathrm{Tr}\, Q D_j\|_{2 + \varepsilon/2} = O(\frac{1}{\sqrt{n}})$ (see [1]). This implies that $z\tilde{b}_j(1 + \omega_j) \sim -1$ and that $\sup_{1 \le j \le n} \|1 + z\tilde{b}_j(1 + \omega_j)\|_{2 + \varepsilon/2} = O(\frac{1}{\sqrt{n}})$. The term $\alpha_j^{-1}\hat{\rho}_j$ being bounded, Minkowski's inequality immediately yields $\|Z_4^n\|_{2 + \varepsilon/2} = O(\frac{1}{\sqrt{n}})$. Let us now look at $Z_5^n$. We have $\hat{\omega}_j = \mathbf{y}_j^T RU Q^{(j)} \mathbf{y}_j$. Perturbation arguments yield:

$$\hat{\omega}_j = \mathbf{y}_j^T RU Q^{(j)} \mathbf{y}_j \sim \mathbf{y}_j^T R^{(j)} U Q^{(j)} \mathbf{y}_j \sim \frac{1}{n} \mathbf{x}_j^T D_j^{1/2} R^{(j)} U Q^{(j)} D_j^{1/2} \mathbf{x}_j,$$

where "$\sim$" loosely means "asymptotically equivalent." Now approximation formulas for quadratic forms yield:

$$\frac{1}{n}\mathbf{x}_j^T D_j^{1/2} R^{(j)} U Q^{(j)} D_j^{1/2} \mathbf{x}_j \sim \frac{1}{n}\operatorname{Tr} R^{(j)} U Q^{(j)} D_j \sim \frac{1}{n}\operatorname{Tr} RUQ D_j.$$

Since $(\alpha_j)^{-1} = -z\tilde{q}_{jj}$, we have:

$$\begin{equation}\frac{1}{N}\sum_{j=1}^n \frac{\tilde{\omega}_j}{\alpha_j} \sim -\frac{z}{N}\sum_{j=1}^n \tilde{q}_{jj}\frac{1}{n}\operatorname{Tr} RUQ D_j. \tag{6.8}\end{equation}$$

On the other hand, straightforward computation based on the mere definition of $b_i$ yields:

$$\begin{equation}\frac{1}{N}\operatorname{Tr}(B^{-1}(z)+zI_N)R(z)UQ(z) = -\frac{1}{n}\sum_{j=1}^n\left(z\tilde{q}_{jj}(z)\frac{1}{N}\operatorname{Tr} D_j R(z)UQ(z)\right) \tag{6.9}\end{equation}$$

[see (B.17) for details]. Combining (6.8) and (6.9) allows to show that $\|Z_5^n\|_{2+\varepsilon/2} = O(\frac{1}{\sqrt{n}})$.  □

PROOF OF LEMMA 6.1.  Since $\frac{1}{N}\operatorname{Tr}(Q(z)-R(z))U = Z_1^n + Z_2^n + Z_3^n + Z_4^n + Z_5^n$, Proposition 6.4 and Minkowski's inequality immediately yields

$$\begin{equation}\mathbb{E}\left|\frac{1}{N}\operatorname{Tr}(Q-R)U\right|^{2+\varepsilon/2} \leq \frac{K_2}{n^{1+\varepsilon/4}}. \tag{6.10}\end{equation}$$

Borel–Cantelli lemma implies that $\frac{1}{N}\operatorname{Tr}(Q(z)-R(z))U \to 0$ almost surely. The assertion $\frac{1}{n}\operatorname{Tr}(\tilde{Q}(z)-\tilde{R}(z))U \to 0$ a.s. can be proved similarly.  □

6.2. *Evaluation of the differences $\frac{1}{N}\operatorname{Tr} R(z) - \frac{1}{N}\operatorname{Tr} T(z)$ and $\frac{1}{n}\operatorname{Tr}\tilde{R}(z) - \frac{1}{n}\operatorname{Tr}\tilde{T}(z)$.* Recall that $d_n$ is defined by $d_n = \max(\frac{n}{N}, \frac{N}{n})$. As $\frac{N}{n} \to c$, $0 < c < +\infty$, it is clear that $\sup_n d_n < d$ where $d < +\infty$. Denote by $\mathcal{D}$ the subset of $\mathbb{C}^+$ defined by:

$$\mathcal{D} = \left\{z \in \mathbb{C}^+, \frac{|z|}{|\operatorname{Im} z|} < 2, |\operatorname{Im} z|^2 > 8\sigma_{\max}^2 \times \max(2d^2\sigma_{\max}^2, d\mathbf{a}_{\max}^2)\right\}.$$

LEMMA 6.6.  *Let $T$ and $\tilde{T}$ be given by Theorem 2.4. For every $z \in \mathcal{D}$, there exists constants $K_3$ and $\tilde{K}_3$ (possibly depending on $z$) such that*

$$\begin{equation}\begin{aligned}\mathbb{E}\left|\frac{1}{N}\operatorname{Tr} R(z) - \frac{1}{N}\operatorname{Tr} T(z)\right|^{2+\varepsilon/2} &\leq \frac{K_3}{n^{1+\varepsilon/4}} \quad and \\ \mathbb{E}\left|\frac{1}{n}\operatorname{Tr}\tilde{R}(z) - \frac{1}{n}\operatorname{Tr}\tilde{T}(z)\right|^{2+\varepsilon/2} &\leq \frac{\tilde{K}_3}{n^{1+\varepsilon/4}}.\end{aligned} \tag{6.11}\end{equation}$$

PROOF.    We only prove the first inequality in (6.11). We use the identity

$$R(z) - T(z) = R(z)(T(z)^{-1} - R(z)^{-1})T(z)$$
$$= R(z)(\Psi(z)^{-1} - B(z)^{-1})T(z) - zR(z)A(\tilde{\Psi}(z) - \tilde{B}(z))A^T T(z).$$

Therefore, $\frac{1}{N}\mathrm{Tr}(R(z) - T(z))$ can be written as

$$(6.12)\quad \frac{1}{N}\mathrm{Tr}(R(z) - T(z)) = \frac{1}{N}\mathrm{Tr}\,R(z)(\Psi(z)^{-1} - B^{-1}(z))T(z)$$
$$- \frac{z}{N}\mathrm{Tr}\,R(z)A(\tilde{\Psi}(z) - \tilde{B}(z))A^T T(z).$$

The first term of the right-hand side of (6.12) is equal to

$$\frac{1}{N}\mathrm{Tr}\,R(z)(\Psi(z)^{-1} - B^{-1}(z))T(z) = \frac{1}{N}\sum_{i=1}^{N}\left(\frac{1}{\psi_i(z)} - \frac{1}{b_i(z)}\right)(TR)_{ii}(z).$$

Let us first prove that the following control holds:

$$(6.13)\qquad\qquad |(TR)_{ii}(z)| \leq \frac{1}{|\mathrm{Im}\,z|^2}.$$

In fact, $|(TR)_{ii}(z)| \leq \|T_{i\cdot}(z)\|\|R_{\cdot i}\|$ where $T_{i\cdot}$ (resp. $R_{\cdot i}$) denotes row number $i$ of $T$ (resp. column number $i$ of $R$). On the other hand, $\|T_{i\cdot}(z)\| \leq \frac{1}{|\mathrm{Im}(z)|}$ and $\|R_{\cdot i}(z)\| \leq \frac{1}{|\mathrm{Im}(z)|}$ by Proposition 5.1 part 3. This proves (6.13). Minkowski's inequality leads to

$$\left\|\frac{1}{N}\mathrm{Tr}(R(z)(\Psi(z)^{-1} - B^{-1}(z))T(z))\right\|_{2+\varepsilon/2}$$
$$\leq \frac{1}{|\mathrm{Im}\,z|^2}\sup_i\left\|\frac{1}{\psi_i(z)} - \frac{1}{b_i(z)}\right\|_{2+\varepsilon/2}.$$

Similarly, the negative of the second term of the right-hand side of (6.12) can be written as

$$\frac{z}{N}\mathrm{Tr}\,R(z)A(\tilde{\Psi}(z) - \tilde{B}(z))A^T T(z) = \frac{z}{N}\sum_{j=1}^{n}(\tilde{\psi}_j(z) - \tilde{b}_j(z))\mathbf{a}_j^T T(z)R(z)\mathbf{a}_j.$$

Using again Minkowski's inequality and the fact that $|\mathbf{a}_j^T T(z)R(z)\mathbf{a}_j| \leq \frac{\mathbf{a}_{max}^2}{|\mathrm{Im}\,z|^2}$, we get:

$$\left\|\frac{z}{N}\mathrm{Tr}\,R(z)A(\tilde{\Psi}(z) - \tilde{B}(z))A^T T(z)\right\|_{2+\varepsilon/2}$$
$$\leq \frac{|z|d\mathbf{a}_{max,n}^2}{|\mathrm{Im}\,z|^2}\sup_j\|\tilde{\psi}_j(z) - \tilde{b}_j(z)\|_{2+\varepsilon/2}.$$

As $\tilde{\psi}_j(z) - \tilde{b}_j(z)$ is equal to $\tilde{\psi}_j(z)(\frac{1}{\tilde{b}_j(z)} - \frac{1}{\tilde{\psi}_j(z)})\tilde{b}_j(z)$, the inequality $\|\tilde{\psi}_j(z) \times \tilde{b}_j(z)\| \le \frac{1}{|\operatorname{Im} z|^2}$ yields

$$\|\tilde{\psi}_j(z) - \tilde{b}_j(z)\|_{2+\varepsilon/2} \le \frac{1}{|\operatorname{Im} z|^2}\left\|\frac{1}{\tilde{\psi}_j(z)} - \frac{1}{\tilde{b}_j(z)}\right\|_{2+\varepsilon/2}.$$

Gathering all the pieces together, we obtain

$$
\begin{aligned}
(6.14) \quad \left\|\frac{1}{N}\operatorname{Tr}(R(z) - T(z))\right\|_{2+\varepsilon/2} &\le \frac{1}{|\operatorname{Im} z|^2}\sup_i\left\|\frac{1}{\psi_i(z)} - \frac{1}{b_i(z)}\right\|_{2+\varepsilon/2} \\
&\quad + \frac{|z|d\mathbf{a}_{\max}^2}{|\operatorname{Im} z|^4}\sup_j\left\|\frac{1}{\tilde{\psi}_j(z)} - \frac{1}{\tilde{b}_j(z)}\right\|_{2+\varepsilon/2}.
\end{aligned}
$$

It is therefore sufficient to show that if $z \in \mathcal{D}$, then

$$(6.15) \qquad \sup_{1 \le k \le N}\left\|\frac{1}{\psi_k(z)} - \frac{1}{b_k(z)}\right\|_{2+\varepsilon/2} \le \frac{K_4}{\sqrt{n}},$$

$$(6.16) \qquad \sup_{1 \le \ell \le n}\left\|\frac{1}{\tilde{\psi}_\ell(z)} - \frac{1}{\tilde{b}_\ell(z)}\right\|_{2+\varepsilon/2} \le \frac{\tilde{K}_4}{\sqrt{n}},$$

where $K_4$ and $\tilde{K}_4$ do not depend on $n$ and $N$, but may depend on $z$.

For this, we use (6.1) in the case where matrix $\tilde{U}$ coincides with $\tilde{U} = \tilde{D}_i$ and constant $\tilde{K}_1 = \sigma_{\max}^2$ (recall that $\sup_{i,n}\|\tilde{D}_i\| < \sigma_{\max}^2$). $\frac{1}{N}\operatorname{Tr}(\tilde{D}_i\tilde{Q}(z))$ can be written as

$$\frac{1}{n}\operatorname{Tr}\tilde{D}_i\tilde{Q}(z) = \frac{1}{n}\operatorname{Tr}\tilde{D}_i\tilde{R}(z) + \tilde{\varepsilon}_i(z)$$

where $\sup_i\|\tilde{\varepsilon}_i(z)\|_{2+\varepsilon/2} \le \frac{\tilde{K}_2^\varepsilon}{\sqrt{n}}$, where $\tilde{K}_2^\varepsilon = \tilde{K}_2^{1/(2+\varepsilon/2)}$ (we denote similarly $K_2^\varepsilon = K_2^{1/(2+\varepsilon/2)}$). Using the very definitions of $\psi_i(z)$ and of $b_i(z)$, we obtain

$$\frac{1}{\psi_i(z)} - \frac{1}{b_i(z)} = \frac{z}{n}\operatorname{Tr}\tilde{D}_i(\tilde{R}(z) - \tilde{T}(z)) + z\tilde{\varepsilon}_i(z).$$

Rewriting $\tilde{R}(z) - \tilde{T}(z)$ as $\tilde{R}(z) - \tilde{T}(z) = \tilde{R}(z)(\tilde{T}(z)^{-1} - \tilde{R}(z)^{-1})\tilde{T}(z)$ and using similar arguments to what precedes [cf. (6.14)], we easily get

$$
\begin{aligned}
&\left\|\frac{1}{\psi_i(z)} - \frac{1}{b_i(z)}\right\|_{2+\varepsilon/2} \\
(6.17) \quad &\le \frac{|z|d\sigma_{\max}^2}{|\operatorname{Im} z|^2}\sup_\ell\left\|\frac{1}{\tilde{\psi}_\ell(z)} - \frac{1}{\tilde{b}_\ell(z)}\right\|_{2+\varepsilon/2} \\
&\quad + \frac{|z|^2\sigma_{\max}^2 d\mathbf{a}_{\max}^2}{|\operatorname{Im} z|^4}\sup_k\left\|\frac{1}{\psi_k(z)} - \frac{1}{b_k(z)}\right\|_{2+\varepsilon/2} + \frac{|z|\tilde{K}_2^\varepsilon}{\sqrt{n}}.
\end{aligned}
$$

Similarly, for each $1 \le j \le n$, we have

$$
(6.18) \quad
\begin{aligned}
&\left\| \frac{1}{\tilde{\psi}_j(z)} - \frac{1}{\tilde{b}_j(z)} \right\|_{2+\varepsilon/2} \\
&\le \frac{|z| d\sigma_{\max}^2}{|\operatorname{Im} z|^2} \sup_k \left\| \frac{1}{\psi_k(z)} - \frac{1}{b_k(z)} \right\|_{2+\varepsilon/2} \\
&\quad + \frac{|z|^2 \sigma_{\max}^2 d\mathbf{a}_{\max}^2}{|\operatorname{Im} z|^4} \sup_\ell \left\| \frac{1}{\tilde{\psi}_\ell(z)} - \frac{1}{\tilde{b}_\ell(z)} \right\|_{2+\varepsilon/2} + \frac{|z| K_2^\varepsilon}{\sqrt{n}}.
\end{aligned}
$$

In order to simplify the notation, we put

$$
\boldsymbol{\alpha}_n = \sup_{1 \le k \le N} \left\| \frac{1}{\psi_k(z)} - \frac{1}{b_k(z)} \right\|_{2+\varepsilon/2}, \qquad \tilde{\boldsymbol{\alpha}}_n = \sup_{1 \le \ell \le n} \left\| \frac{1}{\tilde{\psi}_\ell(z)} - \frac{1}{\tilde{b}_\ell(z)} \right\|_{2+\varepsilon/2}.
$$

Equations (6.17) and (6.18) give immediately

$$
(6.19) \quad
\begin{cases}
\boldsymbol{\alpha}_n \le \dfrac{|z| d\sigma_{\max}^2}{|\operatorname{Im} z|^2} \tilde{\boldsymbol{\alpha}}_n + \dfrac{|z|^2 d\sigma_{\max}^2 \mathbf{a}_{\max}^2}{|\operatorname{Im} z|^4} \boldsymbol{\alpha}_n + \dfrac{|z| \tilde{K}_2^\varepsilon}{\sqrt{n}}, \\[2mm]
\tilde{\boldsymbol{\alpha}}_n \le \dfrac{|z| d\sigma_{\max}^2}{|\operatorname{Im} z|^2} \boldsymbol{\alpha}_n + \dfrac{|z|^2 d\sigma_{\max}^2 \mathbf{a}_{\max}^2}{|\operatorname{Im} z|^4} \tilde{\boldsymbol{\alpha}}_n + \dfrac{|z| K_2^\varepsilon}{\sqrt{n}}.
\end{cases}
$$

As $z \in \mathcal{D}$, we have $\frac{|z|^2 d\sigma_{\max}^2 \mathbf{a}_{\max}^2}{|\operatorname{Im} z|^4} < \frac{1}{2}$. Therefore,

$$
\boldsymbol{\alpha}_n < \frac{2|z| d\sigma_{\max}^2}{|\operatorname{Im} z|^2} \tilde{\boldsymbol{\alpha}}_n + \frac{2|z| \tilde{K}_2^\varepsilon}{\sqrt{n}}.
$$

Plugging this inequality into (6.19), we obtain:

$$
(6.20) \quad \tilde{\boldsymbol{\alpha}}_n < \left( \frac{2|z|^2 d^2 \sigma_{\max}^4}{|\operatorname{Im} z|^4} + \frac{|z|^2 d\sigma_{\max}^2 \mathbf{a}_{\max}^2}{|\operatorname{Im} z|^4} \right) \tilde{\boldsymbol{\alpha}}_n + \frac{\tilde{K}_5}{\sqrt{n}}
$$

for some constant $\tilde{K}_5 > 0$ depending on $z$. It is then easy to check that

$$
\frac{2|z|^2 d^2 \sigma_{\max}^4}{|\operatorname{Im} z|^4} + \frac{|z|^2 d\sigma_{\max}^2 \mathbf{a}_{\max}^2}{|\operatorname{Im} z|^4} < 1
$$

for $z \in \mathcal{D}$. Thus (6.20) implies that for $z \in \mathcal{D}$, $\tilde{\boldsymbol{\alpha}}_n < \frac{\tilde{K}_6}{\sqrt{n}}$ for some constant $\tilde{K}_6$. Similarly, $\boldsymbol{\alpha}_n < \frac{K_6}{\sqrt{n}}$ and (6.15) and (6.16) are established. This, in turn, establishes (6.11). Proof of Lemma 6.6 is complete.  $\square$

6.3. *Proof of Theorem* 2.5.  We are now in position to complete the proof of Theorem 2.5. We first remark that inequality (6.1) for $U = I$ and inequality (6.11) imply

$$
(6.21) \quad \mathbb{E}\left| \frac{1}{N} \operatorname{Tr} Q(z) - \frac{1}{N} \operatorname{Tr} T(z) \right|^{2+\varepsilon/2} \le \frac{K_7}{n^{1+\varepsilon/4}}
$$

for each $z \in \mathcal{D}$. We consider a countable family $(z_k)_{k \in \mathbb{N}}$ with an accumulation point contained in a compact subset of $\mathcal{D}$. Borel–Cantelli lemma and equation (6.21) imply that on a set $\Omega_1$ with probability one, $\frac{1}{N} \mathrm{Tr}(Q(z_k) - T(z_k)) \to 0$ for each $k$. On $\Omega_1$, the functions $z \mapsto \frac{1}{N} \mathrm{Tr}\, Q(z)$ and $z \mapsto \frac{1}{N} \mathrm{Tr}\, T(z)$ belong to $\mathcal{S}(\mathbb{R}^+)$. We denote by $f_n(z)$ the function $f_n(z) = \frac{1}{N} \mathrm{Tr}(Q(z) - T(z))$. On $\Omega_1$, $f_n(z)$ is analytic on $\mathbb{C} - \mathbb{R}^+$. Moreover, $|f_n(z)| \leq \frac{2}{\delta_K}$ for each compact subset $K \subset \mathbb{C} - \mathbb{R}^+$, where $\delta_K$ represents the distance between $K$ and $\mathbb{R}^+$. By the normal family theorem (see, e.g., [21]), there exists a subsequence $f_{n_k}$ of $f_n$ which converges uniformly on each compact subset of $\mathbb{C} - \mathbb{R}^+$ to a function $f^*$ analytic on $\mathbb{C} - \mathbb{R}^+$. But, $f^*(z_k) = 0$ for each $k \in \mathbb{N}$. This implies that $f^*$ is identically 0 on $\mathbb{C} - \mathbb{R}^+$, and that the entire sequence $f_n$ uniformly converges to 0 on each compact subset of $\mathbb{C} - \mathbb{R}^+$. Therefore, almost surely, $\frac{1}{N} \mathrm{Tr}(Q(z) - T(z))$ converges to 0 for each $z \in \mathbb{C} - \mathbb{R}^+$. Proof of Theorem 2.5 is complete.

## APPENDIX A: PROOF OF PROPOSITION 2.2 PART 3

In this section, we prove Proposition 2.2 part 3. We first recall some results that concern the integral representations of some scalar complex functions analytic in the upper half plane $\mathbb{C}^+ = \{z : \mathrm{Im}(z) > 0\}$.

**A.1. The scalar case.** These results can be found in Krein and Nudelman's book [18] and therefore we adopt their notation in this section despite minor interferences with other notation in this article [beware in particular of the difference between class $\mathcal{S}$ functions below and Stieltjes transforms of probability measures denoted in the rest of the article by $\mathcal{S}(\mathbb{R}^+)$]:

THEOREM A.1 ([18], Theorem A.2).   *A function $f(z)$ over $\mathbb{C}^+$ satisfies the following conditions:*

   (i) $f(z)$ *is analytic over* $\mathbb{C}^+$ *and*
   (ii) $\mathrm{Im}\, f(z) \geq 0$ *for* $z \in \mathbb{C}^+$

*if and only if it admits the (unique) representation*

$$\text{(A.1)} \qquad f(z) = a + bz + \int_{-\infty}^{\infty} \left( \frac{1}{t-z} - \frac{t}{1+t^2} \right) \mu(dt),$$

*where $a \in \mathbb{R}$, $b \geq 0$ and $\mu$ is a positive measure such that*

$$\int_{-\infty}^{\infty} \frac{\mu(dt)}{1+t^2} < \infty.$$

The measure $\mu$ can be furthermore obtained by the Stieltjes inversion formula

$$\text{(A.2)} \quad \frac{1}{2}\mu(\{t_1\}) + \frac{1}{2}\mu(\{t_2\}) + \mu((t_1, t_2)) = \frac{1}{\pi} \lim_{\varepsilon \downarrow 0} \int_{t_1}^{t_2} \mathrm{Im}(f(t + i\varepsilon))\, dt.$$

Different names are given to the functions $f(z)$. In Krein and Nudelman's book [18], they are called class $\mathcal{R}$ functions. We will be particularly interested in the following subclass of these functions: A function is said to belong to class $\mathcal{S}$ if it belongs to class $\mathcal{R}$ and if it is furthermore analytic and nonnegative on the negative real axis $(-\infty, 0)$.

THEOREM A.2 ([18], Theorem A.4). *A function $f(z)$ is in class $\mathcal{S}$ if and only if it admits the representation*

$$(A.3) \qquad f(z) = c + \int_0^\infty \frac{\nu(dt)}{t - z},$$

*where $c \geq 0$ and $\nu$ is a positive measure that satisfies*

$$\int_0^\infty \frac{\nu(dt)}{1 + t} < \infty.$$

Class $\mathcal{S}$ functions can also be characterized by the following theorem:

THEOREM A.3 ([18], Theorem A.5). *A function $f(z)$ is in class $\mathcal{S}$ if both $f(z)$ and $zf(z)$ are in class $\mathcal{R}$.*

If $f(z)$ belongs to $\mathcal{S}$, it also admits the representation

$$(A.4) \qquad f(z) = a + bz + \int_{-\infty}^\infty \left( \frac{1}{t - z} - \frac{t}{1 + t^2} \right) \mu(dt).$$

In the following, it is useful to recall the relationships between representations (A.3) and (A.4). The intermediate steps of the proofs of [18], Theorem A.4 and [18], Theorem A.5 give:

PROPOSITION A.4. *The following relations hold:*

– $\mu$ *is carried by $\mathbb{R}^+$ and $\mu = \nu$,*
– $b = 0$,
– $a - \int_0^{+\infty} \frac{t}{1 + t^2}\, d\mu(t) \geq 0$ *and* $c = a - \int_0^{+\infty} \frac{t}{1 + t^2}\, d\mu(t) \geq 0$.

We now address a partial generalization of these results to matrix-valued functions.

**A.2. The matrix case.** A matrix-valued function $F(z)$ on $\mathbb{C}^+$ is said to belong to class $\mathcal{R}$ if $F(z)$ is analytic on $\mathbb{C}^+$ and if $\mathrm{Im}\, F(z) \geq 0$. Recall that matrix $\mathrm{Im}\, F(z)$ is defined as

$$\mathrm{Im}\, F(z) = \frac{1}{2\mathbf{i}} (F(z) - F^*(z)).$$

The generalization of Theorem A.1 can be found, for example, in [8]:

THEOREM A.5 ([8], Theorem 5.4).  *An $n \times n$ function $F(z)$ on $\mathbb{C}^+$ belongs to class $\mathcal{R}$ if and only if it admits the representation*

$$(A.5) \qquad F(z) = A + zB + \int_{-\infty}^{\infty} \left( \frac{1}{t-z} - \frac{t}{1+t^2} \right) \mu(dt),$$

*where $A$ is a Hermitian matrix, $B \geq 0$ and $\mu$ is a matrix-valued nonnegative measure such that*

$$\text{Tr} \int_{-\infty}^{\infty} \frac{\mu(dt)}{1+t^2} < \infty.$$

The proof is based on the corresponding result for the scalar case and the so-called polarization identity ([15], page 263)

$$x^* F(z) y = \tfrac{1}{4} ((x+y)^* F(z)(x+y) - (x-y)^* F(z)(x-y)$$
$$+ \mathbf{i}(x - \mathbf{i}y)^* F(z)(x - \mathbf{i}y) - \mathbf{i}(x + \mathbf{i}y)^* F(z)(x - \mathbf{i}y)).$$

We are now in position to prove Proposition 2.2 part 3. It partly generalizes Theorem A.2 to the matrix case:

THEOREM A.6.  *If a matrix function $F(z)$ satisfies:*
(i) *$F(z)$ and $zF(z)$ are in class $\mathcal{R}$*
    *then, it admits the representation*

$$(A.6) \qquad\qquad F(z) = C + \int_0^{\infty} \frac{\mu(dt)}{t-z},$$

*where $C \geq 0$ and $\mu$ is a matrix-valued nonnegative measure carried by $\mathbb{R}^+$ that satisfies*

$$\text{Tr} \left( \int_0^{\infty} \frac{\mu(dt)}{1+t} \right) < \infty.$$

PROOF.  Theorem A.6 could again be proved using the polarization identity. We however provide the following shorter argument. Assume that $F(z)$ and $zF(z)$ are in class $\mathcal{R}$. Then, $F(z)$ can be written as (A.5). Let $x \in \mathbb{C}^n$. Then, $x^* F(z) x$ is scalar function which belongs to $\mathcal{R}$. The representation given by (A.1) in Theorem A.1 is given by

$$x^* F(z) x = x^* A x + z x^* B x + \int_{-\infty}^{\infty} \left( \frac{1}{t-z} - \frac{t}{1+t^2} \right) x^* \mu x(dt).$$

The quantity $x^* z F(z) x = z x^* F(z) x$ also belongs to belongs to $\mathcal{R}$. Therefore, it can be written as

$$x^* F(z) x = c_x + \int_0^{\infty} \frac{\nu_x(dt)}{t-z},$$

where $c_x \geq 0$ and where $\nu_x$ is a positive measure for which

$$\int_0^\infty \frac{\nu_x(dt)}{1+t} < \infty.$$

Using Proposition A.4, we get immediately that:

– $x^*\mu x$ is carried by $\mathbb{R}^+$ and $\nu_x = x^*\mu x$,
– $x^*Bx = 0$,
– $c_x = x^*Ax - \int_0^\infty \frac{t}{1+t^2}(x^*\mu x)(dt) \geq 0$.

The first item implies that $\mu$ is carried by $\mathbb{R}^+$ and that $\mathrm{Tr}(\int_0^\infty \frac{\mu(dt)}{1+t}) < \infty$. The second item implies that $B = 0$. As $t(1+t^2)^{-1} \leq 2(1+t)^{-1}$ for $t \geq 0$, the finiteness of $\mathrm{Tr}(\int_0^\infty (1+t)^{-1}\mu(dt))$ implies

$$\mathrm{Tr}\left(\int_0^\infty \frac{t}{1+t^2}\mu(dt)\right) < +\infty.$$

Therefore $F(z)$ can be written as

$$F(z) = A - \int_0^\infty \frac{t}{1+t^2}\mu(dt) + \int_0^\infty \frac{\mu(dt)}{t-z}.$$

The third item implies that matrix $C = A - \int_0^\infty \frac{t}{1+t^2}\mu(dt)$ is nonnegative. This completes the proof. $\quad\square$

We finally note that Theorem A.6 can be found in several papers without proof. See, for example, [3], pages 64–65.

## APPENDIX B: PROOF OF PROPOSITION 6.4

In the sequel, $K$ will denote a bounding constant that does not depend on $N$. Its value might change from one inequality to another.

We first recall a useful result on quadratic forms associated with random matrices.

A LEMMA FROM BAI AND SILVERSTEIN (Lemma 2.7 in [1]).    *Let* $\mathbf{x} = (X_1, \ldots, X_n)^T$ *be a vector where the* $X_i$'s *are centered i.i.d. random variables with unit variance. Let* $C$ *be an* $n \times n$ *deterministic complex matrix. Then, for any* $p \geq 2$,

$$(\mathrm{B.1}) \quad \mathbb{E}|\mathbf{x}^T C \mathbf{x} - \mathrm{Tr}\, C|^p \leq K_p((\mathbb{E}|X_1|^4 \,\mathrm{Tr}\, CC^*)^{p/2} + \mathbb{E}|X_1|^{2p} \,\mathrm{Tr}(CC^*)^{p/2}).$$

### B.1. Evaluation of $Z_1^n$, $Z_2^n$ and $Z_3^n$.

PROPOSITION B.1.    *The random variables* $\beta_j$, $\hat{\gamma}_j$ *and* $\hat{\beta}_j$ *satisfy*

$$\sup_{1 \leq j \leq n} \mathbb{E}|\beta_j|^{4+\varepsilon} \leq \frac{K}{n^{2+\varepsilon/2}}, \qquad \sup_{1 \leq j \leq n} \mathbb{E}|\hat{\gamma}_j|^4 \leq \frac{K}{n^2}, \qquad \sup_{1 \leq j \leq n} \mathbb{E}|\hat{\beta}_j|^4 \leq \frac{K}{n^2},$$

*where* $K$ *is a constant that does not depend on* $n$.

PROOF. Let us first establish the inequality for $\beta_j$. We recall that $\mathbf{y}_j$ can be written as $\mathbf{y}_j = \frac{1}{\sqrt{n}} D_j^{1/2} \mathbf{x}_j$ where we recall that $\mathbf{x}_j = (X_{1,j}^n, \ldots, X_{N,j}^n)^T$. Let $\mathbf{v} = D_j^{1/2} Q^{(j)} \mathbf{a}_j$. We can thus write $\beta_j = \frac{1}{\sqrt{n}} \mathbf{v}^T \mathbf{x}_j$. We then have

$$
\begin{aligned}
\mathbb{E}[&|\mathbf{a}_j^T Q^{(j)} \mathbf{y}_j|^{4+\varepsilon} | \mathbf{v}] \\
&= \frac{1}{n^{2+\varepsilon/2}} \mathbb{E}[|\mathbf{x}_j^T \mathbf{v} \mathbf{v}^T \mathbf{x}_j|^{2+\varepsilon/2} | \mathbf{v}] \\
&\leq \frac{2^{2+\varepsilon/2}}{n^{2+\varepsilon/2}} (\mathbb{E}[|\mathbf{x}_j^T \mathbf{v} \mathbf{v}^T \mathbf{x}_j - \mathrm{Tr}\,\mathbf{v}\mathbf{v}^T|^{2+\varepsilon/2} | \mathbf{v}] + |\mathrm{Tr}\,\mathbf{v}\mathbf{v}^T|^{2+\varepsilon/2}) \\
&\overset{(a)}{\leq} \frac{K}{n^{2+\varepsilon/2}} ((\mathbb{E}X_{11}^4 \,\mathrm{Tr}(\mathbf{v}\mathbf{v}^T(\mathbf{v}\mathbf{v}^T)^*))^{1+\varepsilon/4} \\
&\qquad\qquad + \mathbb{E}|X_{11}|^{4+\varepsilon} \,\mathrm{Tr}(\mathbf{v}\mathbf{v}^T(\mathbf{v}\mathbf{v}^T)^*)^{1+\varepsilon/4} + \|\mathbf{v}\|^{4+\varepsilon}) \\
&\leq \frac{K\|\mathbf{v}\|^{4+\varepsilon}}{n^{2+\varepsilon/2}} ((\mathbb{E}X_{11}^4)^{1+\varepsilon/4} + \mathbb{E}|X_{11}|^{4+\varepsilon} + 1),
\end{aligned}
$$
(B.2)

where (a) follows from the independence of $\mathbf{v}$ and $\mathbf{x}_j$ and from Lemma 2.7 in [1]. By noticing that $\|\mathbf{v}\| \leq \|D_j^{1/2}\|_{\mathrm{sp}} \|Q^{(j)}\|_{\mathrm{sp}} \|\mathbf{a}_j\| \leq \frac{\mathbf{a}_{\max}\sigma_{\max}}{|\mathrm{Im}(z)|}$ and by using Assumption A-1, we have the desired result.

Let us now establish the inequality for $\hat{\gamma}_j$. The random variable $\hat{\gamma}_j$ can be written

$$
\begin{aligned}
\hat{\gamma}_j &= \mathbf{y}_j^T R^{(j)} U Q^{(j)} \mathbf{a}_j + \mathbf{y}_j^T (R - R^{(j)}) U Q^{(j)} \mathbf{a}_j \\
&= \mathbf{y}_j^T R^{(j)} U Q^{(j)} \mathbf{a}_j + \mathbf{y}_j^T R^{(j)} (R^{(j)^{-1}} - R^{-1}) R U Q^{(j)} \mathbf{a}_j \\
&= \mathbf{y}_j^T R^{(j)} U Q^{(j)} \mathbf{a}_j \\
&\qquad + \mathbf{y}_j^T R^{(j)} (B^{(j)^{-1}} - B^{-1} - z(A^{(j)} \tilde{B}^{(j)} A^{(j)^T} - A \tilde{B} A^T)) R U Q^{(j)} \mathbf{a}_j \\
&= \hat{\gamma}_{j,1} + \hat{\gamma}_{j,2} - \hat{\gamma}_{j,3} + \hat{\gamma}_{j,4},
\end{aligned}
$$

where

$$
\begin{aligned}
\hat{\gamma}_{j,1} &= \mathbf{y}_j^T R^{(j)} U Q^{(j)} \mathbf{a}_j, \\
\hat{\gamma}_{j,2} &= \mathbf{y}_j^T R^{(j)} (B^{(j)^{-1}} - B^{-1}) R U Q^{(j)} \mathbf{a}_j, \\
\hat{\gamma}_{j,3} &= z \sum_{\substack{\ell=1 \\ \ell \neq j}}^{n} (\tilde{b}_\ell^{(j)} - \tilde{b}_\ell) \mathbf{y}_j^T R^{(j)} \mathbf{a}_\ell \mathbf{a}_\ell^T R U Q^{(j)} \mathbf{a}_j, \\
\hat{\gamma}_{j,4} &= z \tilde{b}_j \mathbf{y}_j^T R^{(j)} \mathbf{a}_j \mathbf{a}_j^T R U Q^{(j)} \mathbf{a}_j.
\end{aligned}
$$
(B.3)

Beginning with the term $\hat{\gamma}_{j,1}$, let $\mathbf{v} = R^{(j)} U Q^{(j)} \mathbf{a}_j$. Recalling that $\mathbf{y}_j = \frac{1}{\sqrt{n}} D_j^{1/2} \times \mathbf{x}_j$ and using the independence of $\mathbf{v}$ and $\mathbf{y}_j$, a standard calculation leads to

$$\mathbb{E}|\hat{\gamma}_{j,1}|^4 \leq \frac{\sigma_{\max}^4}{n^2}(\mathbb{E} X_{11}^4 + 3)\mathbb{E}\|\mathbf{v}\|^4,$$

where $\|\mathbf{v}\|^4 \leq (\|R^{(j)}\|_{\mathrm{sp}}\|U\|_{\mathrm{sp}}\|Q^{(j)}\|_{\mathrm{sp}}\|\mathbf{a}_j\|)^4 \leq \frac{K_1^4 \mathbf{a}_{\max}^4}{(\mathrm{Im}\, z)^8}$. This yields

$$(B.4) \qquad \mathbb{E}|\hat{\gamma}_{j,1}|^4 \leq \frac{\sigma_{\max}^4 K_1^4 \mathbf{a}_{\max}^4}{n^2(\mathrm{Im}\, z)^8}(\mathbb{E} X_{11}^4 + 3) \propto \frac{1}{n^2}.$$

Let us now consider $\hat{\gamma}_{j,2} = \mathbf{y}_j^T R^{(j)}(B^{(j)^{-1}} - B^{-1})RUQ^{(j)}\mathbf{a}_j$. We have

$$(B.5) \qquad \begin{aligned} |\hat{\gamma}_{j,2}| &\leq \|R^{(j)}\|_{\mathrm{sp}}\|B^{(j)^{-1}} - B^{-1}\|_{\mathrm{sp}}\|R\|_{\mathrm{sp}}\|U\|_{\mathrm{sp}}\|Q^{(j)}\|_{\mathrm{sp}}\|\mathbf{a}_j\|\|\mathbf{y}_j\|, \\ &\leq \frac{K_1 \mathbf{a}_{\max}}{(\mathrm{Im}\, z)^3}\|B^{(j)^{-1}} - B^{-1}\|_{\mathrm{sp}}\|\mathbf{y}_j\|. \end{aligned}$$

We now prove that

$$(B.6) \qquad \|B^{(j)^{-1}} - B^{-1}\|_{\mathrm{sp}} \leq \frac{2\sigma_{\max}^2}{n}\left|\frac{z}{\mathrm{Im}\, z}\right|.$$

By applying to $(\Sigma^T \Sigma - zI)^{-1}$ the matrix inversion lemma (see [16], Appendix A, see also [15], Section 0.7.4), we obtain

$$\frac{1}{n}\mathrm{Tr}\,\tilde{D}_i(\Sigma^T\Sigma - zI)^{-1} = x_{j,1} + x_{j,2} + x_{j,3},$$

where

$$\begin{aligned} x_{j,1} &= \frac{1}{n}\mathrm{Tr}\,\tilde{D}_i^{(j)}(\Sigma^{(j)^T}\Sigma^{(j)} - zI)^{-1}, \\ x_{j,2} &= \frac{1}{n}\frac{\mathrm{Tr}\,\tilde{D}_i^{(j)}(\Sigma^{(j)^T}\Sigma^{(j)} - zI)^{-1}\Sigma^{(j)^T}\boldsymbol{\xi}_j\boldsymbol{\xi}_j^T\Sigma^{(j)}(\Sigma^{(j)^T}\Sigma^{(j)} - zI)^{-1}}{-z - z\boldsymbol{\xi}_j^T(\Sigma^{(j)}\Sigma^{(j)^T} - zI)^{-1}\boldsymbol{\xi}_j}, \\ x_{j,3} &= \frac{1}{n}\frac{\sigma_{ij}^2}{-z - z\boldsymbol{\xi}_j^T(\Sigma^{(j)}\Sigma^{(j)^T} - zI)^{-1}\boldsymbol{\xi}_j}. \end{aligned}$$

Definitions (2.7) and (6.2) yield

$$\begin{aligned} \frac{1}{b_i^{(j)}(z)} - \frac{1}{b_i(z)} &= -\frac{z}{n}(\mathrm{Tr}(\tilde{D}_i^{(j)}(\Sigma^{(j)^T}\Sigma^{(j)} - zI)^{-1}) - \mathrm{Tr}(\tilde{D}_i(\Sigma^T\Sigma - zI)^{-1})), \\ &= z(x_{j,2} + x_{j,3}). \end{aligned}$$

We have

$$|x_{j,2}| = \frac{1}{n}\left|\frac{\boldsymbol{\xi}_j^T \Sigma^{(j)}(\Sigma^{(j)T}\Sigma^{(j)} - zI)^{-1}\tilde{D}_i^{(j)}(\Sigma^{(j)T}\Sigma^{(j)} - zI)^{-1}\Sigma^{(j)T}\boldsymbol{\xi}_j}{-z - z\boldsymbol{\xi}_j^T(\Sigma^{(j)}\Sigma^{(j)T} - zI)^{-1}\boldsymbol{\xi}_j}\right|,$$

(B.7)

$$\leq \frac{\|\tilde{D}_i^{(j)}\|_{\mathrm{sp}}^2}{n}\frac{\|(\Sigma^{(j)T}\Sigma^{(j)} - zI)^{-1}\Sigma^{(j)T}\boldsymbol{\xi}_j\|^2}{|z + z\boldsymbol{\xi}_j^T(\Sigma^{(j)}\Sigma^{(j)T} - zI)^{-1}\boldsymbol{\xi}_j|} \overset{\text{(a)}}{\leq} \frac{\sigma_{\max}^2}{n\,\mathrm{Im}\,z},$$

where (a) follows from a singular value decomposition of $\Sigma^{(j)}$. In fact, let $\Sigma^{(j)} = \sum_{\ell=1}^N \nu_\ell \mathbf{u}_\ell \mathbf{v}_\ell^T$ be a singular value decomposition of $\Sigma^{(j)}$ with $\nu_\ell$, $\mathbf{u}_\ell$, and $\mathbf{v}_\ell$ as singular value, left singular vector, and right singular vector respectively. Then,

$$\|(\Sigma^{(j)T}\Sigma^{(j)} - zI)^{-1}\Sigma^{(j)T}\boldsymbol{\xi}_j\|^2 = \sum_{\ell=1}^N \frac{\nu_\ell^2|\mathbf{u}_\ell^T\boldsymbol{\xi}_j|^2}{|\nu_\ell^2 - z|^2},$$

$$\mathrm{Im}(z\boldsymbol{\xi}_j^T(\Sigma^{(j)}\Sigma^{(j)T} - zI)^{-1}\boldsymbol{\xi}_j) = \mathrm{Im}(z)\left(\sum_{\ell=1}^N \frac{\nu_\ell^2|\mathbf{u}_\ell^T\boldsymbol{\xi}_j|^2}{|\nu_\ell^2 - z|^2}\right)$$

which yields (B.7). Furthermore, one has $|x_{j,3}| \leq \frac{\sigma_{\max}^2}{n}\frac{1}{\mathrm{Im}(z)}$. Thus,

(B.8)
$$\left|\frac{1}{b_i^{(j)}(z)} - \frac{1}{b_i(z)}\right| \leq \frac{2\sigma_{\max}^2}{n}\left|\frac{z}{\mathrm{Im}\,z}\right|,$$

which yields (B.6). Plugging this into (B.5), we obtain

$$|\hat{\gamma}_{j,2}| \leq \frac{2\sigma_{\max}^2 K_1 \mathbf{a}_{\max}|z|}{n|\mathrm{Im}\,z|^4}\|\mathbf{y}_j\|.$$

Finally, since $\mathbb{E}\|\mathbf{y}_j\|^4 \leq \sigma_{\max}^4 d^2(\mathbb{E}|X_{11}|^4 + 1)$, we have

(B.9)
$$\mathbb{E}|\hat{\gamma}_{j,2}|^4 \leq \frac{16d^8\sigma_{\max}^{12}K_1^4\mathbf{a}_{\max}^4|z|^4}{|\mathrm{Im}\,z|^{16}n^4}(1 + \mathbb{E}|X_{11}^n|^4) \propto \frac{1}{n^4}$$

for $N$ large enough.

We now deal with

$$\hat{\gamma}_{j,3} = z\sum_{\substack{\ell=1\\ \ell\neq j}}^n (\tilde{b}_\ell^{(j)} - \tilde{b}_\ell)\mathbf{y}_j^T R^{(j)}\mathbf{a}_\ell \mathbf{a}_\ell^T RUQ^{(j)}\mathbf{a}_j.$$

A rough control yields

$$|\hat{\gamma}_{j,3}| \leq |z|\mathbf{a}_{\max}^2\|R\|_{\mathrm{sp}}\|U\|_{\mathrm{sp}}\|Q^{(j)}\|_{\mathrm{sp}}\sum_{\substack{\ell=1\\ \ell\neq j}}^n |\tilde{b}_\ell^{(j)} - \tilde{b}_\ell||\mathbf{y}_j^T R^{(j)}\mathbf{a}_\ell|,$$

$$\leq \frac{|z|\mathbf{a}_{\max}^2 K_1}{|\operatorname{Im} z|^2} \sum_{\substack{\ell=1 \\ \ell \neq j}}^{n} |\tilde{b}_\ell^{(j)} - \tilde{b}_\ell| |\mathbf{y}_j^T R^{(j)} \mathbf{a}_\ell|.$$

We have

$$\begin{aligned}
\tilde{b}_\ell^{(j)}(z) &- \tilde{b}_\ell(z) \\
&= \tilde{b}_\ell^{(j)}(z)\tilde{b}_\ell(z)(\tilde{b}_\ell^{-1}(z) - \tilde{b}_\ell^{(j)-1}(z)) \\
&= z\tilde{b}_\ell^{(j)}(z)\tilde{b}_\ell(z)\frac{1}{n}\operatorname{Tr}(D_\ell((\Sigma\Sigma^T - zI)^{-1} - (\Sigma^{(j)}\Sigma^{(j)T} - zI)^{-1})).
\end{aligned}$$

Since $\tilde{b}_\ell$ and $\tilde{b}_\ell^{(j)}$ belong to $\mathcal{S}(\mathbb{R}^+)$ by Proposition 2.3 part 4(b), their absolute values are bounded above by $|\operatorname{Im} z|^{-1}$. By Lemma 2.6 in [25], we have

$$(B.10) \qquad |\operatorname{Tr} D_\ell((\Sigma\Sigma^T - zI)^{-1} - (\Sigma^{(j)}\Sigma^{(j)T} - zI)^{-1})| \leq \frac{\sigma_{\max}^2}{|\operatorname{Im} z|}.$$

As a result, we obtain

$$(B.11) \qquad |\tilde{b}_\ell^{(j)}(z) - \tilde{b}_\ell(z)| \leq \frac{\sigma_{\max}^2 |z|}{n|\operatorname{Im} z|^3}.$$

Thus,

$$|\hat{\gamma}_{j,3}| \leq \frac{|z|^2\sigma_{\max}^2 \mathbf{a}_{\max}^2 K_1}{n|\operatorname{Im} z|^5} \sum_{\substack{\ell=1 \\ \ell \neq j}}^{n} |\mathbf{y}_j^T R^{(j)} \mathbf{a}_\ell|.$$

By Minkowski's inequality,

$$\mathbb{E}|\hat{\gamma}_{j,3}|^4 \leq \frac{|z|^8\sigma_{\max}^8 \mathbf{a}_{\max}^8 K_1^4}{n^4|\operatorname{Im} z|^{20}} \left( \sum_{\substack{\ell=1 \\ \ell \neq j}}^{n} (\mathbb{E}|\mathbf{y}_j^T R^{(j)} \mathbf{a}_\ell|^4)^{1/4} \right)^4.$$

The terms $\mathbb{E}|\mathbf{y}_j^T R^{(j)} \mathbf{a}_\ell|^4$ can be handled in the same way as for $\hat{\gamma}_{j,1}$. Thus, we obtain

$$(B.12) \qquad \mathbb{E}|\hat{\gamma}_{j,3}|^4 \leq \frac{d^4|z|^8\sigma_{\max}^{12} \mathbf{a}_{\max}^{12} K_1^4}{n^2|\operatorname{Im} z|^{24}} (\mathbb{E}|X_{11}^n|^4 + 3) \propto \frac{1}{n^2}.$$

The last term $\hat{\gamma}_{j,4} = z\tilde{b}_j \mathbf{y}_j^T R^{(j)} \mathbf{a}_j \mathbf{a}_j^T RUQ^{(j)} \mathbf{a}_j$ satisfies

$$|\hat{\gamma}_{j,4}| \leq |z\tilde{b}_j| |\mathbf{a}_j^T RUQ^{(j)} \mathbf{a}_j| |\mathbf{y}_j^T R^{(j)} \mathbf{a}_j|$$

Hence

$$\begin{aligned}
(B.13) \qquad \mathbb{E}|\hat{\gamma}_{j,4}|^4 &\leq \left| \frac{z^4}{(\operatorname{Im} z)^{12}} \right| \mathbf{a}_{\max}^8 K_1^4 \mathbb{E}|\mathbf{y}_j^T R^{(j)} \mathbf{a}_j|^4, \\
&\overset{(a)}{\leq} \left| \frac{z^4}{(\operatorname{Im} z)^{16}} \right| \mathbf{a}_{\max}^{12} K_1^4 \sigma_{\max}^4 (\mathbb{E}|X_{11}^n|^4 + 3)\frac{1}{n^2} \propto \frac{1}{n^2},
\end{aligned}$$

where (a) follows from the fact that the term $\mathbb{E}|\mathbf{y}_j^T R^{(j)} \mathbf{a}_j|^4$ can be handled as $\hat{\gamma}_{j,1}$. Gathering (B.4), (B.9), (B.12), (B.13) and using Minkowski's inequality, we obtain the desired inequality for $\hat{\gamma}_j$.

Using similar arguments, one can prove the same inequality for $\hat{\beta}_j$. Proposition B.1 is proved. □

*The term $Z_1^n$.* Using (6.7), we have
$$Z_1^n = \frac{1}{N} \sum_{j=1}^n \frac{\hat{\gamma}_j}{\alpha_j} \leq \frac{|z|}{N|\operatorname{Im} z|} \sum_{j=1}^n |\hat{\gamma}_j|.$$

Minkowski's inequality and Proposition B.1 yield
$$\mathbb{E}|Z_1^n|^4 \leq \left( \frac{n|z|}{N|\operatorname{Im} z|} \right)^4 \frac{K}{n^2} \propto \frac{1}{n^2}.$$

*The term $Z_2^n$.* We can write $Z_2^n = \frac{1}{N} \sum_{j=1}^n Z_{j,2,1} - \frac{1}{N} \sum_{j=1}^n Z_{j,2,2}$, where
$$Z_{j,2,1} = (1 - z\tilde{b}_j(z)\rho_j)\hat{\beta}_j/\alpha_j \quad \text{and} \quad Z_{j,2,2} = z\tilde{b}_j(z)\hat{\beta}_j\beta_j/\alpha_j.$$

We have $|\rho_j| \leq \|Q^{(j)}(z)\|_{\mathrm{sp}} \|\mathbf{a}_j\|^2 \leq \mathbf{a}_{\max}^2/|\operatorname{Im} z|$. As a consequence,
$$|Z_{j,2,1}| \leq \left| \frac{z}{\operatorname{Im} z} \right| \left( 1 + \mathbf{a}_{\max}^2 \left| \frac{z}{(\operatorname{Im} z)^2} \right| \right) |\hat{\beta}_j|.$$

Proposition B.1 yields: $\mathbb{E}|\hat{\beta}_j|^4 \leq K/n^2$. Therefore, $\mathbb{E}|Z_{j,2,1}|^4 \leq K/n^2$ and Minkowski's inequality implies that $\mathbb{E}|\frac{1}{N} \sum_{j=1}^n Z_{j,2,1}|^4 \leq K/n^2$. In particular, $\mathbb{E}|\frac{1}{N} \sum_{j=1}^n Z_{j,2,1}|^{2+\varepsilon/2} \leq K/n^{1+\varepsilon/4}$ for $\varepsilon \leq 4$ by Lyapunov's inequality.

Let us consider $Z_{j,2,2}$. By the Cauchy–Schwarz inequality,
$$\mathbb{E}|Z_{j,2,2}|^{2+\varepsilon/2} \leq \left| \frac{z}{\operatorname{Im} z} \right|^{4+\varepsilon} (\mathbb{E}|\hat{\beta}_j|^{4+\varepsilon})^{1/2} (\mathbb{E}|\beta_j|^{4+\varepsilon})^{1/2}.$$

We have
$$\mathbb{E}|\hat{\beta}_j|^{4+\varepsilon} \leq \mathbb{E}((\|\mathbf{a}_j\| \|R\|_{\mathrm{sp}} \|U\|_{\mathrm{sp}} \|Q^{(j)}\|_{\mathrm{sp}})^{4+\varepsilon} \|\mathbf{y}_j\|^{4+\varepsilon})$$
$$\leq \left( \frac{\mathbf{a}_{\max} K_1}{(\operatorname{Im} z)^2} \right)^{4+\varepsilon} \mathbb{E}\|\mathbf{y}_j\|^{4+\varepsilon}.$$

By Minkowski's inequality, we have
$$\begin{aligned}
\mathbb{E}\|\mathbf{y}_j\|^{4+\varepsilon} &\leq \frac{\sigma_{\max}^{4+\varepsilon}}{n^{2+\varepsilon/2}} \mathbb{E}\left( \sum_{m=1}^N X_{m,j}^2 \right)^{2+\varepsilon/2} \\
&\leq \frac{\sigma_{\max}^{4+\varepsilon}}{n^{2+\varepsilon/2}} (N^{2+\varepsilon/2} \mathbb{E}|X_{1,1}|^{4+\varepsilon}) \leq K.
\end{aligned}$$
(B.14)

Therefore, $\mathbb{E}|\hat{\beta}_j|^{4+\varepsilon} < K$.

As $\mathbb{E}|\beta_j|^{4+\varepsilon} \le K/n^{2+\varepsilon/2}$ by Proposition B.1, we have $\mathbb{E}|\frac{1}{N}\sum_{j=1}^n Z_{j,2,2}|^{2+\varepsilon/2} \le K/n^{1+\varepsilon/4}$, which yields the desired result.

*The term $Z_3^n$.* We have

$$|Z_3^n| = \left| \frac{1}{N} \sum_{j=1}^n \frac{\hat{\rho}_j}{\alpha_j} (z\tilde{b}_j\beta_j) \right| \le K_1 \frac{\mathbf{a}_{\max}^2 |z|^2}{N|\operatorname{Im} z|^3} \sum_{j=1}^n |\beta_j|.$$

Therefore, Proposition B.1 yields the desired result.

**B.2. Evaluation of $Z_4^n$.** We rely here on Lemma 2.7 in [1] on the quadratic forms recalled at the beginning of this Appendix. Write $Z_4^n = \frac{1}{N}\sum_{j=1}^n Z_{j,4}$ with

$$Z_{j,4} = (1 + z\tilde{b}_j(z)(1+\omega_j))\hat{\rho}_j/\alpha_j.$$

We write $\omega_j = \frac{1}{n}\operatorname{Tr} D_j Q(z) + \varepsilon_{j,1} + \varepsilon_{j,2}$, where

$$\varepsilon_{j,1} = \omega_j - \frac{1}{n}\operatorname{Tr} D_j Q^{(j)}(z) \quad \text{and} \quad \varepsilon_{j,2} = \frac{1}{n}\operatorname{Tr} D_j(Q^{(j)}(z) - Q(z)).$$

Since $\omega_j$ is equal to $\frac{1}{n}\mathbf{x}_j^T D_j^{1/2} Q^{(j)}(z) D_j^{1/2}\mathbf{x}_j$, where $\mathbf{x}_j = (X_{1,j}^n, \ldots, X_{N,j}^n)^T$, (B.1) yields

$$\begin{aligned}
\mathbb{E}|\varepsilon_{j,1}|^{2+\varepsilon/2} &\le \frac{K}{n^{2+\varepsilon/2}}[(\mathbb{E}|X_{11}^n|^4 \operatorname{Tr} D_j^2 Q^{(j)} Q^{(j)*})^{1+\varepsilon/4} \\
&\quad + \mathbb{E}|X_{11}^n|^{4+\varepsilon} \operatorname{Tr}(D_j^2 Q^{(j)} Q^{(j)*})^{1+\varepsilon/4}], \\
&\le \frac{K}{n^{2+\varepsilon/2}}\left[\left(\mathbb{E}|X_{11}^n|^4 \frac{N\sigma_{\max}^4}{(\operatorname{Im} z)^2}\right)^{1+\varepsilon/4} + \mathbb{E}|X_{11}^n|^{4+\varepsilon} \frac{N\sigma_{\max}^{4+\varepsilon}}{|\operatorname{Im} z|^{2+\varepsilon/2}}\right].
\end{aligned}$$

Thanks to Assumption A-1, we therefore have

$$\text{(B.15)} \qquad \mathbb{E}|\varepsilon_{j,1}|^{2+\varepsilon/2} \le \frac{K}{n^{1+\varepsilon/4}}.$$

By Lemma 2.6 in [25], we have

$$\text{(B.16)} \qquad |\varepsilon_{j,2}| \le \frac{\sigma_{\max}^2}{n|\operatorname{Im} z|}.$$

Therefore,

$$\begin{aligned}
|Z_{j,4}| &= |z\hat{\rho}_j\tilde{b}_j(z)/\alpha_j||\varepsilon_{j,1} + \varepsilon_{j,2}| \\
&\le |\hat{\rho}_j|\frac{|z|^2}{|\operatorname{Im} z|^2}(|\varepsilon_{j,1}| + |\varepsilon_{j,2}|) \\
&\le \frac{\mathbf{a}_{\max}^2 |z|^2}{|\operatorname{Im} z|^3}(|\varepsilon_{j,1}| + |\varepsilon_{j,2}|),
\end{aligned}$$

where the last inequality follows from $|\hat{\rho}_j| \leq \|Q^{(j)}(z)\|_{\mathrm{sp}} \|\mathbf{a}_j\|^2 \leq \frac{\mathbf{a}_{\max}^2}{|\mathrm{Im}\,z|}$. Gathering (B.15) and (B.16), we obtain via Minkowski's inequality the desired result.

**B.3. Evaluation of $Z_5^n$.** Recall that $Z_5^n = \frac{1}{N} \sum_{j=1}^{n} \frac{\hat{\omega}_j}{\alpha_j} - \frac{1}{N} \mathrm{Tr}(B^{-1} + zI) R U Q$. We first prove that

$$(B.17) \qquad \frac{1}{N} \mathrm{Tr}(B^{-1} + zI) R U Q = -\frac{1}{n} \sum_{j=1}^{n} z \tilde{q}_{jj}(z) \frac{1}{N} \mathrm{Tr}\, D_j R U Q.$$

This follows from the definition of $b_i(z)$. We have $\frac{1}{b_i(z)} + z = -z\, \mathrm{Tr}\, \tilde{D}_i \tilde{Q}$, thus

$$
\begin{aligned}
B^{-1} + zI &= \mathrm{diag}\left(-\frac{z}{n} \mathrm{Tr}\, \tilde{D}_i \tilde{Q}; 1 \leq i \leq N\right) \\
&= -\frac{z}{n} \mathrm{diag}\left(\sum_{j=1}^{n} \tilde{q}_{jj} \sigma_{ij}^2(n); 1 \leq i \leq N\right) \\
&= -\frac{z}{n} \sum_{j=1}^{n} \tilde{q}_{jj} \, \mathrm{diag}(\sigma_{ij}^2(n); 1 \leq i \leq N) \\
&= -\frac{z}{n} \sum_{j=1}^{n} \tilde{q}_{jj} D_j,
\end{aligned}
$$

which yields (B.17). Since $\alpha_j^{-1} = -z \tilde{q}_{jj}(z)$, we have

$$
\begin{aligned}
Z_5^n &= \frac{1}{N} \sum_{j=1}^{n} \frac{\hat{\omega}_j}{\alpha_j} + \frac{1}{n} \sum_{j=1}^{n} z \tilde{q}_{jj}(z) \frac{1}{N} \mathrm{Tr}\, D_j R U Q \\
&= -\frac{1}{N} \sum_{j=1}^{n} z \tilde{q}_{jj}(z) \left\{ \hat{\omega}_j - \frac{1}{n} \mathrm{Tr}\, D_j R U Q \right\}.
\end{aligned}
$$

$(B.18)$

We now study the asymptotic behavior of $\hat{\omega}_j - \frac{1}{N} \mathrm{Tr}(D_j R U Q)$. Since

$$
\begin{aligned}
&\hat{\omega}_j - \mathbf{y}_j^T R^{(j)} U Q^{(j)} \mathbf{y}_j \\
&= \mathbf{y}_j^T (R - R^{(j)}) U Q^{(j)} \mathbf{y}_j \\
&= \mathbf{y}_j^T R (R^{(j)-1} - R^{-1}) R^{(j)} U Q^{(j)} \mathbf{y}_j \\
&= \mathbf{y}_j^T R (B^{(j)-1} - z A^{(j)} \tilde{B}^{(j)} A^{(j)T} - B^{-1} + z A \tilde{B} A^T) R^{(j)} U Q^{(j)} \mathbf{y}_j,
\end{aligned}
$$

we have

$$
\hat{\omega}_j - \frac{1}{N} \mathrm{Tr}(D_j R(z) U Q(z)) \overset{\triangle}{=} \chi_j,
$$

where

$$\chi_j = \chi_{j,1} + \chi_{j,2} + \chi_{j,3} + \chi_{j,4}$$

with

$$\chi_{j,1} = \mathbf{y}_j^T R^{(j)} U Q^{(j)} \mathbf{y}_j - \frac{1}{N} \operatorname{Tr}(D_j R(z) U Q(z)),$$

$$\chi_{j,2} = \mathbf{y}_j^T R (B^{(j)^{-1}} - B^{-1}) R^{(j)} U Q^{(j)} \mathbf{y}_j,$$

$$\chi_{j,3} = -z \sum_{\substack{\ell=1 \\ \ell \neq j}}^n (\tilde{b}_\ell^{(j)} - \tilde{b}_\ell) \mathbf{y}_j^T R \mathbf{a}_\ell \mathbf{a}_\ell^T R^{(j)} U Q^{(j)} \mathbf{y}_j,$$

$$\chi_{j,4} = z \tilde{b}_j \mathbf{y}_j^T R \mathbf{a}_j \mathbf{a}_j^T R^{(j)} U Q^{(j)} \mathbf{y}_j.$$

As usual we choose $\varepsilon > 0$ that satisfies Assumption A-1. We first handle $\chi_{j,1}$. Using the same arguments as those in Section B.2 to handle

$$\omega_j - \frac{1}{n} \operatorname{Tr} D_j Q = \frac{1}{n} \mathbf{y}_j^T Q^{(j)}(z) \mathbf{y}_j - \frac{1}{n} \operatorname{Tr} D_j Q,$$

we obtain

(B.19)                $$\mathbb{E} |\chi_{j,1}|^{2+\varepsilon/2} \leq \frac{K}{n^{1+\varepsilon/4}}.$$

The random variable $\chi_{j,2}$ satisfies

$$|\chi_{j,2}| \leq \|R\|_{\mathrm{sp}} \|B^{(j)^{-1}} - B^{-1}\|_{\mathrm{sp}} \|R^{(j)}\|_{\mathrm{sp}} \|U\|_{\mathrm{sp}} \|Q^{(j)}\|_{\mathrm{sp}} \|\mathbf{y}_j\|^2,$$

$$\leq \frac{K_1}{|\operatorname{Im} z|^3} \|B^{(j)^{-1}} - B^{-1}\|_{\mathrm{sp}} \|\mathbf{y}_j\|^2 \overset{\text{(a)}}{\leq} \frac{2\sigma_{\max}^2 |z| K_1}{n |\operatorname{Im} z|^4} \|\mathbf{y}_j\|^2,$$

where (a) follows from (B.6). As a result, $\mathbb{E} |\chi_{j,2}|^{2+\varepsilon/2} \leq \frac{K}{n^{2+\varepsilon/2}} E \|\mathbf{y}_j\|^{4+\varepsilon}$. Using (B.14) we obtain

(B.20)                $$\mathbb{E} |\chi_{j,2}|^{2+\varepsilon/2} \leq \frac{K}{n^{2+\varepsilon/2}}.$$

Consider now the random variable $\chi_{j,3}$. Using the upper bound (B.11) for $|\tilde{b}_\ell^{(j)} - \tilde{b}_\ell|$, we have

$$|\chi_{j,3}| \leq \frac{\sigma_{\max}^2 |z|^2}{n |\operatorname{Im} z|^3} \sum_{\substack{\ell=1 \\ \ell \neq j}}^n |\mathbf{y}_j^T R \mathbf{a}_\ell \mathbf{a}_\ell^T R^{(j)} U Q^{(j)} \mathbf{y}_j|.$$

Minkowski and Cauchy–Schwarz inequalities yield

(B.21)   $$\|\chi_{j,3}\|_{2+\varepsilon/2} \leq \left( \frac{\sigma_{\max}^2 |z|^2}{n |\operatorname{Im} z|^3} \right) \left[ \sum_{\substack{\ell=1 \\ \ell \neq j}}^n (\|\mathbf{y}_j^T R \mathbf{a}_\ell\|_{4+\varepsilon} \|\mathbf{a}_\ell^T R^{(j)} U Q^{(j)} \mathbf{y}_j\|_{4+\varepsilon}) \right].$$

Consider first

$$
\begin{aligned}
\text{(B.22)} \qquad \mathbb{E}|\mathbf{y}_j^T R \mathbf{a}_\ell|^{4+\varepsilon} &\le \mathbb{E}(\|R\|_{\mathrm{sp}}^{4+\varepsilon} \|\mathbf{a}_\ell\|^{4+\varepsilon} \|\mathbf{y}_j\|^{4+\varepsilon}) \\
&\le \frac{\mathbf{a}_{\max}^{4+\varepsilon}}{|\operatorname{Im} z|^{4+\varepsilon}} \mathbb{E}\|\mathbf{y}_j\|^{4+\varepsilon} \overset{\text{(a)}}{\le} K,
\end{aligned}
$$

where (a) follows from (B.14). Let us now consider the term $\mathbb{E}|\mathbf{a}_\ell^T R^{(j)} U Q^{(j)} \times \mathbf{y}_j|^{4+\varepsilon} = \frac{1}{n^{2+\varepsilon/2}} \mathbb{E}|\mathbf{v}^T \mathbf{x}_j|^{4+\varepsilon}$ with $\mathbf{v} = D_j^{1/2} Q^{(j)T} U R^{(j)T} \mathbf{a}_\ell$. By a series of inequalities similar to B.2, we obtain $\mathbb{E}|\mathbf{a}_\ell^T R^{(j)} U Q^{(j)} \mathbf{y}_j|^{4+\varepsilon} \le K/n^{2+\varepsilon/2}$. In conclusion, we obtain

$$
\text{(B.23)} \qquad \mathbb{E}|\chi_{j,3}|^{2+\varepsilon/2} \le \frac{K}{n^{1+\varepsilon/4}}.
$$

Finally, the variable $\chi_{j,4}$ can be handled in the same way

$$
\begin{aligned}
\text{(B.24)} \qquad &\mathbb{E}|\chi_{j,4}|^{2+\varepsilon/2} \\
&\le \left|\frac{z}{\operatorname{Im} z}\right|^{2+\varepsilon/2} (\mathbb{E}|\mathbf{y}_j^T R \mathbf{a}_j|^{4+\varepsilon})^{1/2} (\mathbb{E}|\mathbf{a}_j^T R^{(j)} U Q^{(j)} \mathbf{y}_j|^{4+\varepsilon})^{1/2}, \\
&\le \frac{K}{n^{1+\varepsilon/4}}.
\end{aligned}
$$

In order to finish the proof, it remains to gather equations (B.19), (B.20), (B.23) and (B.24) to get

$$
\mathbb{E}|\chi_{j,1} + \chi_{j,2} + \chi_{j,3} + \chi_{j,4}|^{2+\varepsilon/2} \le \frac{K}{n^{1+\varepsilon/4}}.
$$

Plugging this into (B.18) and using Minkowski's inequality, we obtain

$$
\mathbb{E}\left|\frac{1}{N}\sum_{j=1}^n z\tilde{q}_{jj}\chi_j\right|^{2+\varepsilon/2} \le \frac{K}{n^{1+\varepsilon/4}},
$$

which is the desired result.

## APPENDIX C: PROOF OF THEOREM 4.1

**C.1. Proof of the convergence of $C_n(\sigma^2) - \overline{C}_n(\sigma^2)$ to zero.** The proof of the convergence relies on a dominated convergence argument. Denote by $\pi_n$ the probability measure whose Stieltjes transform is $\frac{1}{N}\operatorname{Tr} T_n(z)$; similarly, denote by $\mathbb{P}_n$ the probability measure whose Stieltjes transform is $\frac{1}{N}\operatorname{Tr} Q_n(z)$:

$$
\frac{1}{N}\operatorname{Tr} T_n(z) = \int_0^\infty \frac{\pi_n(d\lambda)}{\lambda - z} \quad \text{and} \quad \frac{1}{N}\operatorname{Tr} Q_n(z) = \int_0^\infty \frac{\mathbb{P}_n(d\lambda)}{\lambda - z}.
$$

The following estimates will be useful.

LEMMA C.1. *The following equalities hold true*

$$\text{(C.1)} \qquad \mathbb{E}\int_0^\infty \lambda \mathbb{P}_n(d\lambda) = \frac{1}{Nn}\sum_{\substack{1\le i\le N\\1\le j\le n}}\sigma_{ij}^2 + \frac{1}{N}\operatorname{Tr} AA^T,$$

$$\text{(C.2)} \qquad \int_0^\infty \lambda \pi_n(d\lambda) = \frac{1}{Nn}\sum_{\substack{1\le i\le N\\1\le j\le n}}\sigma_{ij}^2 + \frac{1}{N}\operatorname{Tr} AA^T.$$

*In particular,*

$$\text{(C.3)} \qquad \sup_n\left(\mathbb{E}\int_0^\infty \lambda \mathbb{P}_n(d\lambda)\right) = \sup_n\left(\int_0^\infty \lambda \pi_n(d\lambda)\right) \le \sigma_{\max}^2 + \mathbf{a}_{\max}^2.$$

PROOF. First notice that

$$\mathbb{E}\int_0^\infty \lambda \mathbb{P}_n(d\lambda) = \frac{1}{N}\mathbb{E}\operatorname{Tr}\Sigma\Sigma^T,$$

which yields immediately equality (C.1). We now compute $\int_0^\infty \lambda \pi_n(d\lambda)$. We first prove that

$$\text{(C.4)} \qquad \int_0^\infty \lambda \pi_n(d\lambda) = \lim_{y\to\infty}\operatorname{Re}\left[-\mathbf{i}y\left(\mathbf{i}y\frac{1}{N}\operatorname{Tr} T_n(\mathbf{i}y)+1\right)\right].$$

Compute

$$-\mathbf{i}y\left(\mathbf{i}y\frac{1}{N}\operatorname{Tr} T_n(\mathbf{i}y)+1\right) = -\mathbf{i}y\left(\int_0^\infty\frac{\mathbf{i}y\pi_n(d\lambda)}{\lambda-\mathbf{i}y}+1\right)$$

$$= -\mathbf{i}y\int_0^\infty\frac{\lambda^2\pi_n(d\lambda)}{\lambda^2+y^2}+y^2\int_0^\infty\frac{\lambda\pi_n(d\lambda)}{\lambda^2+y^2}.$$

Hence,

$$\operatorname{Re}\left[-\mathbf{i}y\left(\mathbf{i}y\frac{1}{N}\operatorname{Tr} T_n(\mathbf{i}y)+1\right)\right] = y^2\int_0^\infty\frac{\lambda\pi_n(d\lambda)}{\lambda^2+y^2}.$$

The monotone convergence theorem yields $\lim_{y\to\infty}y^2\int_0^\infty\frac{\lambda\pi_n(d\lambda)}{\lambda^2+y^2}=\int_0^\infty \lambda\pi_n(d\lambda)$. Equation (C.4) is proved. We now prove

$$\text{(C.5)} \qquad \lim_{y\to\infty}-\mathbf{i}y(\mathbf{i}yT_n(\mathbf{i}y)+I_N) = \frac{1}{n}\operatorname{diag}(\operatorname{Tr}\tilde{D}_i, 1\le i\le N)+AA^T.$$

The mere definition of $T_n$ yields

$$T_n(z) = (\Psi^{-1}-zA\tilde{\Psi}A^T)^{-1} = (I_N-z\Psi A\tilde{\Psi}A^T)^{-1}\Psi.$$

Using the fact that

$$(I_N-z\Psi A\tilde{\Psi}A^T)^{-1} = I_N+z\Psi A\tilde{\Psi}A^T+(z\Psi A\tilde{\Psi}A^T)^2(I_N-z\Psi A\tilde{\Psi}A^T)^{-1},$$

we get that

$$T_n(z) = \Psi + z\Psi A\tilde{\Psi}A^T\Psi + (z\Psi A\tilde{\Psi}A^T)^2 T_n(z).$$

We now compute

(C.6)
$$\begin{aligned}-z&(zT_n(z)+I_N)\\ &= -z(z\Psi(z)+I_N) - z\Psi Az\tilde{\Psi}A^T z\Psi - z\Psi Az\tilde{\Psi}A^T z\Psi Az\tilde{\Psi}A^T T_n(z).\end{aligned}$$

In order to compute the limit of the previous expression when $z = \mathbf{i}y$, $y \to \infty$, recall that $\lim -\mathbf{i}y\Psi(\mathbf{i}y) = \lim -\mathbf{i}yT_n(\mathbf{i}y) = I_N$ and $\lim -\mathbf{i}y\tilde{\Psi}(\mathbf{i}y) = \lim -\mathbf{i}y \times \tilde{T}_n(\mathbf{i}y) = I_n$ whenever $y \to \infty$ by (5.2) and (5.4).

Let us first consider the first term on the right-hand side of (C.6).

$$\begin{aligned}-\mathbf{i}y(\mathbf{i}y\Psi(\mathbf{i}y)+I_N) &= \operatorname{diag}\left(\frac{\mathbf{i}y}{1+(1/n)\operatorname{Tr}\tilde{D}_i\tilde{T}_n(\mathbf{i}y)} - \mathbf{i}y\right)\\ &= \operatorname{diag}\left(\frac{-(\mathbf{i}y/n)\operatorname{Tr}\tilde{D}_i\tilde{T}_n(\mathbf{i}y)}{1+(1/n)\operatorname{Tr}\tilde{D}_i\tilde{T}_n(\mathbf{i}y)}\right)\\ &\xrightarrow[y\to\infty]{} \frac{1}{n}\operatorname{diag}(\operatorname{Tr}\tilde{D}_i, 1\le i\le N).\end{aligned}$$

For the second term on the right-hand side of (C.6) we have

$$(-\mathbf{i}y\Psi(\mathbf{i}y))A(-\mathbf{i}y\tilde{\Psi}(\mathbf{i}y))A^T(-\mathbf{i}y\Psi(\mathbf{i}y))\xrightarrow[y\to\infty]{} AA^T.$$

The third term clearly converges to 0 because $T_n(\mathbf{i}y) \to 0$ when $y \to +\infty$. Equation (C.5) is established. This limit immediately yields

$$\lim_{y\to\infty} -\mathbf{i}y\left(\mathbf{i}y\frac{1}{N}\operatorname{Tr} T_n(\mathbf{i}y)+1\right) = \frac{1}{Nn}\sum_{\substack{1\le i\le N\\ 1\le j\le n}}\sigma_{ij}^2 + \frac{1}{N}\operatorname{Tr} AA^T.$$

Equating this equation with (C.4) implies that $\int\lambda\pi_n(d\lambda)$ is finite and gives its value. Therefore, equation (C.2) is proved. The inequality (C.3) follows immediately from (C.1) and (C.2). Proof of Lemma C.1 is complete.  $\square$

We are now in position to prove that $C_n(\sigma^2) - \overline{C}_n(\sigma^2) \to 0$. Recall that

$$C_n(\sigma^2) = \int_{\sigma^2}^\infty\left(\frac{1}{\omega} - \frac{1}{N}\mathbb{E}\operatorname{Tr} Q_n(-\omega)\right)d\omega.$$

The dominated convergence theorem together with Theorem 2.5 yield

$$\forall\omega > 0 \qquad \frac{1}{\omega} - \frac{1}{N}\operatorname{Tr} T_n(-\omega) - \left(\frac{1}{\omega} - \frac{1}{N}\mathbb{E}\operatorname{Tr} Q_n(-\omega)\right)\xrightarrow[n\to\infty]{} 0.$$

Moreover,

$$\left| \frac{1}{\omega} - \frac{1}{N}\operatorname{Tr} T_n(-\omega) - \left( \frac{1}{\omega} - \frac{1}{N}\mathbb{E}\operatorname{Tr} Q_n(-\omega) \right) \right|$$

$$\leq \left| \frac{1}{\omega} - \frac{1}{N}\operatorname{Tr} T_n(-\omega) \right| + \left| \frac{1}{\omega} - \frac{1}{N}\mathbb{E}\operatorname{Tr} Q_n(-\omega) \right|$$

$$= \left| \int_0^\infty \left( \frac{1}{\omega} - \frac{1}{\lambda + \omega} \right)\pi_n(d\lambda) \right| + \left| \mathbb{E}\int_0^\infty \left( \frac{1}{\omega} - \frac{1}{\lambda + \omega} \right)\mathbb{P}_n(d\lambda) \right|$$

$$\leq \left| \frac{\int_0^\infty \lambda\pi_n(d\lambda)}{\omega^2} \right| + \left| \frac{\mathbb{E}\int_0^\infty \lambda\mathbb{P}_n(d\lambda)}{\omega^2} \right|$$

$$\leq 2\frac{a_{\max}^2 + \sigma_{\max}^2}{\omega^2},$$

which is integrable in $\omega$ over $(\sigma^2, \infty)$. Therefore the dominated convergence theorem yields $C_n(\sigma^2) - \overline{C}_n(\sigma^2) \to 0$ and the first part of Theorem 4.1 is proved.

**C.2. Proof of formula (4.3): some preparation.** Performing the change of variable $\gamma = \frac{1}{\omega}$ in (4.2) yields to the formula

$$\overline{C}(\sigma^2) = \int_0^{1/\sigma^2} \frac{1}{\gamma}\left( 1 - \frac{1}{N}\operatorname{Tr}\frac{1}{\gamma}T\left( -\frac{1}{\gamma} \right) \right)d\gamma.$$

One can check that the integrand is continuous in zero. In fact, Lemma C.1 yields

$$\frac{1}{\gamma}\left( 1 - \frac{1}{N}\operatorname{Tr}\frac{1}{\gamma}T\left( -\frac{1}{\gamma} \right) \right) = \int \frac{1}{\gamma}\left( 1 - \frac{1}{\gamma}\frac{1}{(\lambda + \gamma^{-1})} \right)\pi_n(d\lambda)$$

$$= \int \frac{\lambda\pi_n(d\lambda)}{1 + \lambda\gamma}\xrightarrow[\gamma\to 0]{} \int \lambda\pi_n(d\lambda) < \infty,$$

where $\pi_n$ is the probability measure associated to the Stieltjes transform $\frac{1}{N}\operatorname{Tr} T(z)$ We thus introduce slightly different notation than in the rest of the paper. These notations appear to be more convenient in the forthcoming computations. We denote by

$$S(\gamma) = \frac{1}{\gamma}T\left( -\frac{1}{\gamma} \right), \qquad \theta_i(\gamma) = \frac{1}{\gamma}\psi_i\left( -\frac{1}{\gamma} \right), \qquad \Theta = \operatorname{diag}(\theta_i, 1 \leq i \leq N),$$

$$\tilde{S}(\gamma) = \frac{1}{\gamma}\tilde{T}\left( -\frac{1}{\gamma} \right), \qquad \tilde{\theta}_j(\gamma) = \frac{1}{\gamma}\tilde{\psi}_j\left( -\frac{1}{\gamma} \right), \qquad \tilde{\Theta} = \operatorname{diag}(\tilde{\theta}_j, 1 \leq j \leq n).$$

These notations yield

$$S(\gamma) = (\Theta^{-1}(\gamma) + \gamma A\tilde{\Theta}(\gamma)A^T)^{-1}, \qquad \tilde{S}(\gamma) = (\tilde{\Theta}^{-1}(\gamma) + \gamma A^T\Theta(\gamma)A)^{-1}.$$

$$\theta_i(\gamma) = \frac{1}{1 + (\gamma/n)\operatorname{Tr}\tilde{D}_i\tilde{S}(\gamma)}, \qquad \tilde{\theta}_j(\gamma) = \frac{1}{1 + (\gamma/n)\operatorname{Tr}D_jS(\gamma)}.$$

The general strategy to establish formula (4.3) is to write $\frac{1}{\gamma}(1 - \frac{1}{N}\operatorname{Tr}S(\gamma))$ as the derivative of some well-identified function of $\gamma$. The following quantities will also be of help:

$$\Delta = \frac{1}{n}\operatorname{diag}(\operatorname{Tr}(D_jS, 1 \leq j \leq n)), \qquad \tilde{\Delta} = \frac{1}{n}\operatorname{diag}(\operatorname{Tr}(\tilde{D}_i\tilde{S}, 1 \leq i \leq N)).$$

In the sequel, we use both $f'(\gamma)$ and $\frac{df}{d\gamma}(\gamma)$ for the derivative of $f$.

LEMMA C.2. *The following equality holds true:*

(C.7) $$SA\tilde{\Theta} = \Theta A\tilde{S}.$$

*In particular,*

(C.8) $$\operatorname{Tr}A\tilde{\Theta}A^TS = \operatorname{Tr}A^T\Theta A\tilde{S},$$

(C.9) $$\operatorname{Tr}A\tilde{\Theta}'A^TS = \operatorname{Tr}\tilde{\Theta}'A^T\Theta A\tilde{S}\tilde{\Theta}^{-1}.$$

PROOF. After elementary matrix manipulations (see also [15], Section 0.7.4), we have

$$(I + UV)^{-1} = I - U(I + VU)^{-1}V,$$

which yields afterward

$$(I + UV)^{-1}U = U(I + VU)^{-1}.$$

Let $U = \Theta^{1/2}A$ and $V = \gamma\tilde{\Theta}A^T\Theta^{1/2}$. Then

$$(I + \gamma\Theta^{1/2}A\tilde{\Theta}A^T\Theta^{1/2})^{-1}\Theta^{1/2}A = \Theta^{1/2}A(I + \gamma\tilde{\Theta}A^T\Theta A)^{-1}$$

$$\iff \quad \Theta^{-1/2}S\Theta^{-1/2}\Theta^{1/2}A = \Theta^{1/2}A\tilde{S}\tilde{\Theta}^{-1}$$

$$\iff \quad SA\tilde{\Theta} = \Theta A\tilde{S},$$

which is the desired result. Multiplying by $A^T$, we obtain $A^TSA\tilde{\Theta} = A^T\Theta A\tilde{S}$ which yields (C.8). Multiplying to the left by $\tilde{\Theta}'A^T$ and to the right by $\tilde{\Theta}^{-1}$, we obtain $\tilde{\Theta}'A^TSA = \tilde{\Theta}'A^T\Theta A\tilde{S}\tilde{\Theta}^{-1}$ which yields (C.9). □

LEMMA C.3. *Denote by*

$$I(\gamma) = -\frac{\gamma}{N}\operatorname{Tr}A^T\Theta'A\tilde{S} + \frac{\gamma}{N}\operatorname{Tr}A\tilde{\Theta}A^TS'.$$

*Then the following equality holds true:*

(C.10) $$I(\gamma) = \frac{d}{d\gamma}\left(\frac{1}{N}\operatorname{Tr}(\gamma A\tilde{\Theta}A^TS)\right) - \frac{d}{d\gamma}\left(\frac{1}{N}\log\det\tilde{\Theta}(\tilde{\Theta}^{-1} + \gamma A^T\Theta A)\right).$$

In order to prove Lemma C.3, we shall rely on a differentiation formula. Let $A = A(x)$ be an $N \times N$ matrix, then the equality

(C.11) $$\frac{d}{dx} \log \det A(x) = \operatorname{Tr} A^{-1}(x) \frac{d}{dx} A(x)$$

holds true in the case where $A$ is Hermitian or $A$ is the product of two Hermitian matrices, that is, $A(x) = B(x)C(x)$. We provide here a short proof (in the case of a general square matrix, one can refer to [14], Section 15).

Consider first the case where $A$ is Hermitian and write $A = U\Delta U^*$ where $\Delta(x) = \operatorname{diag}(\lambda_i(x); 1 \le i \le N)$ and $UU^* = U^*U = I$. We have

$$\frac{d}{dx} \log \det A(x) = \frac{d}{dx} \left( \sum_{i=1}^{N} \log \lambda_i(x) \right) = \sum_{i=1}^{N} \frac{\lambda_i'(x)}{\lambda_i(x)}.$$

On the other hand,

$$\operatorname{Tr} A^{-1} \frac{d}{dx} A = \operatorname{Tr} U \Delta^{-1} U^* (U' \Delta U^* + U \Delta' U^* + U \Delta U^{*'})$$

$$= \operatorname{Tr} \Delta^{-1} \Delta' + \operatorname{Tr} U \Delta^{-1} U^* U' \Delta U^* + \operatorname{Tr} U \Delta^{-1} U^* U \Delta U^{*'}$$

$$= \sum_{i=1}^{N} \frac{\lambda_i'(x)}{\lambda_i(x)} + \operatorname{Tr}(U^* U' + U^{*'} U) \stackrel{\text{(a)}}{=} \sum_{i=1}^{N} \frac{\lambda_i'(x)}{\lambda_i(x)},$$

where (a) follows from the fact that $U^* U' + U^{*'} U = 0$ which is obtained by differentiating $U^* U = I$.

In the case where $A(x) = B(x)C(x)$ with $B$ and $C$ Hermitian, we have

$$\frac{d}{dx} \log \det B(x) C(x) = \frac{d}{dx} \log \det B(x) + \frac{d}{dx} \log \det C(x)$$

$$= \operatorname{Tr} B^{-1} B' + \operatorname{Tr} C^{-1} C'$$

$$= \operatorname{Tr} B^{-1} B' C C^{-1} + \operatorname{Tr} C^{-1} B^{-1} B C'$$

$$= \operatorname{Tr} C^{-1} B^{-1} (B' C + B C') = \operatorname{Tr} A^{-1} A',$$

which is the expected result. We are now in position to prove Lemma C.3.

PROOF OF LEMMA C.3. We differentiate the expression

$$\frac{d}{d\gamma} \left( \frac{1}{N} \operatorname{Tr}(\gamma A \tilde{\Theta} A^T S) \right) - \frac{d}{d\gamma} \left( \frac{1}{N} \log \det \tilde{\Theta} (\tilde{\Theta}^{-1} + \gamma A^T \Theta A) \right)$$

$$= \frac{1}{N} \operatorname{Tr}(A \tilde{\Theta} A^T S) + \frac{\gamma}{N} \operatorname{Tr}(A \tilde{\Theta}' A^T S) + \frac{\gamma}{N} \operatorname{Tr}(A \tilde{\Theta} A^T S')$$

$$- \frac{1}{N} \operatorname{Tr}[(\tilde{\Theta} A^T \Theta A)(I_n + \gamma \tilde{\Theta} A^T \Theta A)^{-1}]$$

$$- \frac{\gamma}{N} \operatorname{Tr}[(\tilde{\Theta}' A^T \Theta A)(I_n + \gamma \tilde{\Theta} A^T \Theta A)^{-1}]$$

$$- \frac{\gamma}{N} \operatorname{Tr}[(\tilde{\Theta} A^T \Theta' A)(I_n + \gamma \tilde{\Theta} A^T \Theta A)^{-1}]$$

$$\overset{(a)}{=} I(\gamma) + \frac{1}{N} \operatorname{Tr}(A \tilde{\Theta} A^T S) - \frac{1}{N} \operatorname{Tr}(A^T \Theta A \tilde{S})$$

$$+ \frac{\gamma}{N} \operatorname{Tr}(A \tilde{\Theta}' A^T S) - \frac{\gamma}{N} \operatorname{Tr}[(\tilde{\Theta}' A^T \Theta A)(I_n + \gamma \tilde{\Theta} A^T \Theta A)^{-1}]$$

$$\overset{(b)}{=} I(\gamma) + \frac{\gamma}{N} \operatorname{Tr}(A \tilde{\Theta}' A^T S) - \frac{\gamma}{N} \operatorname{Tr}[(\tilde{\Theta}' A^T \Theta A)(I_n + \gamma \tilde{\Theta} A^T \Theta A)^{-1}]$$

$$\overset{(c)}{=} I(\gamma),$$

where (a) follows from the fact that $(I_n + \gamma \tilde{\Theta} A^T \Theta A)^{-1} = \tilde{S} \tilde{\Theta}^{-1}$, (b) follows from (C.8) and (c) follows from (C.9). $\square$

**C.3. Proof of formula (4.3): the main system of equations.** We are now in position to prove the following:

LEMMA C.4. *Denote by*

$$J(\gamma) = \frac{1}{\gamma} \left( 1 - \frac{1}{N} \operatorname{Tr} S(\gamma) \right),$$

$$J_1(\gamma) = \frac{\gamma}{N} \operatorname{Tr} \tilde{\Delta}' S + \frac{\gamma}{N} \operatorname{Tr} A \tilde{\Theta}' A^T S,$$

$$J_2(\gamma) = \frac{\gamma}{N} \operatorname{Tr} \Delta' \tilde{S} + \frac{\gamma}{N} \operatorname{Tr} A^T \Theta' A \tilde{S}.$$

*Then the following system of equations holds true:*

$$(\text{C.12}) \quad \frac{d}{d\gamma} \left( \frac{1}{N} \log \det(\Theta^{-1} + \gamma A \tilde{\Theta} A^T) \right) = J(\gamma) + J_1(\gamma),$$

$$(\text{C.13}) \quad \frac{d}{d\gamma} \left( \frac{1}{N} \log \det(\tilde{\Theta}^{-1} + \gamma A^T \Theta A) \right) = J(\gamma) + J_2(\gamma),$$

$$(\text{C.14}) \quad \frac{d}{d\gamma} \left( \frac{\gamma}{N} \operatorname{Tr}(\tilde{\Delta} + A \tilde{\Theta} A^T) S \right) = J(\gamma) + J_1(\gamma) \\ + \frac{\gamma}{N} \operatorname{Tr}(\tilde{\Delta} + A \tilde{\Theta} A^T) S'.$$

PROOF. We first give equivalent formulations for the term $J(\gamma)$:

$$J(\gamma) = \frac{1}{\gamma} \left( 1 - \frac{1}{N} \operatorname{Tr} S(\gamma) \right)$$

$$= \frac{1}{\gamma} \left( \frac{1}{N} \operatorname{Tr} S^{-1} S - \frac{1}{N} \operatorname{Tr}(\Theta^{-1} + \gamma A \tilde{\Theta} A^T)^{-1} \right)$$

$$= \frac{1}{\gamma N} \operatorname{Tr}((\Theta^{-1} - I_N)S + \gamma A \tilde{\Theta} A^T S)$$

$$(\text{C.15}) \qquad = \frac{1}{N} \operatorname{Tr} \tilde{\Delta} S + \frac{1}{N} \operatorname{Tr} A \tilde{\Theta} A^T S,$$

$$(\text{C.16}) \qquad \overset{(\text{a})}{=} \frac{1}{N} \operatorname{Tr} \Delta \tilde{S} + \frac{1}{N} \operatorname{Tr} A^T \Theta A \tilde{S},$$

where (a) follows from (C.8) and from the fact that $\operatorname{Tr} \tilde{\Delta} S = \operatorname{Tr} \Delta \tilde{S}$. Consider now

$$\frac{d}{d\gamma} \left( \frac{1}{N} \log \det(\Theta^{-1} + \gamma A \tilde{\Theta} A^T) \right) = \frac{1}{N} \operatorname{Tr}(\Theta^{-1'} + \gamma A \tilde{\Theta}' A^T + A \tilde{\Theta} A^T) S.$$

Easy computation yields $\Theta^{-1'} = \tilde{\Delta} + \gamma \tilde{\Delta}'$ and the previous equality becomes

$$\frac{d}{d\gamma} \left( \frac{1}{N} \log \det(\Theta^{-1} + \gamma A \tilde{\Theta} A^T) \right) = J(\gamma) + J_1(\gamma),$$

where (C.15) has been used to identify $J$. Equation (C.12) is proved. One can prove similarly (C.13) by using (C.16). We now compute

$$\frac{d}{d\gamma} \left( \frac{\gamma}{N} \operatorname{Tr}(\tilde{\Delta} + A \tilde{\Theta} A^T) S \right)$$

$$= \frac{\gamma}{N} \operatorname{Tr} \tilde{\Delta}' S + \frac{\gamma}{N} \operatorname{Tr} A \tilde{\Theta}' A^T S + \frac{1}{N} \operatorname{Tr} \tilde{\Delta} S + \frac{1}{N} \operatorname{Tr} A \tilde{\Theta} A^T S$$

$$+ \frac{\gamma}{N} \operatorname{Tr}(\tilde{\Delta} + A \tilde{\Theta} A^T) S'$$

$$\overset{(\text{a})}{=} J(\gamma) + J_1(\gamma) + + \frac{\gamma}{N} \operatorname{Tr}(\tilde{\Delta} + A \tilde{\Theta} A^T) S',$$

where (a) follows from (C.15). Equation (C.14) is proved. □

## C.4. Proof of formula (4.3): end of the proof.

Eliminating $J(\gamma)$ between (C.13) and (C.14), we end up with

$$\frac{d}{d\gamma} \left( \frac{\gamma}{N} \operatorname{Tr}(\tilde{\Delta} + A \tilde{\Theta} A^T) S \right)$$

$$= J_1(\gamma) + \frac{d}{d\gamma} \left( \frac{1}{N} \log \det(\tilde{\Theta}^{-1} + \gamma A^T \Theta A) \right)$$

$$- \frac{\gamma}{N} \operatorname{Tr} \Delta' \tilde{S} - \frac{\gamma}{N} \operatorname{Tr} A^T \Theta' A \tilde{S} + \frac{\gamma}{N} \operatorname{Tr}(\tilde{\Delta} + A \tilde{\Theta} A^T) S'.$$

Since $-\frac{\gamma}{N} \operatorname{Tr} \Delta' \tilde{S} + \frac{\gamma}{N} \operatorname{Tr} \tilde{\Delta} S' = 0$, we obtain

$$(\text{C.17}) \qquad \frac{d}{d\gamma} \left( \frac{\gamma}{N} \operatorname{Tr}(\tilde{\Delta} + A \tilde{\Theta} A^T) S \right)$$

$$= J_1(\gamma) + \left( \frac{1}{N} \log \det(\tilde{\Theta}^{-1} + \gamma A^T \Theta A) \right)' + I(\gamma).$$

First use (C.10) which expresses $I(\gamma)$ as a derivative, then extract $J_1(\gamma)$ from (C.17) and plug it into (C.12). After some simplifications, we obtain

$$J(\gamma) = \left(\frac{1}{N}\log\det(\Theta^{-1} + \gamma A\tilde{\Theta}A^T)\right)' - \left(\frac{\gamma}{N}\operatorname{Tr}\tilde{\Delta}S\right)' + \left(\frac{1}{N}\log\det\tilde{\Theta}^{-1}\right)',$$

$$\stackrel{\triangle}{=} F'(\gamma).$$

Since $\overline{C}(\sigma^2) = \int_0^{1/\sigma^2} J(\gamma)\,d\gamma$, a mere integration yields $\overline{C}(\sigma^2) = F(\sigma^{-2}) - \lim_{\gamma\to 0} F(\gamma)$. It remains to check that $\lim_{\gamma\to 0} F(\gamma) = 0$ to obtain (4.3). In order to compute the limit of $F$ as $\gamma$ goes to zero, it is sufficient to check the following limits:

$$\lim_{\gamma\to 0} S(\gamma) = I_N, \qquad \lim_{\gamma\to 0} \Theta(\gamma) = I_N \quad \text{and} \quad \lim_{\gamma\to 0} \tilde{\Theta}(\gamma) = I_n.$$

Let us prove the first limit,

$$S(\gamma) = \frac{1}{\gamma}T\left(-\frac{1}{\gamma}\right) \stackrel{(a)}{=} \frac{1}{\gamma}\int \frac{\mu(d\lambda)}{\lambda + \gamma^{-1}} \xrightarrow[\gamma\to 0]{} \mu(\mathbb{R}) = I_N,$$

where $(a)$ follows from Proposition 5.1 part 2. In order to compute the limit involving $\Theta$ and $\tilde{\Theta}$, it is sufficient to note that $\theta_i(\gamma) = \frac{1}{\gamma}\psi_i(-\gamma^{-1})$, to interpret $\psi_i$ as a Stieltjes transform (cf. Proposition 5.1 part 4) and to perform the same computation as for $\frac{1}{\gamma}T(-\gamma^{-1})$. One can compute similarly the limit of $\tilde{\Theta}$. Theorem 4.1 is proved.

**Acknowledgment.** We thank the referee whose very careful reading of the manuscript has helped us to substantially improve the writing of the paper.

## REFERENCES

[1] BAI, Z. D. and SILVERSTEIN, J. W. (1998). No eigenvalues outside the support of the limiting spectral distribution of large-dimensional sample covariance matrices. *Ann. Probab.* **26** 316–345. MR1617051

[2] BAI, Z. D. and SILVERSTEIN, J. W. (2004). CLT for linear spectral statistics of large-dimensional sample covariance matrices. *Ann. Probab.* **32** 553–605. MR2040792

[3] BOLOTNIKOV, V. (1997). On a general moment problem on the half axis. *Linear Algebra and Its Applications* **255** 57–112. MR1433234

[4] BOUTET DE MONVEL, A., KHORUNZHY, A. and VASILCHUK, V. (1996). Limiting eigenvalue distribution of random matrices with correlated entries. *Markov Process. Related Fields* **2** 607–636. MR1431189

[5] DOZIER, R. B. and SILVERSTEIN, J. W. (2007). On the empirical distribution of eigenvalues of large dimensional information-plus-noise type matrices. *J. Multivariate Anal.* To appear.

[6] DUMONT, J., HACHEM, W., LOUBATON, P. and NAJIM, J. (2006). On the asymptotic analysis of mutual information of mimo Rician correlated channels. *Proceedings of the ISCCSP Conference.* Marrakech, Morocco.

[7] DUMONT, J., LOUBATON, P., LASAULCE, S. and DEBBAH, M. (2005). On the asymptotic performance of mimo correlated Ricean channels. In *ICASSP Proceedings* **5** 813–816.

[8] GESZTESY, F. and TSEKANOVSKII, E. (2000). On matrix-valued Herglotz functions. *Math. Nachr.* **218** 61–138. MR1784638

[9] GIRKO, V. L. (1990). *Theory of Random Determinants*. Kluwer, Dordrecht. MR1080966

[10] GIRKO, V. L. (2001). *Theory of Stochastic Canonical Equations*. Kluwer, Dordrecht. MR1887675

[11] GUIONNET, A. and ZEITOUNI, O. (2000). Concentration of the spectral measure for large matrices. *Electron. Comm. Probab.* **5** 119–136. MR1781846

[12] HACHEM, W., LOUBATON, P. and NAJIM, J. (2005). The empirical eigenvalue distribution of a gram matrix: From independence to stationarity. *Markov Process. Related Fields* **11** 629–648. MR2191967

[13] HACHEM, W., LOUBATON, P. and NAJIM, J. (2006). The empirical distribution of the eigenvalues of a Gram matrix with a given variance profile. *Ann. Inst. H. Poincaré Probab. Statist.* **42** 649–670.

[14] HARVILLE, D. A. (1997). *Matrix Algebra from a Statistician's Perspective*. Springer, New York. MR1467237

[15] HORN, R. and JOHNSON, C. (1985). *Matrix Analysis*. Cambridge Univ. Press. MR0832183

[16] KAILATH, T., SAYED, A. H. and HASSIBI, B. (2000). *Linear Estimation*. Prentice Hall, Englewood Cliffs, NJ.

[17] KHORUNZHY, A., KHORUZHENKO, B. and PASTUR, L. (1996). Asymptotic properties of large random matrices with independent entries. *J. Math. Phys.* **37** 5033–5060. MR1411619

[18] KREIN, M. and NUDELMAN, A. (1997). *The Markov Moment Problem and Extremal Problems*. Amer. Math. Soc., Providence, RI. MR0458081

[19] MARČENKO, V. A. and PASTUR, L. A. (1967). Distribution of eigenvalues in certain sets of random matrices. *Mat. Sb. (N.S.)* **72** 507–536. MR0208649

[20] ROZANOV, Y. A. (1967). *Stationary Random Processes*. Holden-Day, San Francisco, CA. MR0214134

[21] RUDIN, W. (1987). *Real and Complex Analysis*. McGraw-Hill, New York. MR0924157

[22] SENGUPTA, A. M. and MITRA, P. (2000). Capacity of multivariate channels with multiplicative noise: I. Random matrix techniques and large-$n$ expansions for full transfer matrices. Available at http://arxiv.org/abs/physics/0010081.

[23] SHLYAKHTENKO, D. (1996). Random Gaussian band matrices and freeness with amalgamation. *Internat. Math. Res. Notices* **20** 1013–1025. MR1422374

[24] SILVERSTEIN, J. W. (1995). Strong convergence of the empirical distribution of eigenvalues of large-dimensional random matrices. *J. Multivariate Anal.* **55** 331–339. MR1370408

[25] SILVERSTEIN, J. W. and BAI, Z. D. (1995). On the empirical distribution of eigenvalues of a class of large-dimensional random matrices. *J. Multivariate Anal.* **54** 175–192. MR1345534

[26] TULINO, A. and VERDÚ, S. (2004). Random matrix theory and wireless communications. *Foundations and Trends in Communications and Information Theory* **1** 1–182.

[27] VERDU, S. and SHAMAI, S. (1999). Spectral efficiency of CDMA with random spreading. *IEEE Trans. Inform. Theory* **45** 622–640. MR1677022

[28] YIN, Y. Q. (1986). Limiting spectral distribution for a class of random matrices. *J. Multivariate Anal.* **20** 50–68. MR0862241

W. HACHEM
SUPÉLEC (ECOLE SUPÉRIEURE D'ELECTRICITÉ)
PLATEAU DE MOULON
3 RUE JOLIOT-CURIE
91192 GIF SUR YVETTE CEDEX
FRANCE
E-MAIL: walid.hachem@supelec.fr

P. LOUBATON
IGM LABINFO
UMR 8049
INSTITUT GASPARD MONGE
UNIVERSITÉ DE MARNE LA VALLÉE
5 BD DESCARTES
CHAMPS SUR MARNE
77454 MARNE LA VALLÉE CEDEX 2
FRANCE
E-MAIL: loubaton@univ-mlv.fr

J. NAJIM
CNRS, TÉLÉCOM PARIS
46 RUE BARRAULT
75013 PARIS
FRANCE
E-MAIL: najim@enst.fr